\documentclass[11pt,letterpaper,reqno]{amsart}
\usepackage[dvipsnames]{xcolor}
\usepackage[unicode]{hyperref}
\usepackage{mathabx} 
\hypersetup{
	linktoc=all,     
	colorlinks = true,
	citecolor = red,
	anchorcolor = green,
	linkcolor = blue}
\numberwithin{equation}{section}
\usepackage{amsmath,amssymb,bm}

\newtheorem{theorem}{Theorem}[section]
\newtheorem{lemma}[theorem]{Lemma}
\newtheorem{proposition}[theorem]{Proposition}
\newtheorem{corollary}[theorem]{Corollary}

\newcommand{\Id}{\mathrm{Id}}
\newcommand{\R}{\mathbb R}

\theoremstyle{definition}

\usepackage{indentfirst}
\usepackage{fancyhdr}
\usepackage{color}
\usepackage{mathrsfs}
\usepackage{fancyhdr}
\theoremstyle{remark}
\newtheorem{remark}[theorem]{Remark}
\newcommand{\pa}{\partial}
\allowdisplaybreaks[3]

\begin{document}

\title[GLW in the critical Besov space of SMCF in $\mathbb{R}^d: d\ge 5$]{Global well-posedness in the critical Besov space of the skew mean curvature flow in $\mathbb{R}^d: d\ge 5$}

\author{Ning-An Lai}
\address{Department of Mathematics, Shaoxing University, Shaoxing 312000, China}
\email{ninganlai@usx.edu.cn(N.-A.Lai)}

\author{Jie Shao}
\address{School of Mathematics and Statistics, Hunan Normal University, Changsha
	410081, China}
\email{shaojiehn@foxmail.com(J. Shao)}

\author{Zexian Zhang}
\address{School of Mathematical Sciences, Fudan University, Shanghai 200433, China}
\email{23110840019@m.fudan.edu.cn(Z. Zhang)}
\author{Yi Zhou}
\address{School of Mathematical Sciences, Fudan University, Shanghai 200433, China}
\email{yizhou@fudan.edu.cn(Y. Zhou)}

%    General info


\begin{abstract}
In this paper we prove small-data global well-posedness for the skew mean curvature flow of codimension-two submanifolds of \(\mathbb R^{d+2}\) (\(d\ge5\)) in the critical Besov space. With harmonic coordinates and Coulomb gauge, the flow is formulated as a quasilinear Schrödinger equation for the complex mean curvature coupled to an elliptic system for the geometric and gauge variables. The main difficulty is to control the frequency interactions at critical regularity, where no derivative margin is available. Our argument combines two complementary spacetime estimates derived from the mass and momentum balance laws: a new div-curl lemma introduced by the fourth author yields a bilinear estimate with a half-derivative gain, providing the key control of low-high interactions; while a quasilinear interaction Morawetz estimate provides critical spacetime bounds for comparable and high-high frequency interactions. These estimates coupled with the Gauss-Codazzi structure of the curvature equations yield the unique global solutions to the gauge-reduced system in the critical Besov space, and improves the previous small-data global regularity results.

\end{abstract}

\maketitle

%\section*{This is an unnumbered first-level section head}
%This is an example of an unnumbered first-level heading.

%% The correct journal style for \specialsection is all uppercase; a known bug
%% in amsart.cls prevents this, so input must be uppercase until it is fixed.
%\specialsection*{This is a Special Section Head}
%\specialsection*{THIS IS A SPECIAL SECTION HEAD}
%This is an example of a special section head%
%%%%%%%%%%%%%%%%%%%%%%%%%%%%%%%%%%%%%%%%%%%%%%%%%%%%%%%%%%%%%%%%%%%%%%%%
%\footnote{Here is an example of a footnote. Notice that this footnote
%text is running on so that it can stand as an example of how a footnote
%with separate paragraphs should be written.
%\par
%And here is the beginning of the second paragraph.}%
%%%%%%%%%%%%%%%%%%%%%%%%%%%%%%%%%%%%%%%%%%%%%%%%%%%%%%%%%%%%%%%%%%%%%%%%
%.

\tableofcontents

\section{Introduction}\label{sec1}

  The
\emph{skew mean curvature flow} (short for SMCF) is a distinguished dispersive geometric flow to describe the evolution of the oriented immersed codimension-two submanifolds.  Let $(M^{d+2},\bar g)$ be an oriented
Riemannian manifold, $\Sigma^d$ be an oriented $d$--dimensional manifold, and
\[
 F:[0,T)\times\Sigma\longrightarrow M
\]
be a smooth one-parameter family of immersions.  We write
$\Sigma_t=F(t,\Sigma)$, denote by $g=F^\ast\bar g$ the induced metric, and by
$H(F)$ the mean curvature vector of $\Sigma_t$.  Since the normal bundle
$N\Sigma_t$ has rank two and is oriented, it admits a canonical complex structure
$J(F)$, given by the positive rotation through $\pi/2$ in each normal plane.  The
SMCF can be written as
\begin{align} \label{smcf}
	\begin{cases}
\displaystyle
		\partial_t F=J(F) \mathbf{H}(F), \quad(t, x) \in [0, T) \times \Sigma, \\
		F(0, x)=F_0(x).
	\end{cases}
\end{align}
Equivalently, and more conveniently for analysis, one may prescribe only the normal
component of the velocity:
\begin{align}  \label{smcf1}
	\begin{cases}
\displaystyle
		(	\partial_t F)^{\perp}=J(F) \mathbf{H}(F), \quad(t, x) \in [0, T)\times \Sigma, \\
		F(0, x)=F_0(x).
	\end{cases}
\end{align}
The latter formulation has a larger gauge group: time-dependent changes of
coordinates on $\Sigma_t$ alter the tangential component of $\partial_tF$ only and
therefore leave \eqref{smcf1} unchanged.

The one-dimensional SMCF in $\R^3$ is actually the classical binormal or vortex-filament
flow.  Its physical origin  goes back to Da Rios \cite{DA}, who derived the localized motion
of a thin vortex filament in an ideal incompressible fluid.  Hasimoto \cite{Ha72} subsequently
showed that, after a geometric transform, the filament equation  changes into a
focusing cubic nonlinear Schr\"odinger equation.  This remarkable correspondence
explains the coexistence of geometric, Hamiltonian and dispersive structures in
SMCF and exhibits solitons and complete integrability in one dimension.

Higher-dimensional SMCF is the natural codimension-two analogue of the vortex-filament
dynamics.  Shashikanth \cite{Sha12} related the two-dimensional SMCF to vortex membranes in
four-dimensional incompressible fluids, and Khesin \cite{Khe12} extended the local induction
picture to codimension-two vortex membranes in arbitrary dimension.  SMCF also
arises in the asymptotic dynamics of vortices for Gross-Pitaevskii and related
Ginzburg-Landau models; and the rigorous connections were established, in particular, by
Lin \cite{Lin98,Lin20} and Jerrard \cite{Jerrard} in important geometric configurations.

From the PDE viewpoint, SMCF is a quasilinear Schr\"odinger-type geometric flow.
Unlike the parabolic mean curvature flow, it is nondissipative and Hamiltonian in
nature.  Its study therefore sits at the intersection of geometric analysis,
dispersive PDE, vortex dynamics and gauge theory.  A related central problem is to
understand whether small perturbations of a flat plane remain globally regular at
the regularity dictated by scaling.  This question is especially delicate since
the Schr\"odinger evolution is coupled to nonlocal elliptic equations for the
metric, second fundamental form and gauge fields.

\subsection{Scaling and gauge freedom}

If $F(t,x)$ solves \eqref{smcf1}, then also
\[
 F_\mu(t,x)=\mu^{-1}F(\mu^2t,\mu x).
\]
Consequently the scaling-critical homogeneous
Sobolev regularity for the immersion is
 $\dot H^{\frac d2+1}$,
whereas, since schematically $H=\Delta_gF$, the mean curvature and the complex
curvature variable $\psi$ admit critical regularity
 $s_c=\frac d2-1$.

Gauge freedom is also central to the analysis.  Formulation \eqref{smcf} is invariant
under time-independent reparametrizations, while \eqref{smcf1} permits
time-dependent reparametrizations.  In addition, the oriented normal frame may be
rotated by a phase.  These two freedoms are used to impose harmonic coordinates on
the evolving surface and the Coulomb gauge on the normal connection in this paper,  under which the geometric system becomes a quasilinear
Schr\"odinger equation coupled to an elliptic system, with no derivative loss at
the level required for the present argument.

\subsection{Known results}

While we strive to present an overview of the pertinent results, a number of relevant publications could remain unaddressed given the bounds of our knowledge.
Da Rios~\cite{DA} introduced the vortex-filament equation, and
Hasimoto~\cite{Ha72} transformed it into the cubic nonlinear Schr\"odinger
equation. Geometric and Hamiltonian aspects of binormal and skew mean curvature motion
were developed in, among others, Gomez~\cite{Go04}, Haller--Vizman~\cite{HV04},
Terng~\cite{Te15}, and Terng--Uhlenbeck~\cite{TU}. Lin~\cite{Lin98,Lin20} and Jerrard~\cite{Jerrard} connected vortex
dynamics for Ginzburg--Landau/Gross--Pitaevskii type models with codimension-two
geometric motion. Shashikanth~\cite{Sha12} and Khesin~\cite{Khe12} related
higher-dimensional SMCF to vortex-membrane dynamics.  Khesin--Yang~\cite{KY21}
also exhibited explicit higher-dimensional examples, including finite-time
singular behavior, showing that unrestricted global regularity cannot hold.

For compact submanifolds, Song--Sun~\cite{SS19} proved local existence
for surfaces in $\mathbb R^4$, and Song~\cite{So21} established local
existence, uniqueness, and continuous dependence in arbitrary
dimensions. Huang--Tataru~\cite{HT22, HT221} exploited the enlarged
gauge freedom of the normal velocity formulation and developed low regularity
local theories.   Li~\cite{Liz21,Liz22} established small-data global stability results for
Euclidean planes, together with scattering and wave-operator results under suitable
regularity assumptions. Huang-Li-Tataru~\cite{HLT24} proved small-data global regularity for the skew mean curvature flow in the harmonic/Coulomb gauge for \(d\ge4\). Their result requires the initial mean curvature to be small in \(H^{s_d}\), with
    \begin{equation}\label{sd}
 s_d\ge3 \quad\text{if }d=4, \qquad s_d> \frac{d+1}{2}+\frac{1}{2(d-1)} \quad\text{if }d\ge5.
\end{equation}
Since the scaling-critical index for the curvature is \(s_c=d/2-1\), this remains two derivatives above scaling in dimension four and more than \(3/2\) derivatives above scaling in higher dimensions. Huang-Tataru \cite{HT25} established local well-posedness of the skew mean curvature flow for large initial data in low-regularity Sobolev spaces for dimensions \(d \geq 2\). More recently, Ifrim--Tataru~\cite{IfT25, IfT24} developed robust methods
for cubic quasilinear Schr\"odinger flows, including low-regularity and global
well-posedness, which provide an important broader perspective on the SMCF problem.

\subsection{The harmonic/Coulomb formulation}
Following \cite{HT22,HLT24}, we use harmonic coordinates on the
evolving submanifold and impose the Coulomb condition on the normal
connection. We briefly recall the derivation of the
Schr\"odinger--elliptic system
\eqref{gaugsmcf}--\eqref{eq:B-elliptic}.

 We  work in the Euclidean ambient space $(\R^{d+2},g_{\R^{d+2}})$ and near a
flat $d$--plane.  Write $F_\alpha=\partial_\alpha F$ and
\[
 g_{\alpha\beta}=F_\alpha\cdot F_\beta,\qquad
 g^{\alpha\gamma}g_{\gamma\beta}=\delta^\alpha_\beta.
\]
Choose an oriented orthonormal normal frame $(v_1,v_2)$ satisfying
$Jv_1=v_2$ and $Jv_2=-v_1$, and introduce
\[
 m=v_1+\sqrt{-1} v_2,\qquad
 \lambda_{\alpha\beta}=\kappa_{\alpha\beta}+\sqrt{-1}\tau_{\alpha\beta},
 \qquad
 \psi=g^{\alpha\beta}\lambda_{\alpha\beta}.
\]
Here $\kappa$ and $\tau$ are the two real components of the second fundamental
form, $\lambda$ is the complex second fundamental form, and $\psi$ is the complex
mean curvature.  The normal connection is represented by
\[
 A_\alpha=\partial_\alpha v_1\cdot v_2,\qquad
 B=\partial_tv_1\cdot v_2,
\]
while
\[
 V_\gamma=\partial_tF\cdot F_\gamma
\]
records the tangential velocity generated by the choice of time-dependent
coordinates.  Thus the geometric evolution \eqref{smcf1} is encoded by the dispersive unknown
$\psi$ together with the elliptically determined geometric fields
$\lambda,g,V,A,B$.
 
 As mentioned above we impose harmonic
coordinates on $\Sigma_t$ and Coulomb gauge condition
\begin{align}
\nabla^\alpha A_\alpha=0.
\label{Coulomb_gauge}
\end{align}
We then may differentiate the basis to obtain
\begin{align} \label{eq_difbase}
	\begin{cases}
	\partial_\alpha F_\beta=\Gamma_{\alpha \beta}^\gamma F_\gamma+\operatorname{Re}\left(\lambda_{\alpha \beta} \bar{m}\right), \\
	\partial_\alpha^A m=-\lambda_\alpha^\gamma F_\gamma,
\end{cases}
\end{align}
where $\partial_\alpha^A=\partial_\alpha+\sqrt{-1} A_\alpha$ and $\nabla_\alpha^A =\nabla_\alpha+\sqrt{-1} A_\alpha$. The spatial compatibility conditions for \eqref{eq_difbase} are the
Gauss--Codazzi--Ricci identities:
\begin{align}
R_{\sigma\gamma\alpha\beta}
=&\,
\operatorname{Re}
\left(
\lambda_{\beta\gamma}\overline{\lambda}_{\alpha\sigma}
-
\lambda_{\alpha\gamma}\overline{\lambda}_{\beta\sigma}
\right),
\label{eq:gauss-complex}\\
\nabla_\alpha^A\lambda_{\beta\gamma}
=&\,
\nabla_\beta^A\lambda_{\alpha\gamma},
\label{eq:codazzi-complex}
\end{align}
and
\begin{align}
\nabla_\alpha A_\beta-\nabla_\beta A_\alpha
=
\operatorname{Im}
\left(
\lambda_\alpha^{\gamma}
\overline{\lambda}_{\beta\gamma}
\right).
\label{eq:ricci-normal}
\end{align}
Contracting the Gauss and Codazzi equations yields
\begin{align}
\operatorname{Ric}_{\gamma\beta}
=
\operatorname{Re}
\left(
\lambda_{\gamma\beta}\overline\psi
-
\lambda_{\gamma\alpha}
\overline{\lambda}^\alpha_{\beta}
\right),
\label{eq:ricci-gauss}
\end{align}
and
\begin{align}
\nabla^{A,\alpha}\lambda_{\alpha\beta}
=
\nabla_\beta^A\psi.
\label{eq:codazzi-contracted}
\end{align}
In particular,
they preserve the exact tensorial structure of the Gauss term
\begin{align}
   (\lambda\wedge\lambda)_{\gamma\sigma}:= \lambda_{\gamma\sigma}\overline\psi
-
\lambda_{\alpha\gamma}
\overline{\lambda}^\alpha_{\sigma}, \label{wed.lam}
\end{align}
which also appears in the metric equation in \eqref{eq_ell}, and will be crucial in the subsequent high--high frequency analysis in Lemma \ref{lem_J11t4}.

We now turn to the time evolution part. Since
\(
\mathbf H
=
\operatorname{Re}
\bigl(
\psi\overline m
\bigr),
\)
the skew mean curvature flow \eqref{smcf1} can be written as
\begin{align}
\partial_tF
&=
J(F)\mathbf H(F)+V^\gamma F_\gamma
\nonumber\\
&=
-\operatorname{Im}
\bigl(
\psi\overline m
\bigr)
+
V^\gamma F_\gamma,
\label{Ft}
\end{align}
were \(V\) represents the tangential component as mentioned  above.

Differentiating \eqref{Ft} with respect to the spatial variables and using
\eqref{eq_difbase}, together with \(m\perp F_\alpha\), lead to the evolution
of the moving frame:
\begin{align}
\left\{
\begin{aligned}
\partial_tF_\alpha
=&\,
-\operatorname{Im}
\left[
\left(
\partial_\alpha^A\psi
-
\sqrt{-1}\lambda_{\alpha\gamma}V^\gamma
\right)
\overline m
\right]
\\
&+
\left[
\operatorname{Im}
\left(
\psi\overline{\lambda}_\alpha^{\gamma}
\right)
+
\nabla_\alpha V^\gamma
\right]F_\gamma,
\\
\partial_t^B m
=&\,
-\sqrt{-1}
\left(
\partial^{A,\alpha}\psi
-
\sqrt{-1}\lambda_\gamma^{\alpha}V^\gamma
\right)F_\alpha,
\end{aligned}
\right.
\label{F_alphat}
\end{align}
where
\(
\partial_t^B
:=
\partial_t+\sqrt{-1}\,B.
\)

The spatial and temporal covariant derivatives satisfy
\[
\left[
\partial_t^B,\partial_\alpha^A
\right]m
=
\sqrt{-1}
\left(
\partial_tA_\alpha-\partial_\alpha B
\right)m.
\]
Using \eqref{eq_difbase} and \eqref{F_alphat} and equating the tangential
and normal components, respectively, we obtain the evolution equation for
the complex second fundamental form,
\begin{equation}\label{dt_lambda}
    \pa_t^B \lambda^\sigma_\alpha + \lambda^\gamma_\alpha\left(\operatorname{Im}(\psi \bar{\lambda}^\sigma_\gamma)+\nabla_\gamma V^\sigma\right) = \sqrt{-1}\nabla_\alpha^A\left(\pa^{A,\sigma}\psi -\sqrt{-1}\lambda_\gamma^\sigma V^\gamma\right),
\end{equation}
as well as the compatibility condition
\begin{equation} \label{dt_A}
    \partial_t A_\alpha-\partial_\alpha B=\operatorname{Re}\left(\lambda_\alpha^\gamma \bar{\partial}_\gamma^A \overline{\psi}\right)-\operatorname{Im}\left(\lambda_\alpha^\gamma \bar{\lambda}_{\gamma \sigma}\right) V^\sigma,
\end{equation}
where  $\overline{\partial}_\gamma^A=\partial_\gamma-\sqrt{-1} A_\gamma$.

Under the Coulomb gauge condition \eqref{Coulomb_gauge}, the evolution equation \eqref{dt_lambda} further reduces to the Schr\"odinger type equation for the complex scalar mean curvature   $\psi$
\begin{align*}
	\begin{cases}
		\sqrt{-1} \partial_t \psi+g^{\alpha \beta} \partial_\alpha \partial_\beta \psi
		= \sqrt{-1}(V-2 A)_\alpha \nabla^\alpha \psi 	-\sqrt{-1} \lambda_\sigma^\gamma \operatorname{Im}\left(\psi \bar{\lambda}_\gamma^\sigma\right)\\
		\qquad\qquad\qquad\qquad\qquad\quad+\left(B+A_\alpha A^\alpha-V_\alpha A^\alpha\right) \psi,\\
		\psi(0)=\psi_0,
	\end{cases}
\end{align*}
where $ g^{\alpha \beta}= h^{\alpha \beta} +\delta_{\alpha \beta} $,   \(\lambda\) is complex-valued, whereas \(g,V,A\), and \(B\) are
real-valued.
And the above equation is equivalent to the following divergence form
\begin{align}  \label{gaugsmcf}
	\begin{cases}
	\sqrt{-1} \partial_t \psi+\partial_\alpha(g^{\alpha \beta}  \partial_\beta \psi)
= \mathcal{N},\\
	\psi(0)=\psi_0,
\end{cases}
\end{align}
with the nonlinear term
\begin{align*}
	 \mathcal{N}:=&\,\sqrt{-1}(V-2 A)_\alpha \nabla^\alpha \psi 	-\sqrt{-1} \lambda_\sigma^\gamma \operatorname{Im}\left(\psi \bar{\lambda}_\gamma^\sigma\right)\\
	 &+\left(B+A_\alpha A^\alpha-V_\alpha A^\alpha\right) \psi+\partial_{\alpha} h^{\alpha \beta} \partial_\beta \psi.\nonumber
\end{align*}

It remains to determine the auxiliary variables
\(
(\lambda,g,V,A,B)
\)
at each fixed time. The first equation below is precisely the div-curl
system obtained from the Codazzi identities
\eqref{eq:codazzi-complex} and
\eqref{eq:codazzi-contracted}. The metric equation follows from the Gauss
equation in harmonic coordinates, the \(A\)-equation follows from the
Ricci equation and the Coulomb condition, while the \(V\)-equation is chosen
so that the harmonic coordinate condition is propagated by the flow.
Altogether, we obtain
\begin{align} \label{eq_ell}
	\left\{\begin{aligned}
	& \nabla_\alpha^A \lambda_{\beta \gamma}-\nabla_\beta^A \lambda_{\alpha \gamma}=0, \quad \nabla_\beta^A \psi=\nabla^{A, \alpha} \lambda_{\alpha \beta}=\nabla_{\alpha}^A \lambda^\alpha_\beta, \\
	& g^{\alpha \beta} \partial_{\alpha \beta}^2 g_{\gamma \sigma}= {\left[-\partial_\gamma g^{\alpha \beta} \partial_\beta g_{\alpha \sigma}-\partial_\sigma g^{\alpha \beta} \partial_\beta g_{\alpha \gamma}+\partial_\gamma g_{\alpha \beta} \partial_\sigma g^{\alpha \beta}\right] } \\
	&\qquad\qquad\quad\,\,+2 g^{\alpha \beta} \Gamma_{\sigma \alpha, \nu} \Gamma_{\beta \gamma}^\nu-2 \operatorname{Re}\left(\lambda_{\gamma \sigma} \overline{\psi}-\lambda_{\alpha \gamma} \bar{\lambda}_\sigma^\alpha\right), \\
	& \nabla^\alpha \nabla_\alpha V^\gamma=-2 \nabla_\alpha \operatorname{Im}\left(\psi \bar{\lambda}^{\alpha \gamma}\right)-\operatorname{Re}\left(\lambda_\sigma^\gamma \overline{\psi}-\lambda_{\alpha \sigma} \bar{\lambda}^{\alpha \gamma}\right) V^\sigma \\
	&\qquad\qquad\quad\,\,+2\big(\operatorname{Im} (\psi \bar{\lambda}^{\alpha \beta})+\nabla^\alpha V^\beta\big) \Gamma_{\alpha \beta}^\gamma, \\
	& \nabla^\gamma \nabla_\gamma A_\alpha= \operatorname{Re} (\psi \bar{\lambda}_\alpha^\sigma-\lambda_\alpha^\beta \bar{\lambda}_\beta^\sigma ) A_\sigma+\nabla^\gamma \operatorname{Im}\left(\lambda_\gamma^\sigma \bar{\lambda}_{\alpha \sigma}\right).
\end{aligned}\right.
\end{align}
Applying \(\nabla^\alpha\) to \eqref{dt_A} and using the gauge conditions
give the elliptic equation for \(B\):
\begin{align}
		\nabla^\gamma\nabla_\gamma B
=&\,
-\nabla^\gamma\left[
\operatorname{Re}\left(\lambda_\gamma^\sigma
\bar\partial^A_\sigma\bar\psi\right)
-
\operatorname{Im}\left(\lambda_\gamma^\sigma
\bar\lambda_{\sigma\beta}\right)V^\beta
\right]
\nonumber\\
&+
\left(
2\operatorname{Im}\left(\psi\bar\lambda^{\beta\gamma}\right)
+\nabla^\beta V^\gamma+\nabla^\gamma V^\beta
\right)\partial_\beta A_\gamma.
\label{eq:B-elliptic}
\end{align}
The metric evolves according to
\begin{equation}
    \pa_t g_{\alpha\beta} =  2 \operatorname{Im}(\psi \bar{\lambda}_{\alpha \beta}) + \nabla_\alpha V_\beta + \nabla_\beta V_\alpha.\label{dt_g}
\end{equation}
Finally, we fix the remaining gauge freedom by imposing
\begin{align}
    \lambda(\infty)=0, \quad g(\infty)=I_d, \quad V(\infty)=0, \quad A(\infty)=0, \quad B(\infty)=0 \label{cond_inf}
\end{align}
in the appropriate averaged sense.

For later estimates it is convenient to use
\[
 u_1=\psi,\qquad u_2=\lambda,\qquad
 u_3=\nabla h,\qquad u_4=V,\qquad u_5=A.
\]
These five variables have compatible scaling and can be treated in a unified
multilinear framework, while $B$ admits a different scaling and will be estimated separately.

\subsection{Balance laws, bilinear estimate and interaction Morawetz estimate}

The main difficulty is not the linear Schr\"odinger evolution but the
frequency interactions generated by the quasilinear metric and the elliptic
variables.  A direct use of Strichartz and energy estimates is insufficient at the
critical endpoint: low--high interactions may lose derivatives, while
high-high interactions may generate dangerous low frequencies.  Our proof
exploits two complementary physical space mechanisms.

For the linear model
\[
 \sqrt{-1}\partial_tu+\Delta u=0,
\]
one has the mass and momentum balance laws.  After fixing a direction, say
$x_\gamma$, and integrating over the transverse variables
$\widehat{x}_{\gamma}\in\R^{d-1}$, these identities become a pair of one-dimensional
balance laws in $(t,x_\gamma)$.  Then applying a new \textbf{div-curl} lemma introduced by Zhou \cite{Zhou1}
to the mass density of a low-frequency wave and the directional momentum density
of a high-frequency wave yields, for $\mu\ll\rho$, the model estimate
\begin{equation}\label{eq:bil}
 \sup_\tau
 \left\|
   \|P_\mu u\|_{L^2_{\widehat{x}_{\gamma}}}
   \|P_{\rho,\gamma}u^\tau\|_{L^2_{\widehat{x}_{\gamma}}}
 \right\|_{L^2_{t,x_\gamma}}
 \lesssim
 \mu^{-s}c_\mu(0)\,
 \rho^{-s-\frac12-\sigma}c_\rho(\sigma).
\end{equation}
The gain $\rho^{-1/2}$ is the decisive feature, which is a genuinely bilinear
spacetime gain and cannot be recovered from energy estimates alone.  Moreover,
the directional decomposition $P_{\rho,\gamma}=P_\rho\Xi_\gamma(D)$ allows one to
trade the high-frequency factor for a derivative in the distinguished direction.
Consequently \eqref{eq:bil} is especially effective for the low--high and
high--low paraproducts that arise in the elliptic equations and commutators.

We should mention that the new \textbf{div-curl} lemma has revealed great potential of application to the regularity problem of dispersive equations, such as extremal hypersurface equation \cite{WZ1}, Schr\"{o}dinger map flow \cite{WZ2}, wave equation \cite{WZ3}, wave map \cite{LZ24}, Schr\"{o}dinger equation \cite{SZ}, mKdV/modified Benjamin-Ono equation \cite{TZ} and the Zakharov system \cite{HTZ25}.

For the quasilinear equation the exact conservation laws acquire the source terms
coming from $\mathcal N$, the metric perturbation $h$, and $\partial_tg$.  A major
part of the proof consists in showing that the geometric structure makes these
errors perturbative under the bootstrap assumptions.  Thus the balance laws are
not merely formal identities: they provide the spacetime mechanism that recovers
the half derivative needed to close the critical frequency summations.

The bilinear estimate mentioned above is directional and is naturally tied to separated
frequencies, which is therefore not appropriate for the high-high
interaction terms.  To treat those terms we adapt the interaction Morawetz estimate of Colliander, Keel, Staffilani, Takaoka and Tao \cite{CKSTT} for semilinear Schr\"{o}dinger equation to the present quasilinear system.  Splitting
$x=(y,z)\in\R^3\times\R^{d-3}$ and using the interaction potential associated with
$a(y)=|y|$, the mass and momentum balance laws yield the coercive term
\[
 -\Delta_y^2|y|=8\pi\delta_0,
\]
and hence the interaction Morawetz type estimate
\begin{equation}\label{eq:IM}
 \|P_\rho u\|_{L^4_{t,y}L^2_z}
 \lesssim
 \rho^{-s+\frac14-\sigma}c_\rho(\sigma).
\end{equation}
This estimate is translation invariant, does not require the directional
multiplier $\Xi_\gamma(D)$, and is therefore well suited to comparable-frequency
and high--high interactions.  The $L^4$ spacetime gain supplies precisely the
integrability needed to control quartic expressions generated by the interaction
of two quadratic geometric terms. Interaction Morawetz estimates are some kind of average in space of classical Morawetz estimates (see \cite{Sta, Col}) and work well for the nonlinear Schr\"{o}dinger equation with general initial data. The first example of this kind of estimate was proved in the famous work \cite{CKSTT} for Schr\"{o}dinger equation in $\mathbb{R}^3$, and then has been developed for other dimensional cases (see \cite{CGT,PV,Vis,VZ} and references therein).

In conclusion, the two estimates are complementary as follows:
{\small
\[
 \boxed{\text{new div-curl lemma}\Longrightarrow \text{ bilinear estimate}\Longrightarrow
 \text{low-high/high-low control}},
\]}
whereas
{\small
\[
 \boxed{\text{interaction Morawetz estimate}
 \Longrightarrow
 \text{high-high/comparable-frequency control}}.
\]}
Together with the energy and elliptic estimates, they replace the derivative room
that would otherwise have to be supplied by working above scaling.

\subsection{Structural features used in the proof}

 We should mention that some geometric structures are also indispensable to obtain our main result. 
 
 First, the Codazzi and contracted Codazzi identities form a covariant div-curl system for the complex second fundamental form \(\lambda\). Through the Hodge identity, this yields an elliptic reconstruction of the schematic form
$$ \lambda=\nabla\Delta^{-1}\nabla\psi +\nabla\Delta^{-1}(\text{quadratic geometric terms}), $$
so that the principal map from \(\psi\) to \(\lambda\) is of order zero. In particular, derivatives occurring in the nonlinear source terms can be traded through the div--curl structure and compensated by the elliptic inverse derivative, rather than accumulating on a single high-frequency factor. 

Second, the Gauss equation forces the quadratic curvature contribution to the metric equation to appear in the special combination
$$ (\lambda\wedge\lambda)_{\gamma\sigma} = \lambda_{\gamma\sigma}\overline{\psi} - \lambda_{\alpha\gamma} \overline{\lambda^\alpha{}_{\sigma}}. $$
Combined with the Codazzi representation of \(\lambda\), this tensorial structure exhibits a cancellation when two comparable high frequencies interact to produce a much lower output frequency. This gain is essential in dimension five, as it compensates for the low-frequency singularity of the inverse elliptic operator in the reconstruction of the metric, thereby allowing the high--high-to-low interactions arising in the weak difference estimates to be summed at the critical regularity.

Third,  the metric evolution has a favorable structure for the interaction Morawetz identity in the quasilinear setting: the time derivative of the principal coefficient is given schematically by
$$ \partial_t g=\underbrace{\nabla V}_{differentiated ~ elliptic}+\underbrace{\lambda\psi}_{quadratic~curvature}+\text{lower-order terms}. $$
Thus the coefficient variation term arising in the momentum balance law is either differentiated elliptic or quadratic in the curvature variables, and is therefore perturbative in the critical bootstrap regime.

\subsection{Main results}

We first state global well-posedness for the Schr\"odinger--elliptic
system.

\begin{theorem}
\label{thm1}
Let $d\ge5$ and $s=s_c=\frac d2-1$.  There exists $\varepsilon>0$ sufficiently small such that, whenever
\begin{equation}\label{inid}
 \|\psi_0\|_{\dot B^s_{2,1}(\R^d)}\le\varepsilon,
\end{equation}
the Schr\"odinger equation \eqref{gaugsmcf}, coupled with the elliptic system  \eqref{eq_ell} and \eqref{eq:B-elliptic}, is globally well-posed in
$\dot B^s_{2,1}(\R^d)$.  The resulting solution also satisfies the corresponding
geometric time-evolution equations \eqref{dt_lambda}, \eqref{dt_A} and \eqref{dt_g}.
\end{theorem}

Turning back to the geometric flow gives the following result.

\begin{corollary}[Global regularity for SMCF]
\label{cor1}
Let $d\ge5$ and $s=\frac d2-1$.  If the initial submanifold is a sufficiently small
perturbation of a Euclidean plane and
\[
 \|\nabla(g_0-\Id)\|_{\dot B^s_{2,1}(\R^d)}
 +\|H_0\|_{\dot B^s_{2,1}(\R^d)}
 \le\varepsilon,
\]
then the SMCF \eqref{smcf1} for maps from $\R^d$ into $\R^{d+2}$ admits a
global solution in the harmonic coordinates and Coulomb gauge.
\end{corollary}

 Theorem~\ref{thm1} will be proved by using bootstrap argument.  Starting
from the local solution in the harmonic coordinate and Coulomb gauge, we first assume enlarged versions of the
critical energy, local-smoothing/elliptic, bilinear and interaction Morawetz bounds.
The nonlinear estimates, together with the geometric elliptic equations then improve
the bootstrap constants, and hence the global regularity follows.

To formulate the bootstrap assumptions,  we use the $\sigma$-frequency envelope  introduced in Tao \cite{Tao01}:
\begin{align}
	\|P_{\rho} \psi_0\|_{L^2} \lesssim &~\rho^{-s-\sigma} c_\rho (\sigma), \label{ini1}\\
	\sum_{\rho} c_\rho (\sigma) \lesssim &~\|\psi_0\|_{\dot{B}^{s+\sigma}_{2,1}}, ~~\sigma=0,2\label{ini2},\\
	c_\nu (\sigma)  \lesssim &~c_{\rho } (\sigma) 2^{\delta|l-k|} \label{ini3},
\end{align}
where $ \nu=2^l, \rho = 2^{k},$ and $ c_\rho (\sigma)$ is a constant depending on $ \rho $ for $  \sigma=0,2. $  We shall refer to the condition \eqref{ini3}
as the slow variation property of the $\sigma$-frequency envelope. In particular,
for $\nu\gtrsim \rho,$ this property gives
\[
       \rho^\delta c_\nu(\sigma)
        \lesssim  \rho^\delta c_\rho(\sigma)2^{\delta (l-k)}
        \lesssim  \nu^\delta c_\rho(\sigma).
\]
The case $\sigma=0$ tracks the critical norm, while $\sigma=2$ propagates two
additional derivatives for smooth approximations and higher-regularity persistence. 

Using the notations of Subsection~\ref{subsec.not}, we impose the
following bootstrap assumptions on $u_i$ and $B$:
\begin{align}
	&\|P_{\rho} u_i(t) \|_{L^2_x} \lesssim  \varepsilon^{-\frac18}\rho^{-s-\sigma} c_\rho (\sigma),\label{sup1}\\
	&\left\|P_{\rho }  u_{\mathfrak{T}}\right\|_{L^2_{t,x } }\lesssim   \varepsilon^{-\frac 18} \rho^{-s-1-\sigma} c_\rho(\sigma),\label{supV} \\
	\sup_{\tau}&\,  \big\| \|P_\mu u_i\|_{L^2_{\widehat{x}_\gamma}}   \  \| P_{\rho,\gamma} u_j^\tau\|_{L^2_{\widehat{x}_\gamma}}  \big\|_{L^2_{t,x_\gamma}} \label{sup2} \\
 &\quad \lesssim \varepsilon^{-\frac14} \mu^{-s} c_\mu(0)   \rho^{-s-\frac 12-\sigma} c_\rho(\sigma), ~	\text{ for } \mu\ll \rho,\nonumber\\
			\sup_J &\, \|  P_\rho u_i \|_{L^4_{t,y_J} L^2_{z_J}}  \lesssim \varepsilon^{-\frac18} \rho^{-s+\frac 14-\sigma} c_\rho(\sigma), \label{sup3}\\
	&\left\|  B_{\rho }\right\|_{L^2_{t,x } }\lesssim  \varepsilon^{-\frac 18} \rho^{-s-\sigma} c_\rho(\sigma)\label{supB2},
\end{align}
where $\varepsilon$ is given by \eqref{inid}, $P_{\rho,\gamma} = P_\rho \Xi_\gamma(D)$ with $\gamma = 1,\dots,d$,  and
\[
\mathfrak{T} \in \{3,4,5\},\quad i,j \in \{1,2,3,4,5\},\quad \sigma \in \{0,2\}.
\]
Let \(\mathfrak I_3:=
\left\{
J=\{j_1,j_2,j_3\}\subset\{1,\ldots,d\}:
j_1<j_2<j_3
\right\},\)
then for each \(J\in\mathfrak I_3\), write
\(
x=(y_{J},z_J)\in\mathbb R^3\times\mathbb R^{d-3}\), the norm in \eqref{sup3} means
\[
\| \cdot \|_{L^4_{t,y_J} L^2_{z_J}} := \big\| \| \cdot \|_{L_{z_J}^2} \big\|_{L_{t,y_{J}}^4}.
\]
 Noting that the norms in the left hand side of \eqref{sup2} and \eqref{sup3} are translation invariant.

Based on the bootstrap assumptions \eqref{sup1}--\eqref{supB2}, we can prove
Theorem \ref{thm1} by  establishing the following a priori estimates for the Schr\"odinger--elliptic system  \eqref{gaugsmcf}--\eqref{eq:B-elliptic}
\begin{align}
	&\|P_{\rho} u_i(t)\|_{L^2_x} \lesssim \rho^{-s-\sigma} c_\rho (\sigma),\label{re1}\\
	&\left\|P_{\rho }   u_{\mathfrak{T} }\right\|_{L^2_{t,x } }\lesssim  \rho^{-s-1-\sigma} c_\rho(\sigma) \label{reV}\\
\sup_{\tau} &\big\| \|P_\mu u_i\|_{L^2_{\widehat{x}_\gamma}}   \| P_{\rho,\gamma} u_j^\tau\|_{L^2_{\widehat{x}_\gamma}}  \big\|_{L^2_{t,x_\gamma}} \label{re2} \\
 &\quad \lesssim \mu^{-s} c_\mu(0)   \rho^{-s-\frac 12-\sigma} c_\rho(\sigma), ~	\text{ for } \mu\ll \rho,\nonumber\\
	 	\sup_J &\ \|  P_\rho u_i  \|_{L^4_{t,y_J}L^2_{z_J}} \lesssim  \rho^{-s+\frac 14-\sigma} c_\rho(\sigma), \label{re3}\\
&\left\|B_{\rho }\right\|_{L^2_{t,x } }\lesssim   \rho^{-s-\sigma} c_\rho(\sigma),\label{reB2}
\end{align}
with the same notations as in  \eqref{sup1}--\eqref{supB2}, which may further be stated as the following theorem.
\begin{theorem} \label{thm_a_apr}
		Suppose $ d\geq 5, s=\frac {d}{2}-1 $ and there exists a constant $ \varepsilon>0 $ sufficiently small such that the initial data $ \psi_0 $ satisfies
	\begin{align*}
		\| \psi_0\|_{\dot{B}_{2,1}^s} \leq \varepsilon.
	\end{align*}	
		Then the critical inequalities \eqref{re1}--\eqref{reB2}, corresponding
	to \(\sigma=0\), hold for all \(t>0\), where
	\(u_1=\psi\), \(u_2=\lambda\), \(u_3=\nabla h\), \(u_4=V\), and
	\(u_5=A\) are solutions of the Schr\"odinger--elliptic system
	\eqref{gaugsmcf}--\eqref{eq:B-elliptic}. If, in addition, the initial
	data admit a frequency envelope \(c_\rho(2)\) satisfying
	\eqref{ini1}--\eqref{ini3} with \(\sigma=2\), then the same estimates
	also hold with \(\sigma=2\).
\end{theorem}

The proof proceeds in five mutually reinforcing stages: nonlinear paraproduct and
commutator estimates; energy estimates for the Schr\"odinger component; the
mass--momentum/div--curl bilinear estimate; the quasilinear interaction Morawetz
estimate; and elliptic estimates for $\lambda,h,V,A,B$.  The resulting bounds
strictly improve the bootstrap assumptions, after which the standard continuity
argument gives global control and the usual approximation/difference argument
yields existence, uniqueness and continuous dependence at the critical Besov
regularity.

\section{Preliminaries} \label{sec2}
\subsection{Notations} \label{subsec.not}
Greek indices
$\alpha,\beta,\gamma$ range over $\{1,\ldots,d\}$ and are summed when repeated.
For $x=(x_1,\ldots,x_d)$, let
\[
 \widehat {x_{\gamma}}
 =(x_1,\ldots,x_{\gamma-1},x_{\gamma+1},\ldots,x_d).
\]
We write $|\nabla|=(-\Delta)^{1/2}$ and
$u^\tau(x)=u(x-\tau)$.  \(\partial_x\) denotes an arbitrary first‑order partial derivative \(\partial_{\alpha}\) with  \(\alpha = 1, \dots, d\), while 
 \(\partial_x^2\) denotes an arbitrary second‑order partial derivative \(\partial_{\alpha} \partial_{\beta}\).   

 Let the spatial Fourier transform $\hat{\phi}(t, \xi)$ be defined by
$$
\mathscr{F}\phi =\hat{\phi}(t, \xi):=\int_{\mathbf{R}^n} e^{-2 \pi i x \cdot \xi} \phi(t, x) \, \mathrm{d}x,
$$
with the inverse transform given by
$$
\mathscr{F}^{-1} u = \check{u}(t, x):=\int_{\mathbf{R}^n} e^{2 \pi i x \cdot \xi} u(t, \xi) \, d\xi.
$$

Next, we introduce the following decomposition which corresponds to decomposing the frequency space according to direction. As in the dyadic partition of unity for the Littlewood–Paley decomposition, there exists $\zeta_{\gamma}(\xi) \in C^\infty ( \mathbb{S}^{d-1})$ supported in $  |\xi| \lesssim|\xi_\gamma|$ such that
$ \sum_{\gamma=1}^{d} \zeta_\gamma (\xi)=1, \xi \in \mathbb{S}^{d-1}$
and
 \begin{equation} 
     \text{Id}=\sum_{\gamma=1}^{d} \Xi_\gamma (D) . \label{decomp_Xi}
 \end{equation} 
 
  Let \(\chi(\xi) \in \mathcal{S}(\mathbb{R}^d)\) be a smooth function supported in \(|\xi| \leq 2\) and equal to \(1\) on \(|\xi| \leq 1\). The Littlewood–Paley operators \(P_\rho\) for \(\rho \in 2^{\mathbb{Z}}\) are defined by  
\begin{align}
&\widehat{P_{\leq \rho} u}(\xi):=q\left( \rho^{-1} \xi\right) \widehat{u}(\xi),\quad P_\rho=P_{\leq \rho}-P_{\leq  \rho/2}\nonumber
\end{align}
with symbol  \(\chi_\rho(\xi) = \chi_1(\rho^{-1}\xi)\), where $q_1(\xi)$ supported in \(\{1/2 \leq |\xi| \leq 2\}\) is given by   
\begin{align}\label{LPD}
	&\chi_1\left(\xi\right)=\chi\left(\xi\right)-\chi \left(2\xi\right).
\end{align}

We denote 
\begin{align*}
 & P_{\sim\rho}: =\sum_{\rho^\prime\sim \rho } P_{\rho^\prime},\ P_{\ll \rho} u:=\sum_{\mu \ll \rho } P_\mu u,\,\\
 &P_{\rho,\gamma} := P_{\rho}\Xi_\gamma (D ), \ P_{\sim \rho,\gamma}:=\sum_{\rho^\prime\sim \rho } P_{\rho^\prime,\gamma},
\end{align*}
where \( \rho^\prime\sim \rho=|\log\rho^\prime-\log\rho|\leq 4 \),
and use the short hand 
\begin{align*}
    u_{\ll \rho} := P_{\ll \rho} u,\ v_{\rho,\gamma} := P_{\rho,\gamma}v.
\end{align*}
We will also use Greek indices $\mu, \nu, \mu^\prime, \dots$  to denote the Littlewood–Paley projectors 
\[
P_\mu, P_\nu,  P_{\mu^\prime}\  \text{etc.}  
\]
 where $  \mu, \nu,  \mu^\prime, \in  2^{\mathbb{Z}} $.

The standard Bony’s decomposition reads as
\begin{align} \label{Bony1}
P_\rho (u_i u_j) =&\,P_{\lesssim \rho } u_i P_{\sim\rho} u_j + P_{\lesssim \rho } u_j P_{\sim\rho } u_i +P_\rho \Big(\sum_{\nu^\prime \sim \nu\gtrsim \rho} P_\nu u_i P_{\nu^\prime} u_j\Big),
\end{align}
where the first term on the right-hand side is referred to as the \textit{low--high} term, the second term as the \textit{high--low} term, and the third term as the \textit{high--high} term.

%which can be written alternatively as
%\begin{align} \label{Bony2pri}
%	P_\rho (u_i u_j) = &  P_{\ll \rho } u_i P_{\sim\rho } u_j +   P_{\ll \rho } u_j P_{\sim \rho  } u_i +P_\rho \Big(\sum_{\nu^\prime \sim \nu \gtrsim \rho}  P_\nu u_i P_{\nu^\prime} u_j\Big), 
%  \end{align}
%by the fact that the difference of the first term on the right hand side of \eqref{Bony1} and \eqref{Bony2pri} can be absorbed by the third term in the right hand side of \eqref{Bony2pri}, and so is the difference between the second terms.
 
With a slight abuse of notation, sometimes we denote
\begin{align} \label{Bony2}
	P_\rho (u_i u_j) = P_{\ll \rho } u_i P_\rho u_j + P_{\ll \rho } u_j P_\rho u_i +P_\rho \Big(\sum_{\nu^\prime \sim \nu \gtrsim \rho}  P_\nu u_i P_{\nu^\prime} u_j\Big),
\end{align}
as the finite sum of $  P_{\sim\rho}=\sum_{\rho^\prime\sim \rho } P_{\rho^\prime} $  would not make a difference in this paper.

Let \(f^*\) denote the decreasing rearrangement of \(f\). The Lorentz
space \(L^{p,q}(\mathbb R^d)\) is equipped with the quasi-norm
\[
\|f\|_{L^{p,q}}
:=
\begin{cases}
\displaystyle
\left(
\int_0^\infty
\big(t^{1/p}f^*(t)\big)^q\frac{\mathrm dt}{t}
\right)^{1/q},
& 1\leq q<\infty,\\[1.2ex]
\displaystyle
\sup_{t>0}t^{1/p}f^*(t),
& q=\infty.
\end{cases}
\]
In dimension three, we shall use the Lorentz--Hölder inequality
\[
\int_{\mathbb R^3}|f(y)g(y)|\,\mathrm dy
\lesssim
\|f\|_{L^{2,1}(\mathbb R^3)}
\|g\|_{L^{2,\infty}(\mathbb R^3)},
\]
as well as the Lorentz--Sobolev and endpoint
Hardy--Littlewood--Sobolev inequalities
\[
\big\||\nabla_y|^{-1/2}f\big\|_{L^{2,1}(\mathbb R^3)}
\lesssim
\|\nabla_yf\|_{L^1(\mathbb R^3)},
\quad
\big\||\nabla_y|^{-3/2}g\big\|_{L^{2,\infty}(\mathbb R^3)}
\lesssim
\|g\|_{L^1(\mathbb R^3)}
\]
 in Section \ref{secInter_ac}.

\subsection{Multilinear expressions}  \label{sec2.1}
Following \cite{IfT24, TaoCMP02}, we introduce the \(L\) notation for multilinear expressions.  
Given scalar functions $(\phi_i)_{1\le i\le S}$, we denote by $L (\phi_1, \ldots, \phi_S)$ any multilinear expression of the form
$$
\begin{aligned}
& L (\phi_1, \ldots, \phi_S )(t, x) \\
& \quad:=\int K(\tau_1, \ldots, \tau_{S }) \, \phi_1^{\tau_1}(t, x ) \cdots \phi_S^{\tau_S}(t, x ) \, \mathrm{d}y_1 \cdots  \mathrm{d} y ,
\end{aligned}
$$
where we denote by $\phi^\tau(x) = \phi(x-\tau)$ the translation and $K$ is a measure with bounded mass. The exact value of $K$ may change from line to line, and we require their mass are uniformly bounded by a sufficiently large constant . The kernel $K$ is never allowed to depend on the index $\alpha$.  
We also extend this notation to the case where $(\phi_i)_{1\le i\le S}$ take values as $d$ dimensional vectors or $d \times d$ matrices by allowing $K$ to be an appropriate tensor. 

We introduce the $L$ notation for two reasons. First, all nonlinear terms in the sequel can be expressed uniformly in this notation. Second, the $L$ notation streamlines our multilinear estimates by allowing us to consider only the case where each term is translated — a reduction that is justified by the translation invariance of the spaces on both sides of the estimate.   For example, if we have the bounds with translation invariant norm $\|\cdot\|_X$ and $\|\cdot\|_{X_i}$ for $i =1 ,2 ,3$:
$$\sup_{\tau_\alpha}\| u_1 u_2^{\tau_2} u_3^{\tau_3} \|_{X}\le C\prod_{i=1}^3 \|u_i\|_{X_i},$$
then for any trilinear form $L$ with integrable kernel we have 
\begin{align*}
  \|L(u_1,u_2,u_3)\|_{X} &\le\sup_{\tau_\alpha} \int K(\tau_1,\tau_2,\tau_3)\|u_1^{\tau_1} u_2^{\tau_2}u_3^{\tau_3}\|_{X}\mathrm{d}y \\
  &=\sup_{\tau_\alpha} \int K(\tau_1,\tau_2,\tau_3)\|u_1  u_2^{\tau_2}u_3^{\tau_3}\|_{X}\mathrm{d}y \\
  &\lesssim C\prod_{i=1}^3 \|u_i\|_{X_i}.
\end{align*}
For more details of \(L\) notation, we refer to \cite{IfT24, TaoCMP02}.

We next list some basic properties of the $L$ notation. The $L$ notation interacts well with Littlewood–Paley operators: 
  \begin{align}
& L (\partial_x P_{\leq \rho} \phi_1, \phi_2, \ldots, \phi_S )= \rho  L (\phi_1, \phi_2, \ldots, \phi_S ), \label{L_formula_1}\\
& \partial_x P_{\leq \rho} L (\phi_1, \phi_2, \ldots, \phi_S )=\rho  L (\phi_1, \phi_2, \ldots, \phi_S ),\label{L_formula_2}.
\end{align}

The $L$ notation is also useful as it enables the conversion between derivatives and multipliers, a point we will revisit frequently in the bilinear estimates. Since $P_\rho $ are homogeneous Littlewood-Paley operators, we have 
  \begin{align}
L (  P_\rho \phi_1, \phi_2, \ldots, \phi_S )=\rho^{-1}  L (\partial_x P_\rho\phi_1, \phi_2, \ldots, \phi_S ). \label{L_formula_3} 
\end{align}  
Moreover, we can rewrite $P_{\rho,\gamma}v $ as 
\begin{align*}
    P_{\rho,\gamma} v = \rho^{-1 }  \widecheck{\chi_{\rho,\gamma}} * (\partial_{\gamma}P_{\rho,\gamma} v), 
\end{align*}
where 
\begin{equation}
\chi_{\rho,\gamma}(\xi) = - \sqrt{-1}  \zeta_\alpha\big(\frac{\xi}{|\xi|}\big) \chi_1\big(\frac{\xi}{\rho}\big) \big(\frac{\xi_1}{\rho}\big)^{-1}.
\end{equation}
Using the fact that $\Xi_\gamma(D)$ are homogeneous operators and   $P_{\rho,\gamma} = P_N \Xi_\gamma (D)$  has Fourier support in $\{|\xi |\lesssim |\xi_\gamma|\}$, we conclude that $ \chi_{\rho, \gamma} \in \mathcal{S}(\mathbb{R}^d)$ and hence
\begin{align}
L (  P_{\rho,\gamma} \phi_1, \phi_2, \ldots, \phi_S )=\rho^{-1}  L (\partial_\gamma P_{\rho,\gamma}\phi_1, \phi_2, \ldots, \phi_S ). \label{L_formula_4} 
\end{align}  
Below we give an example of the application of \eqref{L_formula_4} that is most frequently used in the bilinear estimates in the sequel:
\begin{equation}
\begin{aligned}
       & \sup_{\tau} \left\| \|    u_\mu\|_{L^2_{\widehat{x}_\gamma }}  \| P_{\rho,\gamma} v^\tau\|_{L^2_{\widehat{x}_\gamma}}  \right\|_{L^2_{t,x_\gamma}}\\
       &\lesssim \rho^{-1}\sup_{\tau} \left\| \|    u_\mu\|_{L^2_{\widehat{x}_\gamma}}  \| \pa_\gamma P_{\rho,\gamma} v^\tau\|_{L^2_{\widehat{x}_\gamma}}  \right\|_{L^2_{t,x_\gamma}}.
\end{aligned}
\end{equation}

\subsection{Commutator estimates}
Another setting in which the 
\(L\) notation is especially convenient  arises in commutator expressions, namely
\begin{lemma}[Leibniz rule for $P_\rho$ and $P_{\rho,\alpha}$, \cite{TaoCMP02}]  \label{lem_com}
We have the commutator identity
\[
Q_\rho(f g)=f \cdot Q_\rho g+\rho^{-1} L\left(\partial_x f,  g\right),\ Q_\rho = P_\rho,  P_{\rho,\alpha}.
\]
\end{lemma}  
\begin{proof}
    Since  $P_{\rho,\alpha} = P_{\rho}\Xi_\alpha(D)$ where $\Xi_\alpha(D)$ are zero-th order homogeneous operators, the symbols of $P_\rho$ and $ P_{\rho,\alpha}  $  take the form of $\chi(\rho^{-1}\xi)$. Direct computation shows that 
\begin{align*}
    [Q_\rho, f]g &= \int_{\mathbb{R}^d} \rho^d\widecheck{\chi}(\rho y)(f(x-y)-f(x))g(x-y)dy\\
    & = -\rho^{-1}\int^1_0\int_{\mathbb{R}^d} \rho^d (\rho y) \widecheck{\chi}(\rho y)\cdot \pa_x f(x-sy)g(x-y)dyds.
\end{align*}
The conclusion follows from the rapid decay of $\widecheck{\chi}$.
\end{proof}

Applying Lemma \ref{lem_com} and the Bony's decomposition \eqref{Bony1}, we obtain for $Q_\rho = P_\rho$ or $P_{\rho,\alpha}$ the following commutator estimate:
\begin{align}
&[Q_\rho, f]g = Q_\rho(fg) - f Q_\rho g \label{Bony_com} \\
&= \rho^{-1} L(\partial_x P_{\ll \rho} f, P_{\sim \rho} g) + L(P_{\ll \rho} g, P_{\sim \rho} f) + \sum_{\nu \sim \nu' \gtrsim \rho} P_\rho L(P_\nu f, P_{\nu'} g), \nonumber
\end{align}
where Lemma \ref{lem_com} is used for the terms containing low-frequency components of $f$.

\subsection{New type div-curl lemma}
The following \textbf{div-curl} lemma, which was first introduced by the fourth author \cite{Zhou1}, plays a crucial role in our proof.
\begin{lemma} [div-curl Lemma] \label{lemdiv-curl}
	Suppose that $f^{ij}, i, j=1, 2$ satisfy 	
	$$
	\begin{gathered}
	  \begin{cases}
\displaystyle
			\partial_t f^{11}+\partial_x f^{12}=G^1, \\
			\partial_t f^{21}-\partial_x f^{22}=G^2,\\
		\end{cases}  \\
		f^{11}, f^{12}, f^{21}, f^{22} \rightarrow 0, \text{ as  }x \rightarrow \infty,
	\end{gathered}
	$$
	then it holds that	
	\begin{equation}\label{divcurl}
	\begin{aligned}
		&\int_0^T \int_{-\infty}^{+\infty} f^{11} f^{22}+f^{12} f^{21}  \mathrm{ d} x \mathrm{ d} t \\
        \leq&~ 2\max _{0 \leqslant t \leqslant T}\left\|f^{11}(t)\right\|_{L^1}\cdot\max _{0 \leqslant t \leqslant T}\left\|f^{21}(t)\right\|_{L^1}\\
		&+\left|\int_0^T \int_{-\infty}^{+\infty} \left( \int_{-\infty}^x f^{11} (t,y) \mathrm{ d} y\right)G^2(t,x) \mathrm{ d} x \mathrm{ d} t\right| \\
		&+\left|\int_0^T \int_{-\infty}^{+\infty} \left( \int^{+\infty}_x f^{21} (t,y) \mathrm{ d} y\right)G^1(t,x) \mathrm{ d} x \mathrm{ d} t\right|,
	\end{aligned}
	\end{equation}
	provided that the right side is bounded.
\end{lemma}
\begin{proof}
	 The same computation as in \cite{WZ2} yields
	\begin{align*}
	&	\frac{\partial}{\partial t} \int_{x<y} f^{11} (t,x) f^{21} (t,y) \mathrm{ d} x\mathrm{ d} y+\int_{-\infty}^{+\infty}\left (f^{11} f^{22}+f^{12} f^{21}\right) \mathrm{ d} x \\
	=&\,\int_{-\infty}^{+\infty} \left( \int_{-\infty}^x f^{11} (t,y) \mathrm{ d} y\right)G^2(t,x) \mathrm{ d} x  +\int_{-\infty}^{+\infty} \left( \int^{+\infty}_x f^{21} (t,y) \mathrm{ d} y\right)G^1(t,x) \mathrm{ d} x ,
	\end{align*}
    then \eqref{divcurl} follows by integrating the above inequality with respect to $t$ over $[0, T]$.
	\end{proof} 
 
\subsection{Estimates for the linear Schr\"odinger equation}  \label{sec2.2}
To motivate our estimate for the nonlinear equations, we first demonstrate the corresponding estimates for the linear Schr\"odinger equation to help the readers to better understand our proof, thus, we first use the assumption \eqref{sup1}--\eqref{sup3} to obtain the estimates \eqref{re1}--\eqref{re3}  for 
\begin{equation}	 \label{lsq}
	\begin{cases}
		\sqrt{-1} \pa_t u+ \Delta u=0,\\
		u|_{t=0}=u_0.
	\end{cases}
\end{equation}

\begin{lemma}[Energy estimates]
	If the initial data of \eqref{lsq} satisfies \eqref{ini1}, then \eqref{re1} holds for the solution of \eqref{lsq}.
\end{lemma}
\begin{proof}
	Applying the operator 
	 $ P_{\rho} $ to \eqref{lsq}, we get
	\begin{align} \label{linPrho}
			\sqrt{-1} \pa_t u_\rho+ \Delta   u_\rho=0,
	\end{align}
from which the mass conservation law 
	\[\|  u_\rho\|_{L^2}=	\|P_\rho u_0\|_{L^2} \]
	follows directly, and hence yields \eqref{re1} by combining \eqref{ini1}.
\end{proof}

\begin{lemma}[Bilinear estimate] \label{lem_linbie}
	If $  \mu\ll \rho $, and the solution of \eqref{lsq} satisfies \eqref{sup1}--\eqref{sup3}, then \eqref{re2} also holds, i.e.
	\begin{align} \label{linbie}
		& \sup_{\tau} \left\| \| u_\mu\|_{L^2_{\widehat{x}_\gamma}}  \|   P_{\rho,\gamma } u^\tau\|_{L^2_{\widehat{x}_\gamma}}  \right\|_{L^2_{t,x_\gamma}} \\
		\lesssim &\,\mu ^{-s} c_\mu (0) \rho^{-s-\frac 12-\sigma} c_\rho(\sigma),~~\sigma=0, 2.\nonumber 
	\end{align}
\end{lemma}
\begin{proof}
	Without loss of generality, we may assume $ \gamma =1 $, and define $\phi_\rho=   P_{\rho,1} u^\tau $.

Employing $ P_\mu  $ and $ \Xi_1 (D)P_\rho  $ to the equation \eqref{lsq}, we obtain
\begin{align}
	\sqrt{-1} \pa_t u_\mu + \Delta   u_\mu =0, \label{PMlsq}\\
		\sqrt{-1} \partial_t \phi_\rho+ \Delta \phi_\rho=0 \label{xiPNlsq}
\end{align}
respectively. It is easy to get the balance law of mass to \eqref{PMlsq} and the balance law of  momentum to \eqref{xiPNlsq} as follows
	\begin{align*}
		&\frac 12 	\partial_t \int_{\mathbb{R}^{d-1}} |  u_\mu |^2 \mathrm{d} \widehat{x}_1 -\partial_{1} 	\int_{\mathbb{R}^{d-1}} \mathrm{Im}(u_\mu  \pa_{1}\overline{ u}_{\mu })\mathrm{d} \widehat{x}_1=0,\\
		&\partial_t	\int_{\mathbb{R}^{d-1}} \mathrm{Im}( \phi_\rho  \pa_{1}\overline{\phi}_\rho)\mathrm{d} \widehat{x}_1 -\partial_{1} 	\int_{\mathbb{R}^{d-1}} 2 |\pa_{1}\phi_\rho |^2-\partial_{1}^2  \frac{|\phi_\rho|^2}{2}\ \mathrm{d} \widehat{x}_1=0.
	\end{align*}
Applying the div-curl Lemma \ref{lemdiv-curl} to the above balance  laws and by \eqref{L_formula_4}, we can  integrate by parts to obtain
	\begin{align} \label{linbie1}
		\begin{aligned}
			&\sup_{\tau}\int_{[0,T]\times \mathbb{R}} \int_{\mathbb{R}^{d-1}} |  u_\mu |^2 \mathrm{d} \widehat{x}_1\int_{\mathbb{R}^{d-1}}  |\pa_{1}\phi_\rho |^2\mathrm{d} \widehat{x}_1\mathrm{d} x_1\mathrm{d} t\\
			&+\int_{[0,T]\times \mathbb{R}}	\int_{\mathbb{R}^{d-1}} \mathrm{Re} (  u_\mu  \pa_{1}\overline{ u}_\mu ) \mathrm{d} \widehat{x}_1 \int_{\mathbb{R}^{d-1}} \mathrm{Re}(\phi_\rho\pa_{1}\overline{\phi}_\rho)  \mathrm{d} \widehat{x}_1 \mathrm{d} x_1\mathrm{d} t\\
			&-\int_{[0,T]\times \mathbb{R}}		\int_{\mathbb{R}^{d-1}} \mathrm{Im}(  u_\mu  \pa_{1} \overline{ u}_\mu )\mathrm{d} \widehat{x}_1 \int_{\mathbb{R}^{d-1}}\mathrm{Im}(\phi_\rho\pa_{1}\overline{\phi}_\rho) \mathrm{d} \widehat{x}_1 \mathrm{d} x_1\mathrm{d} t\\
:=&\,\widetilde{\mathrm{I}}+\widetilde{\mathrm{II}}+\widetilde{\mathrm{III}}\\ 
			\lesssim &\, \rho\cdot \max_{0\leq s \leq t} \|  u_\mu (s)\|_{L^2_x}^2  \| \phi_\rho(s)\|^2_{L^2_x}.
		\end{aligned}
	\end{align}
	%It follows from $\phi_N$ has Fourier support in $\{|\xi_1|\lesssim |\xi|\}$, $\mu  \ll N$ and Corollary \ref{cor_shift} that
It follows from \eqref{L_formula_4}  and $\mu  \ll \rho$ that
	\begin{align*}
		&	\left|\widetilde{\mathrm{II}}+\widetilde{\mathrm{III}} \right| \\
		\leq	&\,\int_{[0,T]\times \mathbb{R}} \|  u_\mu  \|_{L^2_{\widehat{x}_1}} \|\pa_{1} u_\mu  \|_{L^2_{\widehat{x}_1}} \|\phi_\rho \|_{L^2_{\widehat{x}_1}} \|\pa_{1}\phi_{\rho } \|_{L^2_{\widehat{x}_1}} \mathrm{d} x_1\mathrm{d} t\nonumber\\
		  \lesssim &\, \mu  \rho^{-1}\int_{[0,T]\times \mathbb{R}} \|  u_\mu  \|^2_{L^2_{\widehat{x}_1}}  \|\pa_{1}\phi_\rho \|_{L^2_{\widehat{x}_1}}^2    \mathrm{d} x_1\mathrm{d} t\nonumber\\
		  \ll&\, \widetilde{\mathrm{I}},
	\end{align*}
	which gives the desired estimates since 
	\begin{align*}
		&\int_{[0,T]\times \mathbb{R}}	\int_{\mathbb{R}^{d-1}}|   u_\mu|^2\mathrm{d} \widehat{x}_1 \int_{\mathbb{R}^{d-1}} |\phi_\rho|^2 \mathrm{d} \widehat{x}_1 \mathrm{d} x_1\mathrm{d} t \\
		\lesssim &\,\rho^{-2} \ 	\widetilde{\mathrm{I}} \\
		\lesssim&\, \rho^{-1} \cdot   \max_{0\leq s \leq t} \|  u_\mu (s)\|_{L^2_x}^2  \| \phi_\rho(s)\|^2_{L^2_x}\\
		\lesssim &\,   \mu ^{-2s} c_\mu ^2(0) \rho^{-2s-1-2\sigma} c_\rho^2(\sigma).
	\end{align*}
\end{proof}

We next introduce  the interaction Morawetz estimate \eqref{re3} for \eqref{lsq}.
\begin{lemma}[Interaction Morawetz estimate]  \label{lem_lininta}
	If the initial data of \eqref{lsq} satisfies \eqref{ini1}, then for the solution of \eqref{lsq}, we have
	\eqref{re3}.
\end{lemma}
\begin{proof}Applying the operator 
	$P_\rho $ to the equation \eqref{lsq}, we get a linear equation as \eqref{linPrho}. Thus
	the same manner as in Appendix \ref{appMora} yields
\begin{align}
				 \|   u_\rho  \|^4_{L^4_{t,y}L^2_{z}} \lesssim &\,    \| u_\rho\|_{L^\infty_t L^2_x}^3  \|\pa_x  u_\rho\|_{L^\infty_t L^2_x}\nonumber\\
	\lesssim&\,  \rho^{-4s+1-4\sigma} c^4_\rho (\sigma), 
\end{align}
	which verifies \eqref{re3}.
\end{proof}

\subsection{The general form for the equations and nonlinear terms}  \label{sec2.3}
\subsubsection{Elliptic equations for $\lambda, h, A$}
As in \cite{HLT24,HT22}, we will study the following simplified elliptic equations instead of \eqref{eq_ell}  for convenience
\begin{align} 
	&\partial_\alpha \lambda_{\alpha \beta}=\partial_\beta \psi+A\cdot \psi+h\cdot \partial_x \lambda+\partial_x h\cdot \lambda,\label{ellis0}\\
    &\partial_\alpha \lambda_{\beta \gamma}-\partial_\beta \lambda_{\alpha \gamma}=A \cdot\lambda+\partial_x h \cdot\lambda,\label{ellis01} \\
	&\Delta h= h \cdot\partial_x^2  h+ \partial_x h\cdot \partial_x h +   \lambda\wedge\lambda+h\cdot\partial_x h\cdot \partial_x h, \label{ellis2}\\
	&\Delta V=  h\cdot \partial_x^2 V+\partial_x h\cdot \partial_x V\label{ellis3}\\
  &\qquad \quad +(\partial_x h)^2\cdot V+\lambda^2(A+V+\partial_x h)+\partial_x (\lambda^2), \nonumber\\
	&\Delta A=  h\cdot \partial_x^2 A+\partial_x h\cdot   \partial_x A+(\partial_x h )^2 \cdot A+\lambda^2(A+\partial_x h)+\partial_x (\lambda^2 ), \label{ellis4} 
\end{align}
where  \(\lambda\wedge\lambda\) is defined in \eqref{wed.lam}. Differentiating \eqref{ellis2} once in the spatial variables and recalling
that \(u_3=\partial_xh\), we obtain
\begin{align}
    \Delta u_3= h\cdot \partial_x^2 u_3+ \partial_x h \cdot\partial_x u_3+\partial_x(\lambda\wedge\lambda)+u_3^3, \label{gen_h}
\end{align}
where  we suppress the term $ h \partial_x   (u_3^2) $, since   Lemma \ref{lemma1.3} gives  \( \|h\|_{L^\infty} \lesssim \varepsilon^{\frac 78} \) and thus \(1+h \sim 1\).
Noting that \eqref{ellis0}  is divergence of $ \lambda, $  while \eqref{ellis01} is curl   of $ \lambda, $ we then can
 use  the Hodge decomposition 
 \[
 \Delta U_\alpha= \sum_\beta \partial_{  \beta}^2 U_\alpha = \sum_\beta \partial_{  \beta} \left( \partial_{  \beta} U_\alpha -\partial_{  \alpha} U_\beta\right)+\sum_\beta \partial_{  \beta} \partial_{  \alpha} U_\beta
 \]
 to get
 \begin{align}
 	\Delta \lambda_{\alpha \beta}= \nabla \left( \partial_\beta \psi+A \psi+h \partial_x \lambda+\partial_x h \lambda\right) +\nabla \times \left( A \lambda+\partial_x h \lambda \right). \label{ellis1}
 \end{align}
As we have unified the unknowns into $u_i,\ i=1, 2, \cdots,$ ($ u_1=\psi, u_2=\lambda,u_3=\partial_x h,  u_4=  V, u_5= A $), the equations \eqref{ellis1}, \eqref{ellis3}, \eqref{ellis4} can be rewritten as 
\begin{align}
	\Delta u_2=&\, \partial_x^2 u_1 +\partial_x (u_1 \cdot u_5)+\partial_x (h \cdot \partial_x u_2)+\partial_x (u_2 \cdot (u_3,u_5)) , \label{u_2ell}\\
	\Delta u_4=&\, h\cdot \partial_x^2  u_4 +u_3\cdot \partial_x u_4+u_3^2\cdot u_4+u_2^2 (u_3+u_4+u_5)+\partial_x (u_2^2),\label{u_4ell}\\
    \Delta u_5 =&\,h \cdot\partial_x^2 u_5+u_3 \cdot \partial_x u_5 +u_3^2\cdot  u_5+u_2^2 (u_3+u_5)+\partial_x (u_2^2),\label{u_5ell}
\end{align}
Since $ \partial_x h, V, A $ admit a similar form of elliptic equation and $ u_i $ has the same estimate assumptions \eqref{sup1}--\eqref{sup3}, we can further simplify system \eqref{u_2ell}, \eqref{gen_h},  \eqref{u_4ell} and \eqref{u_5ell}  into the general forms
\begin{align}
	\Delta u_{2}
	&=\partial_x^2 u_1+h \cdot\partial_x^2 u_{2}+ \partial_x h \cdot\partial_x u_{2}+\partial_x(u_i  u_j), \label{Gen_u2}\\	
	\Delta u_{\mathfrak{T}}
	&=h\cdot\partial_x^2 u_{\mathfrak{T}}+ \partial_x h \cdot \partial_x u_\mathfrak{T}+\partial_x(u_2^2)+u_i u_j u_l,  \label{Gen_ellev}
\end{align}
where $\mathfrak{T}\in\{3,4,5\},\ i,j,l\in\{1,2,3,4,5\} $.
\subsubsection{Equations for $B$}
Turning to the unknown function $ B $, we can use the second equation in the first line of \eqref{eq_ell} and  the above rules   to rewrite the  nonlinear term of  \eqref{eq:B-elliptic} as
\begin{align*}
	\operatorname{Re}\left(\lambda_\alpha^\gamma \bar{\partial}_\gamma^A \overline{\psi}\right)= 	&\,	\operatorname{Re}\left(  \partial_\gamma (\lambda_\alpha^\gamma  \overline{\psi} )\right)- 	\operatorname{Re}\left( \partial_\gamma^A  \lambda_\alpha^\gamma \cdot \overline{\psi}\right)\\
	= 	&\,	\operatorname{Re}\left( \partial_\gamma (\lambda_\alpha^\gamma  \overline{\psi} )\right)- 	\operatorname{Re}\left(\nabla^{A}_\gamma \lambda_{\alpha }^\gamma \cdot \overline{\psi}\right)+u_i u_j u_l\\
	=&\, 	\operatorname{Re}\left( \partial_\gamma (\lambda_\alpha^\gamma  \overline{\psi} )\right)- 	\operatorname{Re}\left(\nabla^{A}_\alpha  \psi  \cdot \overline{\psi}\right)+u_i u_j u_l \\
	=&\, 	\operatorname{Re}\left( \partial_\gamma (\lambda_\alpha^\gamma  \overline{\psi} )\right)- \frac12	\partial_\alpha \left| \psi  \right|^2+u_i u_j u_l.
\end{align*}

and consequently to rewrite \eqref{eq:B-elliptic}  as the  following schematic form
\begin{align}
\Delta B
=&\,
h\cdot\partial_x^2B+u_3\cdot\partial_xB
+\partial_x^2(u_i u_j)+\partial_x(u_i u_j u_l)\nonumber \\
&+\partial_xu_4\cdot\partial_xu_5 +u_i u_j\cdot\partial_xu_l+u_i u_j u_k u_l,
\label{eq:B-schematic}
\end{align}
 where $ i,j,k,l \in \{1,2,3,4,5\}. $  

\subsubsection{Schr\"odinger equations for  \texorpdfstring{$\psi$}{psi}}
We  first derive a schematic elliptic equation for the combination
\[
W^\gamma:=V^\gamma-2A^\gamma,
\qquad
A^\gamma=g^{\gamma\alpha}A_\alpha.
\]
Since \(\nabla g=0\), raising the index in the elliptic equation for
\(A_\alpha\) in \eqref{eq_ell} gives
\begin{align*}
\nabla^\beta\nabla_\beta A^\gamma
=&\,
\operatorname{Re}
\left(
\psi\overline{\lambda^{\gamma\sigma}}
-
\lambda^{\gamma\alpha}
\overline{\lambda_\alpha^\sigma}
\right)A_\sigma
 +
\nabla_\alpha\operatorname{Im}
\left(
\lambda^{\alpha\sigma}
\overline{\lambda^\gamma_\sigma}
\right).
\end{align*}
Subtracting twice this equation from the equation for \(V^\gamma\) in
\eqref{eq_ell}, and using
\(
-\operatorname{Im}
\left(
\psi\overline{\lambda^{\alpha\gamma}}
\right)
=
\operatorname{Im}
\left(
\lambda^{\alpha\gamma}\overline{\psi}
\right),
\)
we obtain
\begin{align*}
\nabla^\beta\nabla_\beta W^\gamma
=&\,
2\nabla_\alpha\operatorname{Im}
\left(
\lambda^{\alpha\gamma}\overline{\psi}
-
\lambda^{\alpha\sigma}
\overline{\lambda^\gamma_\sigma}
\right)
+
2\Gamma_{\alpha\beta}^{\gamma}\nabla^\alpha V^\beta
-
\operatorname{Re}
\left(
\lambda^\gamma_\sigma\overline{\psi}
-
\lambda_{\alpha\sigma}\overline{\lambda^{\alpha\gamma}}
\right)V^\sigma
\\
&-
2\operatorname{Re}
\left(
\psi\overline{\lambda^{\gamma\sigma}}
-
\lambda^{\gamma\alpha}\overline{\lambda_\alpha^\sigma}
\right)A_\sigma
 +
2\Gamma_{\alpha\beta}^{\gamma}
\operatorname{Im}
\left(
\psi\overline{\lambda^{\alpha\beta}}
\right).
\end{align*}
By the definition of the Gauss term in \eqref{wed.lam}, we have
\[
(\lambda\wedge\lambda)^{\alpha\gamma}
=
\lambda^{\alpha\gamma}\overline{\psi}
-
\lambda^{\alpha\sigma}
\overline{\lambda^\gamma_\sigma},
\]
the preceding identity takes the schematic form
\begin{align*}
\nabla^\beta\nabla_\beta W^\gamma
=
2\nabla_\alpha\operatorname{Im}
\left(
(\lambda\wedge\lambda)^{\alpha\gamma}
\right)
+
2\Gamma_{\alpha\beta}^{\gamma}\nabla^\alpha V^\beta
+
u_i u_j u_l.
\end{align*}
Expanding the covariant derivatives in Euclidean coordinates and using
\[
g^{\alpha\beta}
=
\delta^{\alpha\beta}+h^{\alpha\beta},
\qquad
\Gamma=L(g^{-1},\partial_xh),
\]
we obtain schematically
\begin{align}
\Delta W
=&\,
h\cdot\partial_x^2W
+
\partial_x(\lambda\wedge\lambda)
+
\partial_xh\cdot\partial_xV
+
\partial_xh\cdot\partial_xW
+
u_i u_j u_l
\nonumber\\
=&\,
h\cdot\partial_x^2W
+
\partial_x(\lambda\wedge\lambda)
+
\partial_xh\cdot\partial_xu_{\mathfrak T}
+
u_i u_j u_l
\label{eq:W-final-schematic}
\end{align}
 with
\(\mathfrak T\in\{3,4,5\}\).

 With a slight abuse of notation, we write the nonlinearity
\(\mathcal N\) in \eqref{gaugsmcf} in the schematic form
\begin{align} \label{Gen_N}
		\mathcal{N}=&\, (W+u_3) \cdot \partial_x \psi+\psi  u_i u_j +  \psi  B\nonumber\\
		=&\, \mathbf{U} \cdot \partial_x \psi+\psi  u_i u_j +  \psi  B,
	\end{align}
where  $ i,j\in\{1,2,3,4,5\} $, $  B $ is a real function and \(\mathbf{U}:=(W,u_3)\). It follows from \eqref{gen_h} and
\eqref{eq:W-final-schematic} that \(\mathbf U\) satisfies the schematic
elliptic equation
\begin{align}
    \Delta \mathbf{U} =h\cdot\partial_x^2 \mathbf{U}
+
 \partial_xh\cdot\partial_xu_{\mathfrak T}+
\partial_x(\lambda\wedge\lambda)
+
u_i u_j u_l \label{eq_U}
\end{align}
 for \(\mathfrak{T}\in\{3,4,5\}\) and $ i,j,l\in\{1,2,3,4,5\} $.

The term $\psi B$ will be treated separately. Indeed, $B$ has a different
scaling from the variables $u_i$ and satisfies a less favorable elliptic
equation, so a direct product estimate for $Q_\rho(\psi B)$ is not suitable.
Instead, we use the special structure of the Schr\"odinger identities. In the
energy estimate, taking the imaginary part removes the principal term
$B Q_\rho\psi$. In the momentum type estimates, taking the real part allows
us to rewrite the contribution of $B Q_\rho\psi$ as a divergence term plus a
term where the derivative falls on $B$. Thus only the commutator
$[Q_\rho,B]\psi$ and differentiated terms involving $\partial_xB$ need to be
estimated.

Applying the operator $Q_\rho = P_{\rho}, P_{\rho,\gamma}$ to \eqref{gaugsmcf} gives
\begin{equation}\label{Q_psi}
\sqrt{-1}\pa_t \, Q_\rho \psi + \pa_\alpha (g^{\alpha\beta}\pa_\beta \, Q_\rho\psi) = Q_{\rho} \mathcal{N}- \partial_{ \alpha}  ([Q_{\rho}, h^{\alpha \beta}]  \partial_{ \beta}\psi   ),
\end{equation}
 where   
\begin{align}
\label{non_B}
    Q_\rho \mathcal{N} = &\,Q_\rho (\mathbf{U} \cdot \partial_x \psi+\psi  u_i u_j)  + [Q_\rho,B]\psi + B\cdot Q_\rho \psi.
 \end{align}
 Since $B$ is real-valued, the energy identity gives
\begin{equation}
         \operatorname{Im}(Q_\rho \overline{\psi}\cdot  Q_\rho \mathcal{N})= \operatorname{Im}\big( Q_\rho \overline{\psi} \cdot  \big(Q_\rho (\mathbf{U} \cdot \partial_x \psi+\psi  u_i u_j)  + [Q_\rho,B]\psi \big)\label{energy_B}\big),
\end{equation}
where the main term $B|Q_\rho\psi|^2$ disappears.

For the momentum type identities, we have 
 \begin{align}
      \operatorname{Re}\big(\pa_{\alpha'}Q_\rho \overline{\psi}\cdot  Q_\rho \mathcal{N}\big)& = \operatorname{Re}\big( \pa_{\alpha'} Q_\rho \overline{\psi} \cdot  \big(Q_\rho (\mathbf{U} \cdot \partial_x \psi+\psi  u_i u_j)  + [Q_\rho,B]\psi \big) \big)  \nonumber\\
      &\quad + \frac 12\pa_{\alpha'}(B \cdot |Q_\rho \psi|^2) - \frac 12\pa_{\alpha'}B \cdot |Q_\rho \psi|^2,\label{momen_B_1}\\
      \operatorname{Re}(Q_\rho \overline{\psi}\cdot  \pa_{\alpha'}Q_\rho \mathcal{N})& = \operatorname{Re}\big( Q_\rho \overline{\psi} \cdot \pa_{\alpha'}  \big(Q_\rho (\mathbf{U} \cdot \partial_x \psi+\psi  u_i u_j)  + [Q_\rho,B]\psi \big) \big) \nonumber\\
      &\quad +\frac 12 \pa_{\alpha'}(B \cdot |Q_\rho \psi|^2) + \frac 12\pa_{\alpha'}B \cdot |Q_\rho \psi|^2.  \label{momen_B_2}
 \end{align}
Subtracting the two identities, the divergence terms
\(\partial_{\alpha'}(B|Q_\rho\psi|^2)\) cancel and schematically, we obtain
 \begin{equation}\label{momen_B_3}\begin{aligned}
      &\operatorname{Re}(\pa_{\alpha'}Q_\rho \overline{\psi}\cdot  Q_\rho \mathcal{N}) - \operatorname{Re}(Q_\rho \overline{\psi}\cdot  \pa_{\alpha'}Q_\rho \mathcal{N}) \\
      & = \rho\,  L(Q_\rho \psi, Q_\rho (\mathbf{U} \cdot \partial_x \psi+\psi  u_i u_j)+ [Q_\rho , B]\psi )\\
      &\quad  +   \pa_{\alpha'}\big( Q_\rho \overline{\psi} \cdot ([Q_\rho, B]\psi)\big)-\pa_{\alpha'}B\cdot |Q_\rho \psi|^2. 
 \end{aligned}
 \end{equation}
These identities will be used repeatedly below to avoid estimating
$Q_\rho(\psi B)$ directly.

\section{Nonlinear estimates} \label{sec_non_est}

This section collects the nonlinear estimates that will be used repeatedly
in the \textit{energy estimate, bilinear estimate, interaction Morawetz estimate} and \textit{elliptic estimate}.  Such a summary will help to streamline the analysis in sections \ref{secEn}-\ref{sec_ell}.

 The estimates throughout this paper are controlled by the following dyadic summation principle. After applying the relevant estimates and transferring all frequency envelopes to the reference frequency via \eqref{ini3}, it remains only to examine the sign of the residual dyadic power. If the remaining power of a low frequency is positive, then the corresponding low-frequency sum is dominated by its largest admissible frequency, typically the reference frequency \(\rho\). If the remaining power of a high frequency is negative, then the high-frequency tail is dominated by its smallest admissible frequency, again typically \(\rho\). The small loss introduced by \eqref{ini3} is harmless, provided that \(\delta>0\) is chosen sufficiently small.

We first record some direct consequences of the bootstrap assumptions in
Subsection \ref{subsec_non_est1}, then estimate the typical nonlinear
expressions appearing in the main a priori estimates in Subsection
\ref{subsec_non_est2}.

\subsection{Preliminary estimates from the bootstrap assumptions}
\label{subsec_non_est1}

We begin with several estimates that follow directly from the bootstrap
assumption.  They will be used to place low-frequency factors in $L^\infty$-type
norms and to make the dyadic summations converge.

The first lemma gives the basic Besov and $L^\infty$ controls for the
unknowns $u_i$ and the low-frequency part of $B$.

\begin{lemma} \label{lemma1.3}
Assume that  \eqref{sup1} and \eqref{supB2} hold. Then
	\begin{align}
		\| u_i(t) \|_{\dot{B}^s_{2,1}}&\,  \lesssim \varepsilon^{\frac 78} \label{u_Bes},\\
		\||\nabla|^{-1}u_i\|_{L^\infty} 	&\, \lesssim \varepsilon^{\frac 78},  \label{h_infty}\\
      	\| P_\mu u_i\|_{L^\infty} 	&\, \lesssim \varepsilon^{\frac 78} \mu,  \label{u_infty}  \\
   \| P_{\lesssim\rho } u_i\|_{L^\infty} 	&\, \lesssim \varepsilon^{\frac 78} \rho,  \label{u_infty_rho}       \\
  \| B_{\lesssim\rho} \|_{L^2_t L^\infty} &\,\lesssim  \varepsilon^{\frac 78} \rho, \label{B_infty}\\
  \left\| \partial_x B_{\lesssim \rho  }\right\|_{L^2_{t,y }L^\infty_z }  &\, \lesssim  \varepsilon^{\frac 78} \rho^{\frac 12}, \label{Btx2z_inf}
	\end{align}
	where  $ s=\frac {d}{2}-1, i=1,2,3,4,5 $.
\end{lemma}
\begin{proof} By \eqref{ini2} and \eqref{sup1}, it is easy to see
	\begin{align*}
		\| u_i(t) \|_{\dot{B}^s_{2,1}}&= \sum_{\rho} \rho^{ s} \|P_\rho u_i(t) \|_{L^2}\nonumber\\
		&\lesssim \sum_{\rho} c_\rho (0) \varepsilon^{-\frac 18} \lesssim \varepsilon^{-\frac 18 } \| \psi_0\|_{\dot{B}^s_{2,1}} \lesssim \varepsilon^{\frac 78},
	\end{align*}
    which is exactly \eqref{u_Bes}. Bernstein's inequality and \eqref{u_Bes} then give
	\begin{align*}
		\||\nabla|^{-1}u_i\|_{L^\infty} &\lesssim  \sum_{ \rho} \rho^{-1}  \|P_\rho u_i(t) \|_{L_x^\infty}\\
		&\lesssim\sum_{\rho} \rho^{ s} \|P_\rho u_i(t) \|_{L_x^2}=\| u_i(t) \|_{\dot{B}^s_{2,1}}\lesssim \varepsilon^{\frac 78}.
	\end{align*}
Estimates \eqref{u_infty}--\eqref{Btx2z_inf} follow similarly. For instance, Bernstein's inequality and \eqref{supB2} yield
\begin{align*}
  \left\| \partial_x B_{\lesssim \rho  }\right\|_{L^2_{t,y }L^\infty_z } \lesssim&\, \sum_{\mu \ll \rho }  \left\| \partial_x B_{\mu }\right\|_{L^2_{t,y }L^\infty_z }  \nonumber \\
	\lesssim &\,  \sum_{\mu \ll \rho }   \mu^{\frac{d-3}{2}} \mu\left\|   B_{\mu }\right\|_{L^2_{t,x } }   \nonumber\\
	\lesssim &\,  \varepsilon^{-\frac 18}  \sum_{\mu \ll \rho } \mu^{s+\frac 12} \mu^{-s} c_\mu (0)  
 \lesssim    \varepsilon^{\frac 78} \rho^{\frac 12}, 
\end{align*}
which is \eqref{Btx2z_inf}.
\end{proof}	
We next record the complementary mixed-norm consequences of the interaction Morawetz bootstrap bounds \eqref{sup3}.
 \begin{lemma} Assume that  \eqref{sup1}, \eqref{sup2} and \eqref{sup3} hold. Then
     \begin{align}
      & \left\|    P_{\nu }   \psi \right\|_{L^4_{t,y}L^\infty_z}    
       \lesssim \varepsilon^{-\frac18}  \nu^{-\frac 14-\sigma} c_\nu (\sigma), \label{pre.IntM.infty} \\
        &   \left\|    P_{\ll \rho } \partial_x \psi \right\|_{L^4_{t,y}L^\infty_z}    
       \lesssim \varepsilon^{\frac78}   \rho^{ \frac 34}, \label{pre.IntM2}\\
           & \sum_{\nu^\prime \sim \nu \gtrsim \rho}   \left\|	L(P_{ \nu^{ \prime}} u_i, P_{\nu}u_j)\right\|_{L^2_{t,x}}  \lesssim \,  \varepsilon^{\frac34}     \rho^{-s-\sigma} c_{\rho} (\sigma), \label{pre.IntM1}\\
       & \sum_{  \mu \ll \rho} \sum_{  \nu\sim \nu^\prime \gtrsim \rho }  \mu^{-\frac 32} \Big\| \|    	 P_{  \mu } L(P_\nu u_i, 
 P_{\nu'} u_j) \|_{L^\infty_z} \|   P_{ \rho} u_l \|_{L^2_z}     \Big\|_{L^2_{t,y}}  \label{pre.IntM.quad}
 \\ \lesssim &\, \varepsilon^{\frac 12} \rho^{-s-\frac 12-2\sigma} c_\rho^2(\sigma), 
\nonumber
     \end{align} 
 \end{lemma}
\begin{proof}
\eqref{pre.IntM.infty} follows directly from \eqref{sup3}:
    \begin{align*}
        \left\|    P_{\nu }  \psi \right\|_{L^4_{t,y}L^\infty_z}    
		\lesssim&\,     \nu^{\frac{d-3}{2}}  \left\|  \psi_{\nu }   \right\|_{L^4_{t,y}L^2_z} 
    \end{align*}
    Similarly, Bernstein's inequality in the $z$ variables and \eqref{sup3} imply \eqref{pre.IntM2}:
    \begin{align*}
   \left\|    P_{\ll \rho } \partial_x \psi \right\|_{L^4_{t,y}L^\infty_z}    
		\lesssim&\,    \sum_{\mu\ll \rho} \mu \cdot \mu^{\frac{d-3}{2}}  \left\|  \psi_{\mu }   \right\|_{L^4_{t,y}L^2_z} \\
		\lesssim &\,\varepsilon^{-\frac 18}  \sum_{\mu\ll \rho}   \mu^{s+\frac 12}  \mu^{-s+\frac 14} c_\mu(0)  \\
		\lesssim &\, \varepsilon^{\frac78}   \rho^{ \frac 34}.
    \end{align*}

Estimate \eqref{pre.IntM1} follows from H\"older's inequality, \eqref{pre.IntM.infty}, and \eqref{sup3}:
    \begin{align*}
	\sum_{\nu^\prime \sim \nu \gtrsim \rho}   \left\|	L(P_{ \nu} u_i, P_{\nu'}u_j)\right\|_{L^2_{t,x}} \lesssim&\,  	 	\sum_{\nu^\prime \sim \nu \gtrsim \rho}   \left\|	P_{ \nu^{ \prime}} u_i \right\|_{L^4_{t,y}L^\infty_z}   \left\|P_{\nu}u_j\right\|_{L^4_{t,y}L^2_z} \\
	\lesssim &\, \varepsilon^{-\frac 14}	\sum_{\nu^\prime \sim \nu \gtrsim \rho}   \nu^{-\frac 14-\sigma} c_\nu (\sigma) (\nu^{\prime})^{-s+\frac 14} c_{\nu^\prime} (0) \\
%	\lesssim &\, \varepsilon^{\frac34}   \sum_{\nu^\prime \sim \nu \gtrsim \rho}   \nu^{-s-\sigma+\delta} \rho^{-\delta} c_{\rho} (\sigma).
\lesssim&\, \varepsilon^{\frac34}     \rho^{-s-\sigma} c_{\rho} (\sigma),
	\end{align*}
which with \eqref{sup1}, H\"older's inequality, and Bernstein's inequality yields \eqref{pre.IntM.quad}: 
\begin{align*}
    & \sum_{  \mu \ll \rho} \sum_{  \nu\sim \nu^\prime \gtrsim \rho }  \mu^{-\frac 32}  \Big\| \|    	 P_{  \mu } (P_\nu u_i 
 P_{\nu'} u_j) \|_{L^\infty_z} \|   P_{ \rho} u_l  \|_{L^2_z}      \Big\|_{L^2_{t,y}} \\
 \lesssim &\,   \sum_{  \mu \ll \rho} \sum_{  \nu\sim \nu^\prime \gtrsim \rho }   \mu^{-\frac 32} \|    	 P_{  \mu } (P_\nu u_i 
 P_{\nu'} u_j) \|_{L^\infty_x L^2_t} \|   P_{ \rho} u_l  \|_{L^2_x}       \\
 \lesssim &\, \sum_{  \mu \ll \rho}  \sum_{  \nu\sim \nu^\prime \gtrsim \rho }   
 \mu^{ \frac  d2-\frac 32 }  \|    	 P_\nu u_i 
 P_{\nu'} u_j \|_{L^2_{t,x}} \|   P_{ \rho} u_l  \|_{L^2_x}  \\
 \lesssim &\, \varepsilon^{\frac 12} \rho^{-s-\frac 12-2\sigma} c_\rho^2(\sigma).
\end{align*}

\end{proof}

We then collect several mixed-norm estimates derived primarily from the bilinear bootstrap bounds \eqref{sup2}.
\begin{lemma} Assume that  \eqref{sup1}, \eqref{sup2} and \eqref{sup3} hold. Then
    \begin{align}
        &  \sum_{\mu\ll \rho}  \sum_{\gamma=1}^d \mu^{-\frac 12} \big\| \left\|P_{\mu } u_i\right\|_{L^\infty_{\widehat{x}_\gamma}} \left\|P_{  \rho,\gamma}    u_j\right\|_{L^2_{\widehat{x}_\gamma}}\big\|_{L^2_{t,x_\gamma}} \lesssim  
 \varepsilon^{\frac 34}     \rho^{-s-\frac12-\sigma} c_{\rho} (\sigma), \label{pre.bi1}\\
 & \sum_{\mu \ll \rho} \mu^{-1} \rho \sup_\tau\Big\| \|    	L(P_{ \mu} u_i, P_{\ll \mu }u_j)  \|_{L^\infty_z} \|    \psi_{\sim \rho,1} \|_{L^2_z}   \|     \psi _{\rho,1}^\tau\|_{L^2_z}     \Big\|_{L^1_{t,y}}   \nonumber  \\
 \lesssim&\, \varepsilon^{\frac32}    \rho^{-2s  -2\sigma}c_{\rho }^2 (\sigma), \label{pre.bi2} \\
 & \sum_{\mu \ll \rho} \mu^{-\frac 32} \sup_\tau  \Big\| \|    	 P_{  \mu } (u_i u_j) \|_{L^\infty_z} \|   P_{ \rho,1} u_l^\tau \|_{L^2_z}     \Big\|_{L^2_{t,y}}   \nonumber  \\
 \lesssim&\, \varepsilon^{\frac 32} \rho^{-s-\frac12-\sigma} c_\rho(\sigma), \label{pre.bi3} 
    \end{align}
    where we set \(y_1=x_1\) in \(y=(y_1,y_2,y_3)\),  $ i,j,l \in\{1,2,3,4,5\} $.
\end{lemma}

\begin{proof}
  Estimate \eqref{pre.bi1} follows directly from \eqref{sup2} and Bernstein's inequality:
\begin{align*}
 &  \sum_{\mu\ll \rho}  \sum_{\gamma=1}^d \mu^{-\frac 12} \big\| \left\|P_{\mu } u_i\right\|_{L^\infty_{\widehat{x}_\gamma}} \left\|P_{  \rho,\gamma}    u_j\right\|_{L^2_{\widehat{x}_\gamma}}\big\|_{L^2_{t,x_\gamma}}\\
	\lesssim &\, 	 \sum_{\mu\ll \rho}\sum_{\gamma=1}^d \mu^{\frac{d}{2}-1} \big\| \left\|P_{\mu} u_i\right\|_{L^2_{\widehat{x}_\gamma}} \left\|P_{ \rho,\gamma}    u_j\right\|_{L^2_{\widehat{x}_\gamma}}\big\|_{L^2_{t,x_\gamma}}\\
	\lesssim &\,  \sum_{\mu\ll \rho}  \varepsilon^{-\frac 14}\mu^{s-1} c_\mu(0) \rho^{-s-\frac12-\sigma} c_{\rho} (\sigma) \\
	\lesssim &\,\varepsilon^{\frac 34}   \rho^{-s-\frac12-\sigma} c_{\rho} (\sigma),
\end{align*}
where \(\sum_{\mu \ll \rho}   c_{\mu} (0) \lesssim \varepsilon\). Then, the \eqref{pre.bi1} gives   \eqref{pre.bi2} immediately:
\begin{align*}
&  \sum_{\mu \ll \rho} \sup_\tau \mu^{-1}\rho \Big\| \|    	L(P_{ \mu} u_i, P_{\ll \mu }u_j)  \|_{L^\infty_z} \|    \psi_{\sim \rho,1} \|_{L^2_z}   \|     \psi _{\rho,1}^\tau\|_{L^2_z}     \Big\|_{L^1_{t,y}}   \\
\lesssim&  \,\rho   \sum_{\mu \ll \rho}  \mu^{-\frac 12}	\sup_\tau  \Big\| \|    	P_{ \mu} u_i  \|_{L^\infty_{\widehat{x}_1}}        \|     \psi _{\rho,1}^\tau\|_{L^2_{\widehat{x}_1}}     \Big\|_{L^2_{x_1}}    \cdot  \sum_{\omega \ll \mu \ll \rho }	 	\sup_{\tau'}  \Big\| \|    	P_{ \omega} u_j  \|_{L^\infty_{\widehat{x}_1}}        \|     \psi _{\sim\rho,1}^{\tau'}\|_{L^2_{\widehat{x}_1}}     \Big\|_{L^2_{x_1}}\\
 	\lesssim &\, \varepsilon^{\frac32}   \rho^{-2s -2\sigma}c_{\rho }^2 (\sigma),  \nonumber
\end{align*}

To prove \eqref{pre.bi3}, Bony decomposition gives
\begin{align}
     &\sum_{\mu \ll \rho} \mu^{-\frac 32} \sup_\tau  \Big\| \|    	 P_{  \mu } (u_i u_j) \|_{L^\infty_z} \|   P_{ \rho,1} u_l^\tau \|_{L^2_z}     \Big\|_{L^2_{t,y}}    \nonumber \\
     \lesssim &\,  \sum_{  \mu \ll \rho} \mu^{-\frac 32}\sup_\tau  \Big\| \|    	 P_{  \mu } u_i 
 P_{\ll \mu} u_j \|_{L^\infty_z} \|   P_{ \rho,1} u_l^\tau \|_{L^2_z}     \Big\|_{L^2_{t,y}} \label{pre.bi3.t1}\\
 &\,+  \sum_{\mu \lesssim \nu\sim \nu^\prime \ll \rho }\mu^{-\frac 32}\sup_\tau  \Big\| \|    	 P_{  \mu } (P_\nu u_i 
 P_{\nu'} u_j) \|_{L^\infty_z} \|   P_{ \rho,1} u_l^\tau \|_{L^2_z}     \Big\|_{L^2_{t,y}} \label{pre.bi3.t2} \\
 &\,+\sum_{  \mu \ll \rho} \sum_{  \nu\sim \nu^\prime \gtrsim \rho } \mu^{-\frac 32}\sup_\tau  \Big\| \|    	 P_{  \mu } (P_\nu u_i 
 P_{\nu'} u_j) \|_{L^\infty_z} \|   P_{ \rho,1} u_l^\tau \|_{L^2_z}     \Big\|_{L^2_{t,y}}, \label{pre.bi3.t3}
\end{align}
where  \eqref{pre.bi3.t1} can be estimated by \eqref{pre.bi1}  and  \eqref{u_infty_rho} by taking infty norm of \(P_{\ll \mu} u_j \), while   \eqref{pre.bi3.t3} has been estimated in \eqref{pre.IntM.quad}.
As for \eqref{pre.bi3.t2}, we can also use \eqref{pre.bi1}  and  \eqref{u_infty_rho} to obtain
\begin{align*}
     \eqref{pre.bi3.t2} \lesssim &\,  \sum_{\mu \lesssim \nu\sim \nu^\prime \ll \rho }\mu^{d-\frac 52}\sup_\tau  \Big\| \|    	 P_{  \mu } (P_\nu u_i 
 P_{\nu'} u_j) \|_{L^1_{\widehat{x}_1}} \|   P_{ \rho,1} u_l^\tau \|_{L^2_{\widehat{x}_1}}     \Big\|_{L^2_{t,x_1}} \\
 \lesssim &\,  \sum_{\mu \lesssim \nu\sim \nu^\prime \ll \rho }\mu^{2s-\frac 12}\sup_\tau  \Big\| \|    	 P_\nu u_i 
\|_{L^2_{\widehat{x}_1}} \|   P_{ \rho,1} u_l^\tau \|_{L^2_{\widehat{x}_1}}     \Big\|_{L^2_{t,x_1}} \cdot \| P_{\nu'} u_j \|_{L^2_{\widehat{x}_1} L^\infty_{x_1}} \\
\lesssim &\,  \varepsilon^{-\frac 38} \sum_{\mu \lesssim \nu\sim \nu^\prime \ll \rho }\mu^{2s-\frac 12}  \nu^{-s} c_\nu (0 ) \rho^{-s-\frac12-\sigma} c_\rho(\sigma) (\nu')^{-s+\frac 12} c_{\nu'} (0) \\
\lesssim &\, \varepsilon^{\frac 32} \rho^{-s-\frac12-\sigma} c_\rho(\sigma),
\end{align*}
where \(\sum_{\nu \ll \rho}   c_{\nu} (0) \lesssim \varepsilon\).

\end{proof}

\subsection{Typical nonlinear estimates} \label{subsec_non_est2}
We now estimate the nonlinear expressions that arise in the subsequent a priori estimates.

The first lemma establishes space--time estimates for the typical quadratic nonlinearity $u_i u_j$.
\begin{lemma}  \label{lem_l2uiuj}Suppose $ i,j\in\{1,2,3,4,5\} $. Then we have
	\begin{align}\label{lem_l2uiuj_est}
		\left\| P_\rho	(u_i  u_j ) \right\|_{L^2_{t,x}}  \lesssim \varepsilon^{\frac34}\rho^{-s-\sigma} c_\rho(\sigma).
	\end{align}
\end{lemma}
\begin{proof}
The Bony's decomposition yields
\begin{equation}\label{Bon2}
	\begin{aligned}   
		P_{\rho} \left(u_i   u_j \right) 
		=&\,     P_{\ll \rho} u_i P_{\sim \rho }  u_j +	 P_{ \sim \rho } u_i P_{\ll \rho} u_j \\
		&\,+P_\rho \big(\sum_{\nu^\prime \sim \nu  \gtrsim \rho}  	P_{ \nu^{ \prime}} u_i P_{\nu}u_j\big).
	\end{aligned}
 \end{equation}
The low--high terms and high--low terms has been estimated by \eqref{pre.bi1}, while
  the  high--high terms  follows directly from \eqref{pre.IntM1}.

	\end{proof}

We next estimate a delicate low--high interaction arising from the Bony decomposition of the nonlinear terms in the energy estimate. The argument exploits the special Gauss structure of \(\lambda\wedge\lambda\). Its difference analogue will   play an essential role in the weak difference estimates in dimension five; see the analysis of \eqref{uni.critic.d=5} in Section~\ref{sec.wpd.5}.

 \begin{lemma}  \label{lem_J11t4}
Let \(d\geq5\) and \(y_1=x_1\) in \(y=(y_1,y_2,y_3)\). Then
\begin{align} 
 		& \sum_{\mu \ll \rho}  \rho \mu^{-1} \sup_\tau \big\| \| P_{\mu} \left(\lambda\wedge\lambda\right) \|_{L^\infty_z} \|     \psi_{\sim \rho ,1}\|_{L^2_z}   \|     \psi_{\rho,1}^\tau\|_{L^2_z}     \big\|_{L^1_{t,y}}  \nonumber\\
 &\lesssim  \varepsilon^{\frac32} \rho^{-2s -2\sigma}  c_{\rho }^2 (\sigma),\label{est_J11t4}
 	\end{align}
    where \(\lambda\wedge\lambda\) is defined in \eqref{wed.lam}.
 \end{lemma}
 \begin{proof}
Applying Bony's decomposition to \(\lambda\wedge\lambda\), we obtain
\begin{align}
P_\mu(\lambda\wedge\lambda)
=&\,
L\bigl(P_{\ll\mu}\lambda,P_\mu\lambda\bigr)
+
L\bigl(P_\mu\lambda,P_{\ll\mu}\lambda\bigr)  
\nonumber\\
&+ \sum_{\mu\lesssim \nu\sim\nu'\ll\rho} P_\mu L(P_\nu\lambda,P_{\nu'}\lambda)+
\sum_{\nu\sim\nu'\gtrsim\rho}
P_\mu\bigl(\lambda_\nu\wedge\lambda_{\nu'}\bigr).
\label{gauss.bony}
\end{align}

Plugging the first two terms into \eqref{est_J11t4}, we can see it has been bounded by \eqref{pre.bi2}. The third term can be estimated by using \eqref{sup2} in a way similar  to those for \eqref{pre.bi2} as follows:
\begin{align*}
&\sum_{\mu\ll\rho}
\sum_{\mu\lesssim\nu\sim\nu'\ll\rho}
\rho\mu^{-1}
\sup_\tau
\left\|
\|P_\mu L(P_\nu\lambda,P_{\nu'}\lambda)\|_{L_z^\infty}
\|\psi_{\sim\rho,1}\|_{L_z^2}
\|\psi_{\rho,1}^\tau\|_{L_z^2}
\right\|_{L_{t,y}^1}
\\
\lesssim&\, \sum_{\mu\ll\rho}
\sum_{\mu\lesssim\nu\sim\nu'\ll\rho}
\rho\mu^{d-2}
\sup_\tau
\left\|
\|  L(P_\nu\lambda,P_{\nu'}\lambda)\|_{L_{\widehat{x}_1}^1}
\|\psi_{\sim\rho,1}\|_{L_{\widehat{x}_1}^2}
\|\psi_{\rho,1}^\tau\|_{L_{\widehat{x}_1}^2}
\right\|_{L_{t,x_1}^1}\\
\lesssim &\,
\varepsilon^{\frac32}
\rho^{-2s-2\sigma}
c_\rho^2(\sigma).
\end{align*}

As for the fourth term in \eqref{gauss.bony}, we  need to exploit the additional structure  given by \eqref{wed.lam}. For each dyadic frequency \(\nu\), write
\begin{align}
P_\nu
=
\nu^{-1}
\sum_{\beta=1}^d
\partial_\beta\widetilde P_{\nu,\beta},
\label{Pnu.derivative}
\end{align}
where the multipliers \(\widetilde P_{\nu,\beta}\) are uniformly bounded.
By \(u_3=\partial_xh\), the Codazzi equations \eqref{eq:codazzi-complex} imply schematically
\begin{align}
\partial_\beta\lambda_{\gamma\sigma}
&=
\partial_\sigma\lambda_{\gamma\beta}
+
L\bigl((u_3,A),\lambda\bigr),
\label{codazzi.one}
\\
\partial_\alpha\lambda^\alpha_\sigma
&=
\partial_\sigma\psi
+
L\bigl((u_3,A),\lambda\bigr),
\label{codazzi.contracted}
\end{align}
where we write \( (u_3,A)=u_3+A\) for convenience.

For \(\nu\sim\nu'\gtrsim\rho\), we obtain
\begin{align*}
P_\mu(\lambda_\nu\wedge\lambda_{\nu'})
=&\,
\nu^{-1}P_\mu
\left[
\partial_\sigma
\bigl(\widetilde P_{\nu,\beta}\lambda_{\gamma\beta}\bigr)
\overline{\psi_{\nu'}}
-
\partial_\alpha
\bigl(\widetilde P_{\nu,\beta}\lambda_{\beta\gamma}\bigr)
(P_{\nu^\prime}\overline{\lambda^\alpha _{\sigma}})
\right]
\\
&+
\nu^{-1}P_\mu
L\left(
P_\nu[\lambda\cdot(u_3,A)],
P_{\nu'}\lambda
\right).
\end{align*}
Applying the product rule and using
\(\lambda_{\gamma\beta}=\lambda_{\beta\gamma}\), we further get
\begin{align*}
P_\mu(\lambda_\nu\wedge\lambda_{\nu'})
=&\,
\frac{\mu}{\nu}
L(P_\nu\lambda,P_{\nu'}\lambda)
+
\nu^{-1}P_\mu
\left[
\widetilde P_{\nu,\beta}\lambda_{\beta\gamma}
\left(
\partial_\alpha P_{\nu^\prime}\overline{\lambda^\alpha _{\sigma}}
-
\partial_\sigma\overline{\psi_{\nu'}}
\right)
\right]
\\
&+
\nu^{-1}P_\mu
L\left(
P_\nu[\lambda\cdot(u_3,A)],
P_{\nu'}\lambda
\right),
\end{align*}
which combining the contracted Codazzi equation \eqref{codazzi.contracted}
yields
\begin{align}
P_\mu(\lambda_\nu\wedge\lambda_{\nu'})
=&\,
\frac{\mu}{\nu}
L(P_\nu\lambda,P_{\nu'}\lambda)+
\nu^{-1}P_\mu
L\left(
P_\nu[\lambda\cdot(u_3,A)],
\lambda_{\nu'}
\right)
\label{gauss.high.high.identity}
\\
&
 +
\nu^{-1}P_\mu
L\left(
P_{\nu'}[\lambda\cdot(u_3,A)],
\lambda_\nu
\right).
\nonumber
\end{align}
Thus the quadratic high--high interaction gains the factor
\(\mu/\nu\), while the remaining terms are cubic.

Putting  the first term in right hand side of in
\eqref{gauss.high.high.identity} into \eqref{est_J11t4}, Bernstein's inequality and \eqref{sup3} yield
\begin{align*}
 &\sum_{\mu\ll\rho}
\sum_{\nu\sim\nu'\gtrsim\rho}
\rho\mu^{-1}
\sup_\tau
\left\|
\| \frac{\mu}{\nu} \cdot P_\mu L(P_\nu\lambda,P_{\nu'}\lambda)\|_{L_z^\infty}
\|\psi_{\sim\rho,1}\|_{L_z^2}
\|\psi_{\rho,1}^\tau\|_{L_z^2}
\right\|_{L_{t,y}^1}
\\
& \lesssim
\sum_{\mu\ll\rho}
\sum_{\nu\sim\nu'\gtrsim\rho}
\rho\mu^{d-3}\nu^{-1}
\|P_\nu\lambda\|_{L_{t,y}^4L_z^2} \cdot
\|P_{\nu'}\lambda\|_{L_{t,y}^4L_z^2}
 \cdot
\|\psi_{\sim\rho,1}\|_{L_{t,y}^4L_z^2}\cdot
\|\psi_{\rho,1}\|_{L_{t,y}^4L_z^2}
\\
& \lesssim
\varepsilon^{\frac32}
\rho^{-2s-2\sigma}
c_\rho^2(\sigma)
\sum_{\mu\ll\rho}
\left(\frac{\mu}{\rho}\right)^{d-3}\\
&\lesssim  \varepsilon^{\frac32}
\rho^{-2s-2\sigma}
c_\rho^2(\sigma),
\end{align*}
where \(d\geq5\),  and thus  the last sum is bounded.

Finally, the cubic terms we only need to consider \(L\left(
P_{\nu'}(\lambda  u_3),
\lambda_{\nu} \right)\) as follows, as the rest is the same.
\begin{align}
& \sum_{\mu\ll\rho}\sum_{ \nu\sim\nu'\gtrsim\rho }
\rho (\mu \nu)^{-1}
\Big\|
\|P_\mu
L\left(
P_{\nu'}(\lambda  u_3),
\lambda_{\nu}
\right)\|_{L_z^\infty}   
\|\psi_{\sim\rho,1}\|_{L_z^2}
\|\psi_{\rho,1} \|_{L_z^2}
\Big\|_{L_{t,y}^1} \label{cubic.high.high}\\
%&\lesssim	 \sum_{\mu\ll\rho}\sum_{ \nu\sim\nu'\gtrsim\rho }
%\rho\mu^{-1}\nu^{-1}\| P_\muL\left(P_{\nu'}(\lambda\cdot u_3),P_{\nu}\lambda\right)\|_{L^2_{t,y } L^\infty_z}  \cdot\|\psi_{\sim\rho,1}\|_{L_{t,y}^4L_z^2}\cdot\|\psi_{\rho,1}\|_{L_{t,y}^4L_z^2}\nonumber\\
&\lesssim	 \sum_{\mu\ll\rho}\sum_{ \nu\sim\nu'\gtrsim\rho }
\rho\mu^{d-\frac 52}\nu^{-1}\| P_\mu
L\left(
P_{\nu'}(\lambda\cdot u_3),
P_{\nu}\lambda
\right)\|_{L^2_{t } L^1_x}  \cdot
\|\psi_{\sim\rho,1}\|_{L_{t,y}^4L_z^2}\cdot
\|\psi_{\rho,1}\|_{L_{t,y}^4L_z^2} \nonumber\\
&\lesssim  	 \sum_{\mu\ll\rho}\sum_{ \nu\sim\nu'\gtrsim\rho }
\rho\mu^{d-\frac 52}\nu^{-1}\left\|P_{\nu^\prime } ( \lambda \cdot u_3) \right\|_{L^2_{t,x}}  \left\| P_{\nu'} \lambda   \right\|_{L^\infty_{t } L^2_x} \cdot
\|\psi_{\sim\rho,1}\|_{L_{t,y}^4L_z^2}\cdot
\|\psi_{\rho,1}\|_{L_{t,y}^4L_z^2}  \nonumber\\
&\lesssim   \sum_{\mu\ll\rho}\sum_{ \nu\sim\nu'\gtrsim\rho }
\rho\mu^{d-\frac 52}\nu^{-1} \varepsilon^{\frac34}(\nu')^{-s} c_{\nu'}(0)\cdot \varepsilon^{-\frac 18}\nu^{-s}c_{\nu}(0) \cdot  \varepsilon^{-\frac 14} \rho^{-2s+\frac12-2\sigma} c^{2}_{\rho} (\sigma) \nonumber\\
& \lesssim \varepsilon^{2}
\rho^{-2s-2\sigma}
c_\rho^2(\sigma). \nonumber
\end{align}

This completes the proof of this Lemma.

 \end{proof}

An analogue of Lemma~\ref{lem_J11t4} also holds for cubic terms.	
     \begin{lemma}  \label{lem_cubic.term}
 	Let \(y_1=x_1\) in \(y=(y_1,y_2,y_3)\), $ i,j,l\in\{1,2,3,4,5\} $ and \(d \geq 4\). Then we have
	\begin{align}
		\sum_{\mu\ll\rho}	 	\sup_{\tau, \tau^\prime } &\,\rho \mu^{-2}\left\| \left\| 	P_{\mu} \left(u_i u_j u_l \right)\right\|_{L^\infty_z} \left\|     \psi_{\sim\rho ,1}^{\tau } \right\|_{L^2_z}    \|   \psi_{\rho,1} ^{\tau^\prime} \|_{L^2_z}     \right\|_{L^1_{t,y}}\nonumber \\
			 \label{Rem_estJ11t4} & \lesssim   \varepsilon^{\frac 32} \rho^{-2s -2\sigma}  c_{\rho }^2 (\sigma).
	\end{align}
 \end{lemma}
\begin{proof}
   Apply Bony's decomposition to the product \((u_i u_j)u_l\):
 \begin{align}  \label{cubic.bony}
     P_{\mu} \left(u_i u_j u_l \right)  
=&\, P_{\ll\mu}   (u_i   u_j) P_{    \mu} u_l+    P_{\mu} (u_i   u_j) P_{   \ll \mu} u_l\\
&\,+\sum_{ \mu \lesssim \nu^\prime \sim \nu  \ll \rho}P_\mu \big(P_{\nu^\prime } (u_i u_j )P_{\nu}   u_l\big)+\sum_{ \nu^\prime \sim \nu \gtrsim\rho}P_\mu \big(P_{\nu^\prime } (u_i u_j )P_{\nu}   u_l\big).\nonumber
 \end{align}
The low--high term and high--low term can be estimated  similarly, and we will only estimate the high low.  Putting the   low--high term into \eqref{lem_cubic.term}, we can use \eqref{pre.bi1} and \eqref{pre.bi3} to obtain
\begin{align*}
  &  \sum_{\mu\ll\rho}	 	\sup_{\tau, \tau^\prime }  \rho \mu^{-2}\left\| \left\| 	 P_{\ll\mu}   (u_i   u_j) P_{    \mu} u_l\right\|_{L^\infty_z} \left\|     \psi_{\sim\rho ,1}^{\tau } \right\|_{L^2_z}    \|   \psi_{\rho,1} ^{\tau^\prime} \|_{L^2_z}     \right\|_{L^1_{t,y}}\nonumber \\
  \lesssim&\, \rho  \sum_{\omega \ll\mu\ll\rho}	 	\sup_{\tau  }    \mu^{-\frac 32}\big\|  \| 	 P_{\omega}   (u_i   u_j)   \|_{L^\infty_{\widehat{x}_1}}  \|     \psi_{\sim\rho ,1}^{\tau }  \|_{L^2_{\widehat{x}_1}}        \big\|_{L^1_{t,x_1}} \\
  &\cdot \sum_{\mu\ll\rho}	 	\sup_{ \tau^\prime }   \mu^{-\frac 12}\left\|  \| 	P_{\mu}  u_l   \|_{L^\infty_{ \widehat{x}_1}}       \|   \psi_{\rho,1} ^{\tau^\prime} \|_{L^2_{\widehat{x}_1}}     \right\|_{L^1_{t,x_1}}\nonumber \\
  \lesssim &\, \varepsilon^{\frac 32} \rho^{-2s -2\sigma}  c_{\rho }^2 (\sigma).
\end{align*}
Similarly for the third term in \eqref{cubic.bony}, we can put it into  \eqref{lem_cubic.term},  estimate it by using \eqref{sup2} in a way similar  to those for  \eqref{pre.bi1} and \eqref{pre.bi3}  as follows 
\begin{align*}
      &  \sum_{ \mu \lesssim \nu^\prime \sim \nu  \ll \rho} 	\sup_{\tau, \tau^\prime }  \rho \mu^{-2}\left\| \left\| 	 P_\mu \big(P_{\nu^\prime } (u_i u_j )P_{\nu}   u_l\big)\right\|_{L^\infty_z} \left\|     \psi_{\sim\rho ,1}^{\tau } \right\|_{L^2_z}    \|   \psi_{\rho,1} ^{\tau^\prime} \|_{L^2_z}     \right\|_{L^1_{t,y}}\nonumber \\
      \lesssim &\, \sum_{ \mu \lesssim \nu^\prime \sim \nu  \ll \rho} 	\sup_{\tau, \tau^\prime }  \rho \mu^{d-3}\left\| \big\| 	 P_\mu \big(P_{\nu^\prime } (u_i u_j )P_{\nu}   u_l\big)\right\|_{L^1_{\widehat{x}_1}} \left\|     \psi_{\sim\rho ,1}^{\tau } \right\|_{L^2_{\widehat{x}_1}}    \|   \psi_{\rho,1} ^{\tau^\prime} \|_{L^2_{\widehat{x}_1}}     \big\|_{L^1_{t,x_1}} \\
      \lesssim &\,\varepsilon^{\frac 32} \rho^{-2s -2\sigma}  c_{\rho }^2 (\sigma).
\end{align*}

For high high term   for the cubic terms after Bony decomposition, we can deal it in the same manner as \eqref{cubic.high.high}, by noting that \(\sum_{\mu\ll\rho} \mu^{d-\frac 72} \lesssim \rho^{d-\frac 72} \), if \(d \geq 4\).

Combining the four
contributions proves \eqref{Rem_estJ11t4}.

\end{proof}

\begin{remark}  \label{Rem_J11t4}
	In our framework, multiplying the nonlinear term by $ u_j $ has the same effect as taking one derivative. Taking the nonlinear term $u_i u_j u_l$ for example, we have the decomposition
	\begin{align*}
		u_i u_j u_l&= \sum_{\mu_1, \mu_2, \mu_3} P_{\mu_1} u_i P_{\mu_2} u_j P_{\mu_3} u_l.
	\end{align*}
Without loss of generality, we assume $ \mu_1\leq \mu_2 \leq \mu_3, $ and define $\|f\|_{X}$ as some norm of $f$, then schematically, it holds
	\begin{align*}
		\left\| P_{\mu_1} u_i P_{\mu_2} u_j P_{\mu_3} u_l\right\|_{X} \lesssim &\,\| P_{\mu_1} u_i  \|_{L^\infty_x} \left\|P_{\mu_2} u_j P_{\mu_3} u_l\right\|_{X}  \\
		\lesssim  &\,\| P_{\mu_1} u_i  \|_{\dot{B}^{s+1}_{2,1}} 	\left\|P_{\mu_2} u_j P_{\mu_3} u_l\right\|_{X}\\
		\lesssim  &\,\varepsilon^{\frac 78}\mu_2 \left\|P_{\mu_2} u_j P_{\mu_3} u_l\right\|_{X},
	\end{align*}
where we used the fact $ \| P_{\mu_1} u_i  \|_{\dot{B}^{s}_{2,1}} 	\lesssim  \varepsilon^{\frac 78}.$

Thus cubic terms are more favorable than differentiated quadratic terms.  We shall therefore suppress cubic
contributions in later arguments whenever they are covered by the preceding
principle and focus on the more delicate quadratic interactions.

\end{remark}

The following two lemmas are devoted to estimating the typical nonlinear terms involving the combination of the Schr\"odinger solution $ \psi $ and elliptic solution $ u_{\mathfrak{T}}$ for  $\mathfrak{T}\in\{3,4,5\}.  $

The main obstacle in Lemma~\ref{lem_J1} is that the low frequency factor
\(P_{\ll\rho}\mathbf{U}\) is not directly controlled by the bootstrap bounds.
We therefore decompose it dyadically and use the elliptic equation at each
frequency, thereby converting the quadratic interaction into cubic terms with
additional elliptic gain.
\begin{lemma} \label{lem_J1}
For  $ \mathbf{U}=(W, u_3) $ and \(W=V-2A\),  we have
	\begin{align}
	 \left|\int_0^T \int_{\mathbb{R}^{d}} L( P_{\ll \rho } \mathbf{U},   \psi_{\sim \rho }) \cdot\partial_x  \overline{\psi}_{\rho} \mathrm{ d} x\mathrm{ d} t\right|
 \lesssim \varepsilon^{\frac32}\rho ^{-2s-2\sigma}c^2_\rho(\sigma).  \label{est_J1}
	\end{align}
\end{lemma}

\begin{proof} 	 By decomposition \eqref{decomp_Xi},  we have
\begin{align}
       	& \left|\int_0^T \int_{\mathbb{R}^{d}} L( P_{\ll \rho } \mathbf{U},   \psi_{\sim \rho}) \cdot\partial_x  \overline{\psi}_{\rho} \mathrm{ d} x\mathrm{ d} t\right| \label{lem_J1.decom}\\
            \leq &	  \sum_{\gamma,\gamma^\prime=1}^d\left|\int_0^T \int_{\mathbb{R}^{d}} L( P_{\ll \rho } \mathbf{U},   \psi_{\sim \rho ,\gamma}) \cdot\partial_x  \overline{\psi}_{\rho,\gamma^\prime} \mathrm{ d} x\mathrm{ d} t\right|. \nonumber
\end{align}
By symmetry, it suffices to consider the case $\gamma=\gamma'=1$.

Recall equation  \eqref{eq_U} for  $ \mathbf{U} $, namely
\begin{align*}
    \Delta \mathbf{U} =h\cdot\partial_x^2 \mathbf{U}
+
\partial_xh\cdot\partial_xu_{\mathfrak T}+
\partial_x(\lambda\wedge\lambda)
+
u_i u_j u_l,
\end{align*}
where  \(u_3=\partial_x h,  \mathfrak{T}\in\{3,4,5\}\) and $ i,j,l\in\{1,2,3,4,5\} $.
Applying $\Delta^{-1}P_\mu$ to the above equation and using the Fourier transform, Bony's decomposition, and the formula \eqref{L_formula_1}--\eqref{L_formula_3} for    $L$ notation, we obtain
\begin{equation}\label{Bony_U}
    \begin{aligned}  
	P_\mu     \mathbf{U} 
= &\,   L(P_{\leq  \mu} h, P_\mu  \mathbf{U} )+ \mu^{-1} L(P_\mu u_3, P_{\ll \mu } u_{\mathfrak{T}}) \\
&\,+\mu^{-2}\sum_{ \rho\gg\nu\sim \nu^\prime \gtrsim \mu} P_\mu L (P_\nu u_3, \partial_x P_{\nu^\prime}  u_{\mathfrak{T}} ) \\
&\, + \mu^{-2}\sum_{ \nu\sim \nu^\prime \gtrsim \rho}  P_\mu L (P_\nu u_3, \partial_x P_{\nu^\prime}  u_{\mathfrak{T}}) \\
 &\, +\mu^{-1}  P_\mu (\lambda\wedge\lambda)+\mu^{-2}P_\mu (u_iu_ju_l). 
\end{aligned}
\end{equation}
Consequently,
\begin{align}
& \left|\int_0^T \int_{\mathbb{R}^{d}} L( P_{\ll \rho } \mathbf{U},   \psi_{\sim \rho,1}) \cdot\partial_x  \overline{\psi}_{\rho,1} \mathrm{ d} x\mathrm{ d} t\right|\nonumber\\
    		\lesssim&\,   	\sum_{\mu \ll \rho}    \rho       \sup_{\tau} \big\|\left\|P_{\mu} \mathbf{U}  \right\|_{L^\infty_z} \|  \psi_{\sim \rho,1}  \|_{L^2_z}   \|    \psi_{\rho,1}^\tau \|_{L^2_z}   \big\|_{L^1_{t,y}}\label{J11t0}\\
		\lesssim&\,   	\sum_{\mu \ll \rho}   \rho  \|P_{\leq \mu} h\|_{L^\infty_x}  \cdot  \sup_{\tau} \big\|\left\|P_{\mu} \mathbf{U}  \right\|_{L^\infty_z} \|  \psi_{\sim \rho,1}  \|_{L^2_z}   \|    \psi_{\rho,1}^\tau \|_{L^2_z}   \big\|_{L^1_{t,y}}\label{J11t1}\\
	&\, + 	\sum_{\mu \ll \rho} \rho \mu^{-1} \sup_\tau \big\|\left\|L(P_\mu u_3, P_{\ll \mu } u_{\mathfrak{T}}) \right\|_{L^\infty_z} \|    \psi_{\sim \rho,1} \|_{L^2_z}   \|    \psi_{\rho,1}^\tau \|_{L^2_z}   \big\|_{L^1_{t,y}} \label{J11t2}\\
  &\, +	\sum_{\mu \ll \rho} \sum_{\rho \gg \nu\sim \nu^\prime \gtrsim \mu}   \rho \mu^{-2} \sup_{\tau}\big\|\left\| P_\mu L (P_\nu u_3, P_{\nu^\prime} \partial_x u_{\mathfrak{T}} ) \right\|_{L^\infty_z}\label{J11t3_0}\\
 &\quad\quad\cdot\|   \psi_{\sim \rho,1} \|_{L^2_z}   \| \psi_{\rho,1}^\tau\|_{L^2_z}   \big\|_{L^1_{t,y}}\nonumber\\
 &\, +	\sum_{\mu \ll \rho} \sum_{\nu\sim \nu^\prime \gtrsim \rho} \rho \mu^{-2}  \sup_{\tau}\big\|\left\| P_\mu  L (P_\nu u_3, P_{\nu^\prime} \partial_x u_{\mathfrak{T}} ) \right\|_{L^\infty_z}\label{J11t3}\\
 &\quad\quad\cdot\left\|  \psi_{\sim \rho,1} \right\|_{L^2_z}   \left\|   \psi_{\rho,1}^\tau \right\|_{L^2_z}   \big\|_{L^1_{t,y}}\nonumber\\
	&\,  +	\sum_{\mu \ll \rho}  \rho  \mu^{-1} 
 \sup_{\tau} \big\|\left\|P_{\mu}(\lambda\wedge\lambda)\right\|_{L^\infty_z} \cdot\left\|   \psi _{\sim \rho}\right\|_{L^2_z}   \left\| \psi_{\rho,1}^\tau\right\|_{L^2_z}   \big\|_{L^1_{t,y}} \label{J11t4}\\
	&\, + 	\sum_{\mu \ll \rho}  \rho \mu^{-2} \sup_{\tau} \big\|\left\|P_{\mu}(u_i u_j u_k)\right\|_{L^\infty_z} \cdot\|   \psi_{\sim \rho,1}\|_{L^2_z}   \|    \psi_{\rho,1}^\tau\|_{L^2_z}   \big\|_{L^1_{t,y}}, \label{J11t5}
\end{align}
where we let \(y_1=x_1\) in \(y=(y_1,y_2,y_3)\).
It is obvious that  \eqref{J11t1} can be absorbed  by  \eqref{J11t0}, since by \eqref{h_infty} we have $ \|P_{\leq \mu} h\|_{L^\infty_x}\lesssim \varepsilon^{\frac 78} $. 
\eqref{J11t4} has been proved  in Lemma \ref{lem_J11t4} and \eqref{J11t2} is a special case of \eqref{pre.bi2},  while  \eqref{J11t5} can been handled in the way as in Remark \ref{Rem_J11t4}.
Hence we are left with the estimates of \eqref{J11t2}, \eqref{J11t3_0} and \eqref{J11t3}.

Observing that the decomposition \eqref{decomp_Xi} and  Bernstein's inequality give
\begin{align*}
    \eqref{J11t3_0}  
 \lesssim &\, 	\sum_{\mu \ll \rho } \sum_{\rho \gg \nu\sim \nu^\prime \gtrsim \mu}   \rho \mu^{d-3} \sup_{\tau } \big\|\left\| L (P_\nu u_3, P_{\nu^\prime} \partial_x u_{\mathfrak{T}} ) \right\|_{L^1_{\widehat{x}_1}}\nonumber\\
 &\cdot \|     \psi_{\sim \rho ,1}  \|_{L^2_{\widehat{x}_1}}   \|  \psi_{\rho,1}^{\tau } \|_{L^2_{\widehat{x}_1}}   \big\|_{L^1_{t,x_1}}, 
\end{align*}
we can see  \eqref{J11t3_0} can be estimated just as \eqref{J11t2}.

For \eqref{J11t3},  we can use \eqref{sup1}, \eqref{supV} and \eqref{sup3} to obtain
\begin{align*}
\eqref{J11t3} \lesssim &\,	\sum_{\mu \ll \rho} \sum_{\nu\sim \nu^\prime \gtrsim \rho }      \rho  \mu^{-2}  \left\| P_\mu (L (P_\nu u_3, P_{\nu^\prime} \partial_x u_{\mathfrak{T}} )) \right\|_{L^1_{t}L^\infty_{x}} \cdot\big \|  \| \psi_{\sim\rho }  \|_{L^2_z}   \left\|   \psi_{\rho} \right\|_{L^2_z}   \big\|_{L^\infty_{t} L^1_{y}}\nonumber\\
 \lesssim &\, 	\sum_{\mu \ll \rho} \sum_{\nu\sim \nu^\prime \gtrsim \rho }    \rho  \mu^{d-2} \nu^\prime  \left\|  P_\nu   u_3 \right\|_{L^2_{t,x}} \left\|  P_{\nu^\prime}   u_{\mathfrak{T}} \right\|_{L^2_{t,x}}  \cdot \left\|  \psi_{\sim \rho } \right\|_{L^\infty_t L^2_x}   \left\|    \psi_{\rho} \right\|_{L^\infty_t L^2_x}    \nonumber\\
 \lesssim &\, \varepsilon^{-\frac{1}{2}} 	\sum_{\mu \ll \rho} \sum_{\nu\sim \nu^\prime \gtrsim \rho } \sum_{\rho^\prime \sim \rho} \rho  \mu^{d-2} \nu^\prime \nu^{-s-1} c_{\nu}(0) (\nu^\prime)^{-s-1} c_{\nu^\prime}(0)\\
 &\cdot (\rho^\prime)^{-s-\sigma} c_{\rho^\prime} (\sigma)  \rho^{-s-\sigma} c_\rho (\sigma)\\
 \lesssim &\,\varepsilon^{\frac 32} \rho^{-2s-2\sigma} c_\rho^2 (\sigma).
\end{align*}
This completes the proof.
	\end{proof}

    Different from Lemma \ref{lem_J1}, the estimates in Lemma \ref{lem_J2} can be controlled directly by   
bootstrap assumption of \(u_{\mathfrak{T}}\).
\begin{lemma} \label{lem_J2}
For $\mathfrak{T}\in\{3,4,5\}  $, we have
	\begin{align}
		   \| L( P_{\sim \rho}  u_{\mathfrak{T}} ,  \partial_x   \psi_{\ll \rho } ) \cdot  \overline{\psi}_{\rho}   \|_{L^1_{t,x}}& \lesssim \varepsilon^{\frac58} \rho ^{-2s-2\sigma}c^2_\rho(\sigma),  \label{est_J2} \\
				\sum_{\nu^{\prime} \sim \nu \gtrsim \rho }  \| L(P_{\nu}  u_{\mathfrak{T}},   \psi_{\nu^\prime } ) \cdot \partial_x  \overline{\psi}_{\rho}   \|_{L^1_{t,x}} &\lesssim \varepsilon^{\frac58}\rho ^{-2s-2\sigma}c^2_\rho(\sigma).  \label{est_J3}
	\end{align}
\end{lemma}
\begin{proof}  The estimates  \eqref{est_J2} and  \eqref{est_J3}  follows from \eqref{pre.IntM2} and \eqref{pre.IntM.infty} respectively.  To be specific, we can use  \eqref{supV}, \eqref{sup3}, \eqref{pre.IntM2} and H\"older's inequality   to obtain
	\begin{align*}
		&\,    \left\| L( P_{\sim \rho}  u_{\mathfrak{T}} , \partial_x     \psi_{\ll \rho } ) \cdot   \overline{\psi}_{\rho}  \right\|_{L^1_{t,x}}\\
		\lesssim&\,      \left\| P_{\sim \rho}  u_{\mathfrak{T}}  \right\|_{L^2_{t,x}}\left\|    P_{\ll \rho } \partial_x  \psi \right\|_{L^4_{t,y}L^\infty_z} \left\|    \psi_{\rho}  \right\|_{L^4_{t,y}L^2_z} \\
		\lesssim &\,\varepsilon^{-\frac 14}  \sum_{\rho^\prime \sim  \rho}   (\rho^\prime)^{-s-1-\sigma} c_{\rho^\prime} (\sigma)\cdot \varepsilon^{\frac78}   \rho^{ \frac 34}\cdot\rho^{-s+\frac 14-\sigma} c_\rho (\sigma)\\
		\lesssim &\, \varepsilon^{\frac58}\rho ^{-2s-2\sigma}c^2_\rho(\sigma)
	\end{align*}
and by \eqref{pre.IntM.infty} and the   slow variation property of the $\sigma$-frequency envelope \eqref{ini3}
\begin{align*}
&\sum_{\nu^{\prime} \sim \nu \gtrsim \rho } \left\| L(P_{\nu}  u_{\mathfrak{T}}, \psi_{\nu^\prime } ) \cdot \partial_x  \overline{\psi}_{\rho}  \right\|_{L^1_{t,x}} \\
	 \lesssim&\, \sum_{\nu^{\prime} \sim \nu \gtrsim \rho } \rho \left\| P_{\nu}  u_{\mathfrak{T}}  \right\|_{L^2_{t,x}}\left\|   \psi_{\nu^\prime }  \right\|_{L^4_{t,y}L^2_z}  \|   \psi_{\rho}\|_{L^4_{t,y}L^\infty_z}\\
\lesssim &\,\varepsilon^{-\frac 14} \sum_{\nu^{\prime} \sim \nu \gtrsim \rho } \rho   \nu^{-s-1-\sigma} c_\nu (\sigma) \cdot(\nu^\prime)^{-s+\frac 14} c_{\nu^\prime}(0)  \cdot  \varepsilon^{-\frac18}  \rho^{-\frac 14} c_\rho (\sigma)\\
\lesssim &\, \varepsilon^{\frac58}\rho ^{-2s-2\sigma}c^2_\rho(\sigma).
\end{align*} 
	\end{proof}

We next establish an $L_x^2$ estimate for the typical term $ u_i u_j$, which will be used in the elliptic estimates.
\begin{lemma} \label{lem_l2space_u2}
	For $ i,j =1,2,3,4,5 $, we have
	\begin{align}
\left\| P_\rho (u_i  u_{j} )\right\|_{L^2_x}\lesssim  \varepsilon^{\frac 34}    \rho^{-s+1 -\sigma} c_\rho (\sigma). \label{est_L2nlt_type1}
	\end{align}
\end{lemma}
\begin{proof}
	By the Bony's decomposition in  \eqref{Bony2},  it suffices to estimate the low--high and high--high contributions, since   \(u_i u_j\) is symmetry.
	By H\"older's inequality,  Bernstein's inequality \eqref{h_infty} and \eqref{sup1}, we see that
	\begin{align*}  
		\| P_{\ll \rho} u_i P_\rho   u_{j}\|_{L^2_x}& 
		\lesssim \|P_{\ll \rho} u_i\|_{L^\infty_x} \cdot  \| P_\rho u_{j}\|_{L^2_x}\nonumber\\
		&\lesssim \varepsilon^{\frac 78} \cdot \rho  \varepsilon^{-\frac 18} \rho^{-s -\sigma} c_\rho (\sigma)\nonumber \\
		&\lesssim  \varepsilon^{\frac 34}    \rho^{-s + 1-\sigma} c_\rho (\sigma) \nonumber
	\end{align*}
	and  
	\begin{align*} 
		&\Big\| P_\rho\big(\sum_{ \nu\sim \nu^{\prime} \gtrsim \rho} P_\nu u_i   P_{\nu^\prime}   u_{j}\big)\Big\|_{L^2_x} \\
		 \lesssim &\sum_{ \nu\sim \nu^{\prime} \gtrsim \rho}\rho^{\frac d2} \| P_\nu u_i   P_{\nu^\prime}   u_{j} \|_{L^1_x} \nonumber\\
		\lesssim &\sum_{ \nu\sim \nu^{\prime} \gtrsim \rho}\rho^{\frac d2} \| P_\nu u_i \|_{L^2_x} \|   P_{\nu^\prime}   u_{j} \|_{L^2_x} \nonumber \\
		\lesssim &\,\varepsilon^{-\frac 14}  \sum_{ \nu\sim \nu^{\prime} \gtrsim \rho}\rho^{\frac d2} \nu^{-s } c_\nu (0)  \cdot   (\nu^\prime)^{-s -\sigma} c_{\nu^\prime} (\sigma) \nonumber\\
		\lesssim & \,\varepsilon^{\frac 34}\   \rho^{-s +1-\sigma} c_\rho (\sigma),\nonumber
	\end{align*} 
    where we have used \eqref{ini3}.
Combining the preceding bounds proves \eqref{est_L2nlt_type1}.
	\end{proof}

The next lemma treats the nonlinear terms involving $B$.
\begin{lemma}
For \(i,j\in\{1,2,3,4,5\},\) we have
\begin{align}
     \rho^{-1} \int_0^T \int_{\mathbb{R}^{d}}    L(  \partial_xB_{\ll \rho }, P_{\sim\rho } u_i ) \cdot   P_{\rho} u_j  \,\mathrm{ d} x \mathrm{ d} t& \lesssim \varepsilon^{\frac58}\rho ^{-2s-2\sigma}c_{\rho }^2 (\sigma),
 \label{est_EnBt1}\\
	    \int_0^T \int_{\mathbb{R}^{d}}    L( B_{ \sim \rho }, P_{\ll \rho  } u_i )\cdot P_{\rho} u_j  \,\mathrm{ d} x \mathrm{ d} t &\lesssim\varepsilon^{\frac58}\rho ^{-2s-2\sigma}c_{\rho }^2 (\sigma),\label{est_EnBt2}\\
	 \sum_{\nu\sim \nu^\prime\gtrsim \rho }   \int_0^T \int_{\mathbb{R}^{d}}  P_\rho L(B_{ \nu},P_{\nu^\prime }u_i)\cdot P_{\rho} u_j \, \mathrm{ d} x \mathrm{ d} t &\lesssim \varepsilon^{\frac58}\rho ^{-2s-2\sigma}c_{\rho }^2 (\sigma). \label{est_EnBt3}
\end{align}
\end{lemma}

\begin{proof}    
The \eqref{est_EnBt1}  follows from
  H\"older's inequality,  \eqref{sup3} and \eqref{Btx2z_inf} directly:
\begin{align*}
& \rho^{-1} \int_0^T \int_{\mathbb{R}^{d}}    L(  \partial_xB_{\ll \rho }, P_{\sim\rho } u_i ) \cdot   P_{\rho} u_j  \,\mathrm{ d} x \mathrm{ d} t\\	
 \lesssim &\,    \rho^{-1} \left\| \partial_x B_{\ll \rho}\right\|_{L^2_{t,y }L^\infty_z } \left\|P_{\sim\rho }u_i \right\|_{L^4_{t,y } L^2_z }  \left\|P_{\rho} u_j \right\|_{L^4_{t,y } L^2_z }\nonumber \\
 \lesssim  &\, \varepsilon^{\frac58}\rho ^{-2s-2\sigma}c_{\rho }^2 (\sigma).  
\end{align*}	
Similarly for \eqref{est_EnBt3},  
\begin{align*}
&\sum_{\nu\sim \nu^\prime\gtrsim \rho }   \int_0^T \int_{\mathbb{R}^{d}}  P_\rho L(B_{ \nu},P_{\nu^\prime }u_i)\cdot P_{\rho} u_j \, \mathrm{ d} x \mathrm{ d} t\\
    \lesssim &\,   \sum_{\nu\sim \nu^\prime\gtrsim \rho }  \left\|B_{ \nu}\right\|_{L^2_{t,x} } \left\| P_{\nu^\prime}u_i \right\|_{L^4_{t,y } L^2_z }  \left\|P_{\rho}u_j \right\|_{ L^4_{t,y }L^\infty_z}\\
	\lesssim&\,\varepsilon^{-\frac{3}{8}}  \sum_{\nu\sim \nu^\prime\gtrsim \rho }     \nu^{-s} c_{\nu} (0) (\nu^\prime)^{-s+\frac 14-\sigma} c_{\nu^\prime} (\sigma) \rho^{-s+\frac 14-\sigma} c_\rho (\sigma)\\
	\lesssim&\, \varepsilon^{\frac58}\rho ^{-2s-2\sigma}c_{\rho }^2 (\sigma),
\end{align*}
where we have used \eqref{pre.IntM.infty} and \eqref{ini3}.
Turning to \eqref{est_EnBt2}, it  follows from the decomposition \eqref{decomp_Xi},  and \eqref{pre.bi1} that
\begin{align*}
&  \int_0^T \int_{\mathbb{R}^{d}}    L( B_{ \sim \rho }, P_{\ll \rho  } u_i )\cdot P_{\rho} u_j  \,\mathrm{ d} x \mathrm{ d} t \\	
 \lesssim &\,    \| P_{\sim \rho } B\|_{L^2_{t,x}} \cdot  \sum_{\gamma=1}^d  \left\| \|P_{\ll \rho  } u_i \|_{L^\infty_{\widehat{x}_\gamma }}  \|P_{\rho, \gamma}  u_j\|_{L^2_{\widehat{x}_\gamma}}\right\|_{L^2_{t,x_\gamma}}\\
	\lesssim&\, \varepsilon^{\frac58}\rho ^{-2s-2\sigma}c_{\rho }^2 (\sigma).
\end{align*}
This completes the proof of  the lemma.
\end{proof}

\section{Energy estimate}\label{secEn}

In this section, we are going to establish the energy estimate for the nonlinear Schr\"odinger equation \eqref{gaugsmcf} coupled with the elliptic system  \eqref{eq_ell} and \eqref{eq:B-elliptic}. Letting  $ Q_\rho = P_\rho $ in the equation \eqref{Q_psi} implies
\begin{align} \label{Plamdeq}
	&\sqrt{-1} \pa_t \psi_{\rho} +  \partial_{ \alpha} (g^{\alpha \beta}  \partial_{ \beta} \psi_{\rho} ) =P_{\rho} \mathcal{N}- \partial_{ \alpha}  ([P_{\rho}, h^{\alpha \beta}]  \partial_{ \beta}\psi   ):= \mathcal{I}_\rho,
\end{align}
 multiplying which with $ \overline{\psi}_{\rho} $, taking the imaginary part and then integrating over $ [0,T]\times \mathbb{R}^d $ yield
\begin{align}
	&\|\psi_\rho(t)\|^2_{L^2_x} - \|P_\rho\psi_0\|^2_{L^2_x}\nonumber\\
	=&\,  \operatorname{Im}\int_0^T \int_{\mathbb{R}^{d}}  \overline{\psi}_{\rho}\cdot  P_{\rho} \mathcal{N}    \mathrm{ d} x \mathrm{ d} t \label{enenltn}\\
	&\,- \operatorname{Im}\int_0^T   \int_{\mathbb{R}^{d}}    \partial_{\alpha}( [P_{\rho}, h^{\alpha\beta}] \partial_{ \beta} \psi)  \cdot \overline{\psi}_{\rho}    \mathrm{ d} x \mathrm{ d} t. \label{enenltc}
\end{align}
For the second nonlinear term \eqref{enenltc},  we can use  integration by parts and Bony's decomposition of commutator \eqref{Bony_com} to obtain
\begin{align}
	&- \mathrm{Im}\int_0^T \int_{\mathbb{R}^{d}}   \partial_{\alpha} ( [P_{\rho}, h^{\alpha\beta}] \partial_{ \beta} \psi  ) \cdot \overline{\psi}_{\rho} \,  \mathrm{ d} x\mathrm{ d} t\nonumber\\
	= &\,\mathrm{Im}\int_0^T \int_{\mathbb{R}^{d}}    ( [P_{\rho}, h^{\alpha \beta} ] \partial _{\beta} \psi ) \partial _{\alpha} \overline{\psi}_{\rho} \,  \mathrm{ d} x\mathrm{ d} t\nonumber\\
	=   &\,\int_0^T \int_{\mathbb{R}^{d}} L ( \partial_x P_{\ll \rho }  h^{\alpha \beta},   \psi_{\sim \rho })\cdot \partial_{\alpha} \overline{\psi}_{\rho} \,\mathrm{ d} x\mathrm{ d} t\label{ene_ht1}\\
	&\quad +    \int_0^T \int_{\mathbb{R}^{d}} L (\partial_{ \beta} P_{\ll \rho  }   \psi, P_{\sim \rho } h^{ \alpha \beta})\cdot \partial_{\alpha} \overline{\psi}_{\rho}\, \mathrm{ d} x\mathrm{ d} t\label{ene_ht2}\\
	&\quad +\sum_{\nu^{\prime} \sim \nu \gtrsim \rho }  \int_0^T \int_{\mathbb{R}^{d}} P_{\rho} L( P_{\nu} h^{\alpha \beta},  \partial _{ \beta}   
 \psi_{\nu^\prime })\cdot\partial _{\alpha} \overline{\psi}_{\rho} \,\mathrm{ d} x\mathrm{ d} t\label{ene_ht3}.
\end{align}
Given that $ h^{\alpha \beta} =g^{\alpha \beta} -I_d= (I_d+h_{\alpha\beta})^{-1}-I_d = -h_{\alpha\beta}+o((h_{\alpha\beta})^2)$ and $ h_{\alpha \beta} $ satisfies \eqref{ellis2}, then the estimate of 	\eqref{ene_ht1} is completely similar to that in Lemma \ref{lem_J1}.  Meanwhile, the estimates for \eqref{ene_ht2} and \eqref{ene_ht3} follow from Lemma \ref{lem_J2} for $u_{\mathfrak{T}}= \partial_x h.$

%\begin{align*}
%    &-\int_0^T \int_{\mathbb{R}^{d}} \mathrm{Im}\left( \partial_{\alpha} \left(\left[P_{\rho}, h^{\alpha \beta}\right] \partial _{ \beta} \psi \right)\overline{\psi}_{\rho} \right) \mathrm{ d} x\mathrm{ d} t\nonumber\\
	%=&\,\int_0^T \int_{\mathbb{R}^{d}}  \mathrm{Im}\left(\left(\left[P_{\rho}, h^{\alpha \beta}\right] \partial _{\beta} \psi \right) \partial _{\alpha} P_{\rho}\overline{\psi} \right) \mathrm{ d} x\mathrm{ d} t\nonumber\\
	%=&\,\mathrm{Im}\int_0^T \int_{\mathbb{R}^{d}}\sum_{\rho^\prime \sim  \rho} \left[P_{\rho}, P_{\ll \rho^\prime } h^{\alpha \beta}\right] \partial _{ \beta} P_{\rho^\prime} \psi \partial_{\alpha}P_{\rho}  \overline{\psi} \mathrm{ d} x\mathrm{ d} t \\
	%&\,+\mathrm{Im}\int_0^T \int_{\mathbb{R}^{d}}\sum_{\rho^\prime \sim  \rho} \left[P_{\rho}, P_{\rho^\prime} h^{ \alpha \beta}\right] P_{\ll \rho^\prime }  \partial_{ \beta} \psi \partial_{\alpha} \overline{\psi}_{\rho} \mathrm{ d} x\mathrm{ d} t \\
	%&\,+\mathrm{Im}\int_0^T \int_{\mathbb{R}^{d}} \sum_{\nu^{\prime} \sim \nu \gtrsim \rho }\left [P_{\rho}, P_{\nu} h^{\alpha \beta}\right] \partial _{ \beta}  P_{\nu^\prime } \psi\partial _{\alpha} \overline{\psi}_{\rho} \mathrm{ d} x\mathrm{ d} t.
%\end{align*}

For the first nonlinear term \eqref{enenltn}, it follows from \eqref{energy_B} that
\begin{align}
	\eqref{enenltn}
	=&\,  \mathrm{Im} \int_0^T \int_{\mathbb{R}^{d}}    P_{\rho} \left( \mathbf{U} \cdot\partial_x \psi\right) \overline{\psi}_{\rho}  \mathrm{ d} x \mathrm{ d} t\label{enet1}\\
	&\,+\mathrm{Im} \int_0^T \int_{\mathbb{R}^{d}}    P_{\rho} \left( \psi  u_i u_j  \right) \overline{\psi}_{\rho}  \mathrm{ d} x \mathrm{ d} t\label{enet2}\\
	&\,+\mathrm{Im} \int_0^T \int_{\mathbb{R}^{d}}     [P_{\rho},B] \psi  \cdot  \overline{\psi}_{\rho}  \mathrm{ d} x \mathrm{ d} t. \label{enet3}
\end{align}
\eqref{enet1} can be estimated exactly the same as \eqref{enenltc}, while the estimate of \eqref{enet2} is similar to \eqref{enet1} by noting Remark \ref{Rem_J11t4}.  
Thus, it suffices to estimate
 \eqref{enet3} and one can use \eqref{Bony_com} to obtain
\begin{align}
\eqref{enet3}
	=&\,\rho^{-1}   \int_0^T \int_{\mathbb{R}^{d}}    L(\partial_x B_{\ll \rho },  \psi_{\sim \rho } ) \cdot   \overline{\psi}_{\rho}  \,\mathrm{ d} x \mathrm{ d} t \label{EnBt1}\\
	&\,+     \int_0^T \int_{\mathbb{R}^{d}}    L( B_{ \sim \rho }, \psi_{\ll \rho  } )\cdot \overline{\psi}_{\rho}  \,\mathrm{ d} x \mathrm{ d} t\label{EnBt2}\\
	&\,+ \sum_{\nu\sim \nu^\prime\gtrsim \rho }   \int_0^T \int_{\mathbb{R}^{d}}  P_\rho L(B_{ \nu}, \psi_{\nu^\prime })\cdot  \, \mathrm{ d} x \mathrm{ d} t, \label{EnBt3}
\end{align}
which have been estimated by \eqref{est_EnBt1}--\eqref{est_EnBt3} with the choice \(u_i=u_j=\psi\).

The combination of estimates for \eqref{EnBt1}--\eqref{EnBt3} implies
\begin{align} \label{En_estB}
\eqref{enet3}\lesssim \rho ^{-2s-2\sigma}c_{\rho }^2 (\sigma).
\end{align}
Together with the estimates for \eqref{enet1}, \eqref{enet2}, and \eqref{enenltc}, this implies \eqref{re1} for $ i=1 $.

\section{Bilinear estimate} \label{sectBi}
This section is devoted to establishing the bilinear estimate \eqref{re2}. Without loss of generality, we may always set $ \gamma =1 $ and $\tau=0.$
Replacing $ P_\rho $ in \eqref{Plamdeq} with $ P_\mu $, we obtain
\begin{align} 
	&\sqrt{-1} \pa_t \psi_\mu +\sum_{\alpha,\beta=1}^d \partial_{\alpha} (g^{\alpha \beta} \partial_{ \beta}  \psi_\mu )=P_\mu \mathcal{N}- \partial_{\alpha}  ( [P_\mu, h^{\alpha \beta} ]  \partial_\beta  \psi ):= \mathcal{I}_\mu,	\label{PMeq}
\end{align}
 Similarly, 
 setting  $ Q_\rho =P_{\rho ,1}  $ in \eqref{Q_psi} and $\phi_\rho = \psi_{\rho, 1}$    gives 
 \begin{align} 
	&\sqrt{-1} \partial_t \phi_\rho  +  \partial_{\alpha}  (g^{\alpha \beta}  \partial_{\beta}\phi_\rho  ) = P_{\rho,1} \mathcal{N} -  \partial_{\alpha}  ( [P_{\rho,1}, h^{\alpha\beta} ]   \partial_{ \beta} \psi  )  
 	:= \mathcal{G}_\rho. \label{phiNeq}
\end{align}

Following \cite{SZ}, we can derive the momentum balance law for \eqref{phiNeq} as
\begin{equation} \label{consm2}\begin{aligned}
		&-\partial_t \int_{\mathbb{R}^{d-1}}  \mathrm{Im}\left( \phi_{\rho } \partial_{1} \overline{\phi}_\rho\right) \mathrm{d} \widehat{x}_1 \\
	&+\partial_{1}\Big(    \int_{\mathbb{R}^{d-1}} 2\mathrm{Re} (\partial_{1} \phi_\rho\,   g^{1\beta} \partial_{\beta} \overline{\phi}_{\rho})-\partial_{1 }\mathrm{Re} ( \phi_\rho\,    g^{1\beta} \partial_{\beta} \overline{\phi}_{\rho} ) \mathrm{d} \widehat{x}_1 \Big)\\
	=&\, -\int_{\mathbb{R}^{d-1}}  \mathrm{Re} (\partial_{\alpha} \phi_\rho   \partial_1 h^{\alpha \beta} \partial_\beta \overline{\phi}_{\rho} )  \mathrm{d} \widehat{x}_1  \\
	&\,  + \int_{\mathbb{R}^{d-1}}  \mathrm{Re} (\mathcal{G}_\rho \partial_{1} \overline{\phi}_\rho) -\mathrm{Re} ( \overline{\phi }_{\rho}\partial_{1} \mathcal{G}_\rho )\mathrm{d} \widehat{x}_1 ,      
\end{aligned}
\end{equation}
while the mass balance law for \eqref{PMeq} is given by 
\begin{equation}\label{mas}
\begin{aligned}
	&\frac 12\partial_t \int_{\mathbb{R}^{d-1}}   \left|P_\mu \psi\right|^2  	\mathrm{d} \widehat{x}_1+  \partial_{1} \int_{\mathbb{R}^{d-1}} \mathrm{Im}(   \overline{ \psi}_\mu g^{1\beta}   \partial_{\beta} \psi_\mu  ) \mathrm{d} \widehat{x}_1  \\
	=&\,\int_{\mathbb{R}^{d-1}} \mathrm{Im} (  \mathcal{I}_\mu  \overline{ \psi_\mu }) \mathrm{d} \widehat{x}_1.
\end{aligned}
\end{equation}
Employing the div-curl Lemma \ref{lemdiv-curl} for system \eqref{consm2} and \eqref{mas} yields
\begin{align}  \label{Bie}
	& \mathrm{I}+\mathrm{II}+\mathrm{III}\\
&\lesssim 	\rho \max_{0\leq t \leq T} \| \phi_\rho (t, \cdot) \|_{  L_x^2}^2  \max_{0\leq t \leq T} \|    \psi _\mu(t,\cdot )\|_{ L_x^2}^2 \label{biNt0}\\
	&\,+\Big|\int_{[0,T]\times \mathbb{R}}  \left( \int_{-\infty}^{x_1} \int_{\mathbb{R}^{d-1}}  \mathrm{Im}\left( \phi_{\rho } \partial_{1} \overline{\phi}_\rho\right) 
 \mathrm{d} \widehat{x}_1 \mathrm{ d} y\right) \label{biNt1}\\
 &\qquad\qquad\qquad  \cdot\int_{\mathbb{R}^{d-1}} \mathrm{Im} \left( \overline{ \psi}_\mu  \mathcal{I}_\mu    \right) \mathrm{d} \widehat{x}_1 \mathrm{ d} x_1 \mathrm{ d} t\Big| \nonumber\\
	&\,+\Big|\int_{[0,T]\times \mathbb{R}}  \left( \int^{+\infty}_{x_1}   \left\|\psi_\mu\right\|^2_{L^2_{\widehat{x}_1}}   \mathrm{ d} y\right )\label{non.bit1}\\
 &\qquad\qquad\qquad  \cdot\int_{\mathbb{R}^{d-1}}  \mathrm{Re}\left(\partial_{\alpha} \phi_\rho   \partial_1 h^{\alpha \beta} \partial_\beta \overline{\phi}_{\rho}\right)  \mathrm{d} \widehat{x}_1\mathrm{ d} x_1 \mathrm{ d} t \Big|\nonumber\\
	&\,+\Big|\int_{[0,T]\times \mathbb{R}}  \left( \int^{+\infty}_{x_1}   \left\|\psi_\mu\right\|^2_{L^2_{\widehat{x}_1}}   \mathrm{ d} y\right) 
 \Big|\label{non.bit2}\\
	&\qquad\qquad\qquad  \cdot \int_{\mathbb{R}^{d-1}} \mathrm{Re} \left(\overline{\phi}_\rho \partial_{1} \mathcal{G}_\rho - \mathcal{G}_\rho \partial_{1} \overline{\phi}_\rho\right)   \mathrm{d} \widehat{x}_1 \mathrm{ d} x_1 \mathrm{ d} t \Big|, \nonumber
\end{align}
where
\begin{align*}
	\mathrm{I}&=\int_{[0,T]\times \mathbb{R}} \int_{\mathbb{R}^{d-1}} |  \psi_\mu  |^2 	\mathrm{d} \widehat{x}_1 \int_{\mathbb{R}^{d-1}} \mathrm{Re} (\partial_{1} \phi_\rho   g^{1\beta} \partial_{\beta} \overline{\phi}_{\rho} ) \mathrm{d} \widehat{x}_1 \mathrm{ d} x_1 \mathrm{ d} t,\\
	\mathrm{II}&=\int_{[0,T]\times \mathbb{R}} \int_{\mathbb{R}^{d-1}} \mathrm{Re} (    \psi_\mu  \partial_{1}  \overline{\psi}_\mu   )	\mathrm{d} \widehat{x}_1 \int_{\mathbb{R}^{d-1}} \mathrm{Re} (\phi_\rho \,   g^{1\beta} \partial_{\beta} \overline{\phi}_{\rho} ) \mathrm{d} \widehat{x}_1 \mathrm{ d} x_1 \mathrm{ d} t,\\
	\mathrm{III}&=-\int_{[0,T]\times \mathbb{R}}  \int_{\mathbb{R}^{d-1}}  \mathrm{Im} ( \phi_{\rho } \partial_{1} \overline{\phi}_\rho ) \mathrm{d} \widehat{x}_1 \int_{\mathbb{R}^{d-1}}\mathrm{Im} ( \overline{\psi}_\mu g^{1\beta}    \partial_{\beta}\psi_\mu    ) \mathrm{d} \widehat{x}_1 \mathrm{ d} x_1 \mathrm{ d} t.
\end{align*}

Noting $  g^{\alpha \beta}= h^{\alpha \beta} +\delta_{\alpha \beta}  $, we can  use \eqref{L_formula_4} to obtain
\begin{align*}
	\mathrm{I}=&\, \int_{[0,T]\times \mathbb{R}} \int_{\mathbb{R}^{d-1}} |\psi_\mu |^2  	\mathrm{d} \widehat{x}_1 \int_{\mathbb{R}^{d-1}} |\partial_{1} \phi_\rho |^2 \mathrm{d} \widehat{x}_1 \mathrm{ d} x_1 \mathrm{ d} t\nonumber\\
	&\,	+ \int_{[0,T]\times \mathbb{R}} \int_{\mathbb{R}^{d-1}} |\psi_\mu |^2  	\mathrm{d} \widehat{x}_1 \int_{\mathbb{R}^{d-1}} \mathrm{Re} (\partial_{1} \phi_\rho   h^{1 \beta} \partial_{\beta} \overline{\phi}_{\rho} ) \mathrm{d} \widehat{x}_1 \mathrm{ d} x_1 \mathrm{ d} t\nonumber\\
\gtrsim &\,  \rho^2  \sup_{\tau}\int_{[0,T]\times \mathbb{R}} \int_{\mathbb{R}^{d-1}} |\psi_\mu |^2  \mathrm{d} \widehat{x}_1 \int_{\mathbb{R}^{d-1}} |  \phi_\rho^\tau |^2 \mathrm{d} \widehat{x}_1 \mathrm{ d} x_1 \mathrm{ d} t\\
	&- C\varepsilon^{\frac78} \rho^2 \sup_{\tau}\int_{[0,T]\times \mathbb{R}} \int_{\mathbb{R}^{d-1}} |\psi_\mu |^2  	\mathrm{d} \widehat{x}_1 \int_{\mathbb{R}^{d-1}} | \phi_\rho^\tau |^2 \mathrm{d} \widehat{x}_1 \mathrm{ d} x_1 \mathrm{ d} t \\
    \gtrsim &\,   \rho^2 \sup_{\tau} \int_{[0,T]\times \mathbb{R}} \int_{\mathbb{R}^{d-1}} |\psi_\mu  |^2  	\mathrm{d} \widehat{x}_1 \int_{\mathbb{R}^{d-1}} |  \phi_\rho^\tau |^2 \mathrm{d} \widehat{x}_1 \mathrm{ d} x_1 \mathrm{ d} t
\end{align*}
by taking \(L^\infty\) norm of \(h^{1 \beta} \)  and using \eqref{h_infty}.

Similarly, for $\mathrm{II}$ and $ \mathrm{III}$ we have
\begin{align*}
	\mathrm{II} +
	\mathrm{III}  
    \lesssim  
	 &\,    \rho\mu  \sup_{\tau} \int_{[0,T]\times \mathbb{R}} \int_{\mathbb{R}^{d-1}} |\psi_\mu^\tau  |^2  	\mathrm{d} \widehat{x}_1 \int_{\mathbb{R}^{d-1}} |  \phi_\rho |^2 \mathrm{d} \widehat{x}_1 \mathrm{ d} x_1 \mathrm{ d} t
\end{align*}
where \( \| g^{\alpha \beta} \|_{L^\infty}  \lesssim 1.\)   By the assumption \( \mu \ll \rho\), we finally  obtain 
\begin{align} \label{Bie1}
	\mathrm{I}+\mathrm{II}+\mathrm{III}\gtrsim  \rho^2 \int_{[0,T]\times \mathbb{R}} \int_{\mathbb{R}^{d-1}}   \left|\psi_\mu\right|^2 	\mathrm{d} \widehat{x}_1 \int_{\mathbb{R}^{d-1}} | \phi_\rho |^2 \mathrm{d} \widehat{x}_1 \mathrm{ d} x_1 \mathrm{ d} t.
\end{align} 

We now come to the  estimates of the nonlinear terms \eqref{biNt1}--\eqref{non.bit2}.
For \eqref{biNt1}, we have 
\begin{align*}
    \eqref{biNt1} \lesssim&\, \left\|\phi_{\rho } \partial_{1} \overline{\phi}_\rho \right\|_{L^\infty_tL^1_x }  \left\|\mathrm{Im} \left(  \mathcal{I}_\mu \overline{ \psi}_\mu\right) \right\|_{L^1_{t,x}}\\
    \lesssim &\,  \rho \left\|\phi_{\rho } \right\|_{L^\infty_t L^2_x }^2  \left\|\mathrm{Im} \left(  \mathcal{I}_\mu \overline{ \psi}_\mu\right) \right\|_{L^1_{t,x}} \\
    \lesssim &\, \mu ^{-2s} c_\mu ^2(0)  \rho^{-2s+1-2\sigma} c_\rho^2(\sigma),
\end{align*}
where 
 \( \left\|\phi_{\rho } \right\|_{L^\infty_t L^2_x }^2 \lesssim \left\|\psi_{\rho } \right\|_{L^\infty_tL^2_x }^2 \)
 and  we have established in Section \ref{secEn} that
 \begin{align*}
 \left\|\mathrm{Im} \left(  \mathcal{I}_\mu \overline{ \psi}_\mu\right) \right\|_{L^1_{t,x}}   \lesssim  \varepsilon^{\frac 58} \mu ^{-2s} c_\mu ^2(0).
 \end{align*}

Regarding \eqref{non.bit1},   the frequency  of $ \partial_1  h^{\alpha \beta}  $ in \eqref{non.bit1} are less than  the sum of those for $ \phi_\rho , \overline{\phi}_\rho, \psi_\mu.  $ Hence, we can rewrite and estimate \eqref{non.bit1} as follows:  
\begin{align*}
	 \eqref{non.bit1} =&\,\Big|\int_{[0,T]\times \mathbb{R}}  \left( \int^{+\infty}_{x_1}   \left\|\psi_\mu\right\|^2_{L^\infty_t L^2_{\widehat{x}_1}}   \mathrm{ d} y\right )\nonumber\\
 &\qquad \qquad \cdot\int_{\mathbb{R}^{d-1}}  \mathrm{Re}\left(\partial_{\alpha} \phi_\rho   \partial_1 P_{\lesssim \rho}h^{\alpha \beta} \partial_\beta \overline{\phi}_{\rho}\right)  \mathrm{d} \widehat{x}_1\mathrm{ d} x_1 \mathrm{ d} t \Big|\\
	 \lesssim &\, \left\|\psi_\mu\right\|_{L^2_{x}}^2 \cdot \|\partial_{\alpha} \phi_\rho   P_{\lesssim \rho}\partial_1 h^{\alpha \beta} \partial_\beta \overline{\phi}_{\rho} \|_{L^1_{t,x}}.
\end{align*}
Then \eqref{non.bit1} can be estimated in exactly the same way as \eqref{est_J1} and \eqref{est_J3}.

Moving to \eqref{non.bit2}, we use \eqref{momen_B_3} to obtain
 \begin{align}
   \eqref{non.bit2}  \lesssim &\,     \left\|\psi_\mu\right\|^2_{L^\infty_t L^2_{x}}    \cdot \Big( \rho\,  \|  L( \phi_\rho, P_{\rho,1} (\mathbf{U} \cdot \partial_x \psi+\psi  u_i u_j)+ [P_{\rho,1} , B]\psi )\|_{L^1_{t,x}} \nonumber\\
 & \qquad \qquad      \quad + \|  \pa_x \phi_\rho \cdot \partial_{\alpha} \big([P_{\rho,1}, h^{\alpha\beta} ]   \partial_{ \beta} \psi \big) \|_{L^1_{t,x}}\Big) \label{non.bit2_1}  \\
    & \label{non.bit2_4}\, +\Big| \int_{[0,T]\times \mathbb{R}}  \left( \int^{+\infty}_{x_1}   \left\|\psi_\mu\right\|^2_{L^2_{\widehat{x}_1}}   \mathrm{ d} y\right) \\
  &\qquad \qquad\qquad \cdot  \int_{\mathbb{R}^{d-1}}   \partial_1( [P_{\rho,1} , B]\psi)-\pa_{1}B \cdot |\phi_\rho |^2 \mathrm{d} \widehat{x}_1 \mathrm{ d} x_1 \, \mathrm{ d} t\Big|. \nonumber
 \end{align}
For \eqref{non.bit2_1}, by \eqref{Bony_com} and the same approach as in Section \ref{secEn}, we have 
\begin{equation*} 
   \begin{aligned}
    & \|L \big( \phi_\rho , P_{\rho,1} (\mathbf{U} \cdot \partial_x \psi+ \psi  u_i u_j)+[P_{\rho,1}, B]\psi+\partial_{\alpha} ([P_{\rho,1}, h^{\alpha\beta} ]   \partial_{ \beta} \psi)  \big )   \|_{L^1_{t,x}} \\
    & \lesssim   \varepsilon^{\frac 58} \rho^{-2s-2\sigma} c_\rho^2(\sigma),
\end{aligned} 
\end{equation*}
 which gives  
\begin{align*}
  \eqref{non.bit2_1}  \lesssim   \mu ^{-2s} c_\mu ^2(0)  \rho^{-2s+1-2\sigma} c_\rho^2(\sigma).
\end{align*}

Turning to \( \eqref{non.bit2_4} \), we shall do integration by part for the first term and  it is sufficient to bound 
 \begin{align}
	&\Big|\int_{[0,T]\times \mathbb{R}}  \left( \int^{+\infty}_{x_1}   \left\|\psi_\mu\right\|^2_{L^2_{\widehat{x}_1}}\mathrm{ d} y\right)\cdot \int_{\mathbb{R}^{d-1}} \partial_{1}  B_{\lesssim \rho}   \cdot  | \phi_\rho |^2  \mathrm{d} \widehat{x}_1\, \mathrm{ d} x_1 \mathrm{ d} t\Big|\label{biBt1_1}\\
    &\,+ \Big|  \int_{[0,T]\times \mathbb{R}} \|\psi_{\mu}\|_{L^2_{\hat{x}_1}}^2 \int_{\mathbb{R}^{d-1}}\phi_\rho [P_{\rho,1},B] \psi \,  \mathrm{d}\widehat{x}_1 \,\mathrm{d}x_1 \mathrm{d} t\Big|, \label{biBt1_15}
%	&\,+ \Big|  \int_{[0,T]\times \mathbb{R}} \|\psi_{\mu}\|_{L^2_{\hat{x}_1}}^2 \int_{\mathbb{R}^{d-1}}  B_{\lesssim \rho}\cdot |\phi_\rho|^2 \mathrm{d}\widehat{x}_1 \,\mathrm{d}x_1 \mathrm{d} t\Big|. \label{biBt1_2}
\end{align}
where we    have used the fact  that 
 the frequencies   of $ B $ in \eqref{biBt1_1}  are less than  the sum of those for $ \phi_\rho , \overline{\phi}_\rho, \psi_\mu  $.

The estimate of \eqref{biBt1_1} is analogous to that of \eqref{EnBt1},
since
\[
    \sup_{x_1}
    \int_{x_1}^{+\infty}
        \|\psi_\mu\|_{L^2_{\widehat x_1}}^2\,dy
    \le
    \|\psi_\mu\|_{L^2_x}^2.
\]
The estimate of \eqref{biBt1_15} is analogous to that of \eqref{enet3},
after taking the \(L^\infty_{x_1}\) norm of
\(\|\psi_\mu\|_{L^2_{\widehat x_1}}\).  
%For \eqref{biBt1_2},  it follows from  Bernstein inequality, \eqref{sup2} and \eqref{B_infty}  that
%\begin{align*}
	%\eqref{biBt1_2} \lesssim&\,    \left\|  B_{\lesssim \rho}   \right\|_{L^2_t L^\infty_x} \cdot  \left\| \psi_\mu  \right\|_{ L^\infty_{x_1} L^2_{\widehat{x}_1}}  \big\|  P_{\rho,1} \psi   \big\|_{ L^2_{x_1} L^2_{\widehat{x}_1}} \cdot \big\|\left\| \psi_\mu  \right\|_{L^2_{\widehat{x}_1}}  \left\| P_{\rho,1} \psi  \right\|_{L^2_{\widehat{x}_1}}\big\|_{L^2_{t,x_1}}\\
	%\lesssim &\,   \varepsilon^{\frac 78}\rho \cdot \varepsilon^{-\frac 18 } \mu ^{-s+\frac 12}  c_\mu(0)  \cdot \varepsilon^{-\frac 18 } \rho^{-s-\sigma  } c_\rho(\sigma) \cdot \varepsilon^{-\frac 14 } \mu ^{-s} \rho^{-s-\sigma -\frac 12} c_\mu (0) c_\rho(\sigma)\\
	%\lesssim&\,   \mu ^{-2s} \rho^{-2s-2\sigma +1} c_\mu ^2(0) c_\rho^2(\sigma).
%\end{align*}
Combining these  estimates   implies
\begin{align*}
	\eqref{non.bit2_4}\lesssim&\,   \mu ^{-2s} \rho^{-2s+1-2\sigma } c_\mu ^2(0) c_\rho^2(\sigma),
\end{align*}
which closes the estimates of \eqref{non.bit2}.

Since  \( \eqref{biNt0} \) can be bounded directly by energy estimates, 
 combining \eqref{Bie1} and the estimates of \eqref{Bie}, \eqref{biNt1}--\eqref{non.bit2}, we finally get
\begin{align*}
	&\int_{[0,T]\times \mathbb{R}} \int_{\mathbb{R}^{d-1}}    |\psi_\mu |^2  	\mathrm{d} \widehat{x}_1 \int_{\mathbb{R}^{d-1}} | \phi_\rho |^2 \mathrm{d} \widehat{x}_1 \,\mathrm{ d} x_1 \mathrm{ d} t \\
	\lesssim &\mu ^{-2s} \rho^{-2s-2\sigma -1} c_\mu ^2(0) c_\rho^2(\sigma),
\end{align*}
which is the bilinear estimate \eqref{re2} for $ i,j=1. $

\section{Interaction Morawetz estimate}  \label{secInter_ac}
We now establish the interaction Morawetz estimate \eqref{re3} for 
$i =1$. Without loss of generality, we set $y_J = (x_1,x_2,x_3)$ and $z_J = ( x_4,\cdots, x_d)$ in \eqref{re3}. 
 
 From \eqref{Plamdeq}, we obtain
\begin{align} \label{tPleq}
	&	\sqrt{-1} \pa_t  \psi_\rho +  \partial_{\alpha} (g_{\ll \rho}^{\alpha \beta} \pa_\beta\psi_\rho)   = \mathcal{I}_\rho -   \partial_{ \alpha} (P_{\gtrsim  \rho }  h^{\alpha \beta}\pa_\beta \psi_{\rho} )
	:=\mathcal{I}_\rho^{h}, 
\end{align}
where using \eqref{dt_g} we have 
\begin{equation} \label{Prho.dt_g}
     P_{\ll \rho}\pa_t g_{\alpha \beta}  =P_{\ll \rho}( \pa_x u_4 + u_i u_j).
\end{equation}

Inheriting the notations from the Appendix, we apply Theorem \ref{thm_im} to \eqref{tPleq} and derive the interaction Morawetz estimate:
\begin{align}
	& \label{intaet0} \| \psi_\rho  \|^4_{L^4_{t,y}L^2_{z}}  \\
 &  \lesssim  \rho \|\psi_\rho\|_{L^\infty_t L^2_x}^2 (\|\psi_\rho\|_{L^\infty_t L^2_x}^2 +    \|\operatorname{Im}(\overline{\psi}_\rho \mathcal{I}_\rho^{h}) \|_{L^1_t L^1_x})\label{intaet1} \\
		&\quad + \rho\|\psi_\rho\|_{L^\infty_t L^2_x}^2(\| L(\psi_\rho,\psi_\rho,  P_{\ll \rho} \pa_t g),\  L(\pa_x h_{\lesssim \rho}, \psi_\rho ,\pa_x \psi_\rho)\|_{L^1_{t,x}})\label{intaet2} \\
        &\quad +\int_0^T\Big| \int a_k(y-y')\,   |\psi_\rho(y,z)|^2 \, \operatorname{Re}\big(\mathcal{I}_\rho^{h}  \partial^k \overline{\psi}_\rho - \overline{\psi}_\rho \partial^k\mathcal{I}_\rho^{h} \big)(y',z')  \Big| \mathrm{d}t\label{intaet25}\\
	 & \quad  + \int_0^T\Big|\int(|\nabla_y|^{-2}\pa_y  |\psi_\rho(y,z)|^2  \cdot (\pa_x h_{\ll \rho} \cdot |\psi_\rho|^2)(y,z')\, \Big|\mathrm{d}t\label{intaet3} \\
		&\quad    +\int_0^T\Big|\int   |\nabla_y|^{-2} |\psi_\rho(y,z)|^2 \cdot (\pa_x h_{\ll \rho}  \cdot \pa_x|\psi_\rho|^2)(y,z')\, \Big| \mathrm{d}t\label{intaet4} \\
	&\quad    + \int_0^T\Big|\int(|\nabla_y|^{-1}|\psi_\rho(y,z)|^2 \cdot( \pa_x h_{\ll \rho}  \cdot |\psi_\rho|^2)(y,z')\Big|\mathrm{d}t \label{intaet5}\\
	&\quad  +\int_0^T \Big|  \int R^4_y |\psi_\rho(y,z)|^2 \cdot (h_{\ll \rho}  \cdot|\psi_\rho|^2)(y,z'))\Big|, \mathrm{d}t\label{intaet6}
\end{align}
where \(k\) ranges from 1 to 3.
The term \eqref{intaet1} is estimated as in Sections~\ref{secEn} and \ref{sectBi}, while \eqref{intaet2} is treated as in \eqref{ene_ht1} of Section~\ref{secEn}.

For \eqref{intaet25},  using \eqref{momen_B_3} and integration by parts gives
\begin{align}
    &\eqref{intaet25}\nonumber\\  
  & \lesssim\rho\|\psi_\rho\|_{L^\infty_t L^2_x}^2 \big(\| L(\psi_\rho, P_\rho (\mathbf{U} \cdot \partial_x \psi+\psi  u_i u_j)+ [Q_\rho , B]\psi  ) \|_{L^1_{t,x}}\label{intea250}\\
  &\qquad \qquad \qquad\quad    +\|L\big(\psi_\rho, \pa_\alpha[P_\rho, h^{\alpha\beta }]\pa_\beta \psi  \big)\|_{L^1_{t,x}}\big) \nonumber\\
    &\quad +\int_0^T \Big| \int (\pa a )(y-y')\, |\psi_\rho(y,z)|^2 \label{intea251}\\
    &\qquad \qquad\qquad \cdot \Big[\pa_x h_{\ll \rho} \cdot \big( \psi_\rho \cdot [P_\rho, B]\psi) + \pa_\alpha[P_\rho, h^{\alpha\beta }]\pa_\beta \psi \big)\Big](y',z') \Big|\mathrm{d}t\nonumber\\
    &\quad + \int_0^T \Big| \int (\pa^2 a )(y-y')\, |\psi_\rho(y,z)|^2    \label{intea253}\\
    &\qquad \qquad\quad  \cdot \big( \psi_\rho \cdot [P_\rho, B]\psi + \psi_\rho \cdot \pa_\alpha[P_\rho, h^{\alpha\beta }]\pa_\beta \psi \big) (y',z') \Big|\mathrm{d}t.\nonumber
\end{align}
Term \eqref{intea250} is estimated similarly to \eqref{intaet1} and \eqref{intaet2}; the same applies to \eqref{intea251} after taking the $L^\infty_x $ norm of $\pa_x h_{\ll \rho}$.
 
 Modulo zeroth-order Riesz transforms, which are bounded on the relevant Lorentz spaces, convolution with \(\partial^2a\) is represented by \(|\nabla_y|^{-2}\).  Factoring
\(
|\nabla_y|^{-2}
=
|\nabla_y|^{-1/2}|\nabla_y|^{-3/2},
\)
we estimate \eqref{intea253} as follows:
\begin{align*}
\eqref{intea253}
&\lesssim
\int_0^T\int
\big\|
|\nabla_y|^{-1/2}|\psi_\rho(y,z)|^2
\big\|_{L_y^{2,1}}
\\
&\qquad\quad\cdot
\big\|
|\nabla_y|^{-3/2}
\Big(
\psi_\rho\cdot[P_\rho,B]\psi
+
\psi_\rho\cdot
\partial_\alpha
\big([P_\rho,h^{\alpha\beta}]\partial_\beta\psi\big)
\Big)(y,z')
\big\|_{L_y^{2,\infty}}
\,\mathrm dz\,\mathrm dz'\,\mathrm dt
\\
&\lesssim
\int_0^T\int
\big\|
\nabla_y|\psi_\rho(y,z)|^2
\big\|_{L_y^1}
\\
&\qquad\quad\cdot
\big\|
\Big(
\psi_\rho\cdot[P_\rho,B]\psi
+
\psi_\rho\cdot
\partial_\alpha
\big([P_\rho,h^{\alpha\beta}]\partial_\beta\psi\big)
\Big)(\cdot,z')
\big\|_{L_y^1}
\,\mathrm dz\,\mathrm dz'\,\mathrm dt
\\
&\lesssim
\rho\|\psi_\rho\|_{L_t^\infty L_x^2}^2
\big\|
\psi_\rho\cdot[P_\rho,B]\psi
+
\psi_\rho\cdot
\partial_\alpha
\big([P_\rho,h^{\alpha\beta}]\partial_\beta\psi\big)
\big\|_{L_{t,x}^1}
\\
&\lesssim
\varepsilon^{\frac38}
\rho^{-4s+1-4\sigma}c_\rho^4(\sigma),
\end{align*}
where \(L^{p,q}\) denotes the Lorentz space, and we have used the
Lorentz--Hölder, Lorentz--Sobolev, and endpoint
Hardy--Littlewood--Sobolev inequalities stated in
Subsection~\ref{subsec.not}.

The term \eqref{intaet6} can be estimated by taking the $L^\infty_x$ norm of $h$ and is absorbed into \eqref{intaet0} by the smallness condition.

It remains to estimate \eqref{intaet3}--\eqref{intaet5}. We present the
details only for \eqref{intaet3}, since the other two terms are treated
similarly. By the Lorentz--Hölder, Lorentz--Sobolev, and endpoint
Hardy--Littlewood--Sobolev inequalities, we have
\begin{align*}
\eqref{intaet3}
&\lesssim
\int_0^T\int
\big\|
|\nabla_y|^{-1/2}\partial_y
|\psi_\rho(y,z)|^2
\big\|_{L_y^{2,1}}
\\
&\qquad\quad\cdot
\big\|
|\nabla_y|^{-3/2}
\big(
\partial_xh_{\ll\rho}\cdot|\psi_\rho|^2
\big)(y,z')
\big\|_{L_y^{2,\infty}}
\,\mathrm dz\,\mathrm dz'\,\mathrm dt
\\
&\lesssim
\int_0^T\int
\big\|
\nabla_y\partial_y|\psi_\rho(y,z)|^2
\big\|_{L_y^1}
\big\|
\big(
\partial_xh_{\ll\rho}\cdot|\psi_\rho|^2
\big)(\cdot,z')
\big\|_{L_y^1}
\,\mathrm dz\,\mathrm dz'\,\mathrm dt
\\
&\lesssim
\rho^2
\|\psi_\rho\|_{L_t^\infty L_x^2}^2
\big\|
L(\partial_xh_{\ll\rho},\psi_\rho,\psi_\rho)
\big\|_{L_{t,x}^1}
\\
&\lesssim
\rho
\|\psi_\rho\|_{L_t^\infty L_x^2}^2
\big\|
L(\partial_xh_{\ll\rho},
  \psi_\rho,\partial_x\psi_\rho)
\big\|_{L_{t,x}^1}
\\
&\lesssim
\varepsilon^{\frac38}
\rho^{-4s+1-4\sigma}c_\rho^4(\sigma),
\end{align*}
 where the last inequality then follows from Lemma~\ref{lem_J1}.  This completes the proof.

\section{Elliptic estimates}  \label{sec_ell}
In this section, we are going to  prove \eqref{reV} and \eqref{re1}–\eqref{re3} for \( i \neq 1 \) or \( j \neq 1 \).  
We first prove \eqref{reV}, namely 
\begin{align*}
	\left\|P_{\rho}   u_{\mathfrak{T} }\right\|_{L^2_{t,x } }\lesssim  \rho^{-s-1-\sigma} c_\rho(\sigma), \mathfrak{T} =3,4,5.
\end{align*}
 Similar to \eqref{Bony_U}, by the formula \eqref{L_formula_1}--\eqref{L_formula_3} for    $L$ notation,  we can apply \(P_\rho\) to \eqref{Gen_ellev} to obtain
\begin{align*}
		 	P_{\rho}   u_{\mathfrak{T}}
= &\,   L(h_{\ll  \rho}, P_\rho u_{\mathfrak{T}} )+ \rho^{-1} L(P_\rho u_3, P_{\ll \rho } u_{\mathfrak{T}})\nonumber\\
&\,+ \rho^{-2}\sum_{ \nu\sim \nu^\prime \gtrsim \rho}  P_\rho L (P_\nu u_3,  \partial_x P_{\nu^\prime} u_{\mathfrak{T}} )  +\rho^{-2} \partial_x P_\rho (u_2^2)+\rho^{-2}P_\rho(u_i u_j u_l).
\end{align*}
Taking \(L^\infty_{t,x}\) norm of \(h_{\ll  \rho}\), \( P_{\ll \rho } u_{\mathfrak{T}}\), we estimate the first and second terms on the right-hand side of the above equation directly via \eqref{h_infty} and \eqref{supV}. 
By Remark \ref{Rem_J11t4}, the last term can be estimated in the same way as the second-last term, while the second-last term can be estimated using Lemma \ref{lem_l2uiuj}. The estimate for the remaining third term is straightforward:
\begin{align*}
&\rho^{-2}  \sum_{ \nu\sim \nu^\prime \gtrsim \rho} \left\|P_\rho L (P_\nu u_3, P_{\nu^\prime} \partial_x u_{\mathfrak{T}} )\right\|_{L^2_{t,x } }  \\
 \lesssim&\, \rho^{\frac{d}{2}-2}  \sum_{ \nu\sim \nu^\prime \gtrsim \rho} \left\|  L (P_\nu u_3, P_{\nu^\prime} \partial_x u_{\mathfrak{T}} )\right\|_{L^2_{t } L^1_x } \\
 \lesssim &\,  \rho^{\frac{d}{2}-2}  \nu^\prime \sum_{ \nu\sim \nu^\prime \gtrsim \rho} \left\| P_\nu u_3 \right\|_{L^2_{x }  }    \left\| P_{\nu^\prime} u_{\mathfrak{T}}  \right\|_{L^2_{t,x }  }  \\
 \lesssim &\, \rho^{-s-1-\sigma} c_\rho(\sigma).
\end{align*}

Regarding  \eqref{re1}--\eqref{re3} for $ i\neq 1 $ or $ j\neq 1,$  the case   $ i=2 \text{ or } j=2 $  can be proved similarly with other cases, thus we can further unify the elliptic equations \eqref{Gen_u2} and \eqref{Gen_ellev} into a
more   general equation
\begin{align} \label{Gen_ell_e}
	\Delta u_{\mathfrak{n}}
	&=\partial_x^2 u_1+h \cdot \partial_x^2 u_{\mathfrak{n}}+ \partial_x h \cdot \partial_x u_\mathfrak{n}+ \partial_x(u_l   u_k),
\end{align} 
where 
$$ \mathfrak{n}\in \{2,3,4,5\} ,\   l,k\in \{1,2,3,4,5\} . $$ 
The cubic terms \( u_k u_l u_m \) have been omitted, as they can always be estimated similarly to \( u_k \cdot \partial_x u_l \) using Remark \ref{Rem_J11t4}.

By the property of \(L\) notation and Bony's decomposition, one can employ $\Delta^{-1} P_\rho  $ to \eqref{Gen_ell_e} to  get
\begin{align*}
   & 	\|P_{\rho} u_{\mathfrak{n}}(t) \|_{L^2}\lesssim    \| L(h_{\ll  \rho}, P_{ \sim \rho} u_{\mathfrak{n}} )\|_{L^2}+  \rho^{-1} \|L(P_{\sim \rho } u_3, P_{\ll\rho } u_{\mathfrak{n}})\|_{L^2}\\
    &\,  + \sum_{\nu\sim \nu^\prime \gtrsim \rho}  \|P_\rho L (P_\nu u_3, P_{\nu^\prime} \partial_x u_{\mathfrak{n}} )\|_{L^2}+  \rho^{-1} \|P_\rho\left ( u_l   u_{k}\right)\|_{L^2}  +\| P_\rho u_1 \|_{L^2},
\end{align*}
from which we immediately obtain \eqref{re1} for $ i=2,3,4,5 $, noting that   that the first term on the right-hand side of the above inequality can be estimated by taking \(L^\infty\)  norm of  \( h_{\ll  \rho}\) and using \eqref{h_infty}. The term \( \| P_\rho  u_1 \|_{L^2} \) has already been estimated in Section \ref{secEn}, while the remaining terms can be handled similarly to those in Lemma \ref{lem_l2space_u2}. 

Let us come to the bilinear estimate \eqref{re2}, i.e.
\begin{align}
	&\sup_{\tau} \left\| \|P_\mu u_{i}\|_{L^2_{\widehat{x}_\gamma}}  \|  P_{\rho, \gamma} u_{j}^\tau\|_{L^2_{\widehat{x}_\gamma}}  \right\|_{L^2_{t,x_\gamma}}  \label{bieI_neq1} \\
	&\lesssim \mu^{-s} c_\mu(0) \rho^{-s-\frac 12-\sigma} c_\rho(\sigma), \quad  	\text{ for } ~~\mu\ll \rho, \nonumber
\end{align}
for $i\neq 1 $ or $ j\neq 1. $   Without loss of generality, we may set  $ \gamma=1 $. 

We first consider the case \(j \neq 1\). Denoting index $j$ by $\mathfrak{n}$ in this case, by the properties of the \(L\) notation and Bony's decomposition, we  apply \( P_{\rho,1} \) to \eqref{Gen_ell_e} to obtain
 \begin{align} 
     \Delta 	P_{\rho,1}   u_{\mathfrak{n}}=& \, \, \partial_x^2 P_{\rho, 1}u_1 +  P_{\rho,1 } ( P_{\ll \rho} h \cdot  \partial_x^2P_{ \sim \rho} 
 u_{\mathfrak{n}}+ P_{ \sim \rho} h \cdot  \partial_x^2P_{\ll\rho  } 
 u_{\mathfrak{n}}) \nonumber \\
	&\, + P_{\rho,1 } ( P_{\ll \rho} u_3 \cdot  \partial_x P_{ \sim \rho} 
 u_{\mathfrak{n}}+ P_{  \sim \rho} u_3 \cdot  \partial_xP_{\ll\rho  } 
 u_{\mathfrak{n}})  \label{P_Nn}\\
 	&\, + \rho  P_{\rho,1 } ( L(P_{\ll \rho} u_l,  P_{ \sim \rho} 
 u_{k})+ L(P_{  \sim \rho} u_l,   P_{\ll\rho  } 
 u_{k})) \nonumber\\
 &\,+   \sum_{\nu \sim \nu^\prime \gtrsim \rho}P_{\rho,1 }    L (  P_\nu u_3,   \partial_x P_{\nu^\prime } 
 u_{\mathfrak{n}})+ \rho \sum_{\nu \sim \nu^\prime \gtrsim \rho} P_{\rho,1 }   L (P_{\nu} u_l,     P_{\nu^\prime } 
u_k)   \nonumber\\
\nonumber=&\,  \partial_x^2 P_{\rho, 1}u_1+      h_{\ll \rho} \cdot \partial_x^2 P_{\rho , 1} u_{l}   +      \rho   L( P_{\ll \rho} u_l, P_{ \sim \rho} u_{k}  ) \\
 &\,+  \sum_{\nu \sim \nu^\prime \gtrsim \rho}P_{\rho,1 }    L (  P_\nu u_3,   \partial_x P_{\nu^\prime } 
 u_{\mathfrak{n}}) + \rho \sum_{\nu \sim \nu^\prime \gtrsim \rho}  P_{\rho,1 }   L (P_{\nu} u_l,     P_{\nu^\prime } 
u_k)  , \nonumber
 \end{align}
where \( u_3 =\partial_x h\) and the Lemma \ref{lem_com} gives
\begin{align*}
     [P_{\rho,1}, P_{\ll \rho} h]\partial_x^2 P_{ \sim \rho } u_{l}  =  \rho L( P_{\ll \rho} u_3,   P_{ \sim \rho } u_{l}).   
\end{align*}

 Taking \(u_{\mathfrak n}\) as the second factor in
\eqref{bieI_neq1}, substituting \eqref{P_Nn}, and applying H\"older's
and Bernstein's inequalities together with
\eqref{sup1}--\eqref{sup3} and \eqref{h_infty}, we obtain
\begin{align}
	& \rho^2  \sup_{\tau} \big\| \|P_\mu u_{i}\|_{L^2_{\widehat{x}_1}}  \| P_{\rho,1} u_{\mathfrak{n}}^{\tau}\|_{L^2_{\widehat{x}_1}}  \big\|_{L^2_{t,x_1}}\nonumber\\
    \lesssim &\,  \rho^2 \sup_{\tau} \big\| \|P_\mu u_{i}\|_{L^2_{\widehat{x}_1}}  \| P_{\rho,1} u_1^\tau\|_{L^2_{\widehat{x}_1}}  \big\|_{L^2_{t,x_1}}\nonumber\\
   &\, +  \| P_{\ll \rho} h\|_{L^\infty_x}    \rho^2  \sup_{\tau} \big\| \|P_\mu u_{l}\|_{L^2_{\widehat{x}_1}}  \| P_{\rho,1} u_{k}^{\tau}\|_{L^2_{\widehat{x}_1}}  \big\|_{L^2_{t,x_1}}\nonumber\\
    &+   \rho    \|P_\mu u_{i}\|_{L^2_{\widehat{x}_1} L^\infty_{x_1}}  \cdot    \|  L( P_{\ll \rho} u_l, P_{\sim \rho} u_{k}  )  \big\|_{L^2_{t,x}} \nonumber\\
       &\, +\sum_{\nu \sim \nu^\prime \gtrsim \rho}   \|P_\mu u_{i}\|_{L^2_{\widehat{x}_1} L^\infty_{x_1}}  \cdot    \| P_{\rho,1} L (P_{\nu} u_3,   \partial_x  P_{\nu^\prime }
u_{\mathfrak{n}}) \big\|_{L^2_{t,x}} \nonumber\\ 
&\, +\sum_{\nu \sim \nu^\prime \gtrsim \rho}  \|P_\mu u_{i}\|_{L^2_{\widehat{x}_1} L^\infty_{x_1}}  \cdot  \rho   \| P_{\rho,1} L (P_{\nu} u_l,     P_{\nu^\prime }
u_k) \big\|_{L^2_{t,x}} \nonumber\\ 
	\lesssim &\, \rho^2 \sup_{\tau} \big\| \|P_\mu u_{i}\|_{L^2_{\widehat{x}_1}}  \| P_{\rho,1} u_1^\tau\|_{L^2_{\widehat{x}_1}}  \big\|_{L^2_{t,x_1}}+\varepsilon^{\frac58}\mu^{-s} c_\mu(0) \rho^{-s+\frac 32 -\sigma} c_\rho(\sigma). \label{bieI_neq2}
\end{align}
Indeed, \eqref{pre.bi1} gives
\begin{align*}
    &  \rho     \|P_\mu u_{i}\|_{L^2_{\widehat{x}_1} L^\infty_{x_1}}  \cdot    \|  L( P_{\ll \rho} u_l, P_{\sim \rho} u_{k}  )  \big\|_{L^2_{t,x}}\\
    \lesssim & \rho  \mu^{\frac 12} \|P_\mu u_{i}\|_{L^2_{x}} \cdot \varepsilon^{\frac 34}\sum_{\omega\ll \rho}  \omega^{\frac 12}   \rho^{-s-\frac12-\sigma} c_{\rho} (\sigma) \\
    \lesssim &\, \varepsilon^{\frac58}\mu^{-s} c_\mu(0) \rho^{-s+\frac 32 -\sigma} c_\rho(\sigma).
   \end{align*}
   For the first high--high contribution, H\"older's and Bernstein's
inequalities, together with \eqref{sup1}, \eqref{supV}, and
\eqref{ini3}, yield
   \begin{align}
     & \sum_{\nu \sim \nu^\prime \gtrsim \rho}   \|P_\mu u_{i}\|_{L^2_{\widehat{x}_1} L^\infty_{x_1}}  \cdot    \| P_{\rho,1} L (P_{\nu} u_3,   \partial_x  P_{\nu^\prime }
u_{\mathfrak{n}})  \|_{L^2_{t,x}}  \label{ell.bi.t1} \\
\lesssim &\, \sum_{\nu \sim \nu^\prime \gtrsim \rho}   \|P_\mu u_{i}\|_{L^2_{\widehat{x}_1} L^\infty_{x_1}}  \cdot   \rho^{\frac d2} \|  L (P_{\nu} u_3,   \partial_x  P_{\nu^\prime }
u_{\mathfrak{n}})  \|_{L^2_{t} L^1_x}  \nonumber\\
\lesssim &\, \sum_{\nu \sim \nu^\prime \gtrsim \rho}  \mu^{\frac12} \|P_\mu u_{i}\|_{L^2_x}   \cdot   \rho^{\frac d2} \|  P_{\nu} u_3  \|_{L^2_{t,x}  }    \|    \partial_x  P_{\nu^\prime }
u_{\mathfrak{n}} \|_{ L^\infty_t  L^2_x} \nonumber \\
\lesssim &\, \sum_{\nu \sim \nu^\prime \gtrsim \rho}  \mu^{\frac12-s} c_\mu(0) \cdot \rho^{s+1} \nu^{-s-1-\sigma} c_\nu(\sigma) \cdot (\nu')^{-s+1}   c_{\nu'} (0)\nonumber\\
\lesssim &\, \varepsilon^{\frac58}\mu^{-s} c_\mu(0) \rho^{-s+\frac 32 -\sigma} c_\rho(\sigma) \nonumber
   \end{align}
 Similarly, for the second high--high contribution, we have
  \begin{align}
   &\sum_{\nu \sim \nu^\prime \gtrsim \rho} \|P_\mu u_{i}\|_{L^2_{\widehat{x}_1} L^\infty_{x_1}}  \cdot   \rho  \| P_{\rho,1} L (P_{\nu} u_l,     P_{\nu^\prime }
u_k) \big\|_{L^2_{t,x}} \label{ell.bi.t2} \\
 \lesssim &\sum_{\nu \sim \nu^\prime \gtrsim \rho} \mu^{\frac 12} \rho^{\frac{d-1}{2}} \|P_\mu u_{i}\|_{L^2_{x}  }  \cdot    \| P_{\rho,1} L (P_{\nu} u_l,   P_{\nu^\prime }
u_k) \big\|_{L^2_{t,y} L^1_z} \nonumber\\
   \lesssim &\,\sum_{\nu \sim \nu^\prime \gtrsim \rho}  \mu^{\frac 12} \rho^{\frac{d-1}{2}}   \|P_\mu u_{i}\|_{L^2_{x}} \cdot \| P_{\nu} u_l\|_{ L^4_{t,y} L_z^2} \| P_{\nu^\prime }
u_k\|_{L^4_{t,y} L_z^2 } \nonumber\\
 \lesssim & \,\varepsilon^{-\frac 38} \rho^{\frac {d-1}2} \sum_{\nu \sim \nu^\prime \gtrsim \rho} \mu^{-s+\frac 12} c_\mu(0)   \nu^{-s+\frac 14-\sigma} c_{\nu } (\sigma)  (\nu^\prime) ^{-s+\frac 14} c_{\nu^\prime } (0) \nonumber\\
 \lesssim &\, \varepsilon^{\frac58}\mu^{-s} c_\mu(0) \rho^{-s+\frac 32 -\sigma} c_\rho(\sigma). \nonumber
\end{align}

Since we have already proved the bilinear estimate \eqref{re2} for \(i=j=1\), applying \eqref{bieI_neq2} yields \eqref{bieI_neq1} for \(i=1\) and \(j\neq 1\).  Consequently,  it suffices to establish \eqref{bieI_neq1} for the case \(i\neq 1\).

We now consider \eqref{bieI_neq1} when 
 \(i\neq 1\). Denoting $i = \mathfrak{n}^\prime$ in this case, 
by the property of \(L\) notation and Bony's decomposition,
 we can   apply \( P_{\mu} \) to \eqref{Gen_ell_e} to get
\begin{align}
	 \Delta 	P_{\mu}   u_{\mathfrak{n}^\prime}=&  \,  \partial_x^2 P_{\mu}u_1+   \mu^2  L( h_{\ll \mu},   P_{\sim \mu} u_{l})+   \mu  L( P_{\ll \mu} u_l, P_{\sim \mu} u_{k}  )  \label{P_NJ}  \\
 &\, +  \sum_{\rho \gg \nu \sim \nu^\prime \gtrsim \mu} \nu^\prime P_{\mu }    L (  P_\nu u_3,     P_{\nu^\prime } 
 u_{\mathfrak{n}})+ \sum_{\rho \gg \nu \sim \nu^\prime \gtrsim \mu} 
 \mu P_{\mu} L (P_{\nu} u_l,     P_{\nu^\prime }
u_k)  \nonumber\\
 &\, + \sum_{\nu \sim \nu^\prime \gtrsim \rho}P_{\mu }    L (  P_\nu u_3,   \partial_x P_{\nu^\prime } 
 u_{\mathfrak{n}})+ \sum_{\nu \sim \nu^\prime \gtrsim \rho} 
 \mu P_{\mu} L (P_{\nu} u_l,     P_{\nu^\prime }
u_k)  \nonumber
\end{align}
for \(l,k\in \{1,2,3,4,5\}\).   Replacing \( i\) in \eqref{bieI_neq1} with \( \mathfrak{n}^\prime \),  we can plug the above equality into  \eqref{bieI_neq1}, then use 
   H\"older's inequality, Bernstein's inequality, \eqref{h_infty}, \eqref{sup1}, \eqref{sup2}, \eqref{sup3},   \eqref{bieI_neq1} for \(i=1\) and similar  analysis for \eqref{ell.bi.t1}-\eqref{ell.bi.t2} to obtain
   \begin{align*}
       & \mu^2  \sup_{\tau} \big\| \|P_\mu u_{\mathfrak{n}^\prime}\|_{L^2_{\widehat{x}_1}}  \| P_{\rho,1} u_{j}^{\tau}\|_{L^2_{\widehat{x}_1}}  \big\|_{L^2_{t,x_1}}\nonumber\\
    \lesssim &\,  \mu^2 \sup_{\tau} \big\| \|P_\mu u_{1}\|_{L^2_{\widehat{x}_1}}  \| P_{\rho,1} u_j^\tau\|_{L^2_{\widehat{x}_1}}  \big\|_{L^2_{t,x_1}}\nonumber\\
    &\, +   \left(  \mu^2 \| h_{\ll \mu}  \|_{L^\infty_x}+\mu  \| P_{\leq \mu} u_i\|_{L^\infty_x} \right) \cdot \sup_{\tau} \big\| \|P_{\sim \mu} u_{l}\|_{L^2_{\widehat{x}_1}}  \| P_{\rho,1} u_{j}^{\tau}\|_{L^2_{\widehat{x}_1}}  \big\|_{L^2_{t,x_1}}\\
    &\, + \sum_{\rho \gg \nu \sim \nu^\prime \gtrsim \mu}  \mu^{\frac{d-1}2}  (\nu'+\mu)\sup_{\tau} \big\| \|L(P_{\nu} 
 u_l,P_{\nu^\prime} u_{k})\|_{L^1_{\widehat{x}_1}}  \| P_{\rho,1} u_{j}^{\tau}\|_{L^2_{\widehat{x}_1}}  \big\|_{L^2_{t,x_1}} \\
 &\, +\sum_{ \nu \sim \nu^\prime \gtrsim \rho} \nu' \| P_{\mu} L (P_{\nu} u_3,    P_{\nu^\prime }
u_\mathfrak{n}) \|_{L^2_{t,x}} \cdot \|  P_{\rho,1} u_{j} \|_{L^\infty_{t,x_1}L^2_{x_1} }\\
 &\, +\sum_{ \nu \sim \nu^\prime \gtrsim \rho} \mu \| P_{\mu} L (P_{\nu} u_l,    P_{\nu^\prime }
u_k) \|_{L^2_{t,x}} \cdot \|  P_{\rho,1} u_{j} \|_{L^\infty_{t,x_1}L^2_{x_1} }\\
\lesssim &\, \varepsilon^{\frac58}\mu^{-s+2} c_\mu(0) \rho^{-s-\frac 12 -\sigma} c_\rho(\sigma)\\
&\,+ \varepsilon^{-\frac{3}{8}}    \mu^2 c_\mu(0) \cdot \varepsilon^{-\frac18}  \mu^{-s} c_{\mu}(0) \rho^{-s-\frac 12 -\sigma} c_\rho(\sigma)\\
&\, + \sum_{\rho \gg \nu \sim \nu^\prime \gtrsim \mu}  \mu^{\frac{d-1}2} \nu^{\frac 32}    \| P_{\nu} 
 u_l \|_{L^\infty_t L^2_{x}} \cdot \sup_{\tau} \big\| \| P_{\nu^\prime} u_{k}\|_{L^2_{\widehat{x}_1}}  \| P_{\rho,1} u_{j}^{\tau}\|_{L^2_{\widehat{x}_1}}  \big\|_{L^2_{t,x_1}}  \\
  &\,+ \sum_{  \nu \sim \nu^\prime \gtrsim \rho} \nu'  \mu^{\frac{d}{2}} \|   L (P_{\nu} u_3,    P_{\nu^\prime }
u_\mathfrak{n}) \|_{L^2_{t} L^1_x} \cdot  \rho^{\frac12}\|  P_{\rho} u_{j} \|_{L^2_{x} } \\
 &\,+ \sum_{  \nu \sim \nu^\prime \gtrsim \rho}   \mu^{\frac{d-1}{2}} \|   L (P_{\nu} u_l,    P_{\nu^\prime }
u_k) \|_{L^2_{t,y} L^1_z} \cdot  \rho^{\frac12}\|  P_{\rho} u_{j} \|_{L^2_{x} } \\
\lesssim &\, \varepsilon^{\frac58}\mu^{-s+2} c_\mu(0) \rho^{-s-\frac 12 -\sigma} c_\rho(\sigma)\\
&\,+\varepsilon^{-\frac 38} \sum_{\rho \gg \nu \sim \nu^\prime \gtrsim \mu}  \mu^{\frac{d-1}2} \nu^{\frac 32} \nu^{-s} c_\nu( 0)  (\nu^\prime)^{-s} c_{\nu^\prime}(0) \rho^{-s-\frac 12 -\sigma} c_\rho(\sigma) \\
&\, +\varepsilon^{-\frac 38}\sum_{  \nu \sim \nu^\prime \gtrsim \rho} \mu^{\frac{d}{2}}  \nu^{-s-1} c_{\nu }(0) (\nu^\prime)^{-s+1} c_{\nu^\prime }(0)  \rho^{-s+\frac 12-\sigma} c_{\rho} (\sigma)\\
&\, +\varepsilon^{-\frac 38}\sum_{  \nu \sim \nu^\prime \gtrsim \rho} \mu^{\frac{d-1}{2}}  \nu^{-s+\frac 14} c_{\nu }(0) (\nu^\prime)^{-s+\frac 14} c_{\nu^\prime }(0)  \rho^{-s+\frac 12-\sigma} c_{\rho} (\sigma)\\
\lesssim &\, \varepsilon^{\frac58}\mu^{-s+2} c_\mu(0) \rho^{-s-\frac 12 -\sigma} c_\rho(\sigma).
   \end{align*}
 This completes the proof  of \eqref{bieI_neq1} for \(i\neq 1\) or \(j\neq 1\).

Finally, we prove the interaction Morawetz estimates \eqref{re3} for $  i=\mathfrak{n} \neq 1 $, namely
\begin{align*}
	  \left\|  P_\rho u_{\mathfrak{n}}   \right\|_{L^4_{t,y}L^2_{z} }\lesssim \rho^{-s+\frac 14-\sigma} c_\rho(\sigma).
\end{align*}
Replacing $ \mu$ in \eqref{P_NJ} with $ \rho  $   gives
\begin{align}
	 \Delta 	P_{\rho}   u_{\mathfrak{n}}=&  \,  \partial_x^2 P_{\rho}u_1+  \rho^2  L(P_{\leq \rho} h,   P_{\sim \rho} u_{l})+   \rho   L( P_{\ll \rho} u_l, P_{\sim \rho} u_{k}  )  \label{P.math.n}  \\
 &\, + \sum_{\nu \sim \nu^\prime \gtrsim \rho}P_{\rho }    L (  P_\nu u_3,   \partial_x P_{\nu^\prime } 
 u_{\mathfrak{n}})+  \sum_{  \nu \sim \nu^\prime \gtrsim \rho} 
 \rho  P_{\rho} L (P_{\nu} u_l,     P_{\nu^\prime } u_k), \nonumber
\end{align}
which together with \eqref{sup3} for $  i=  1 $, Bernstein  inequality and H\"older's inequality yields
\begin{align*}
    \rho^2 \left\|  P_\rho u_{\mathfrak{n}}   \right\|_{L^4_{t,y}L^2_{z} } \lesssim &\,  \rho^2\left\|  P_\rho u_{1}   \right\|_{L^4_{t,y}L^2_{z} } +    \rho^2 \|   h_{\ll \rho}\|_{L^\infty_x}   \left\|  P_{\sim \rho} u_{l}   \right\|_{L^4_{t,y}L^2_{z} }\\
    &\,+ \rho   \| P_{\ll \rho} u_l\|_{L^\infty_x}   \left\|  P_{\sim \rho} u_{k}   \right\|_{L^4_{t,y}L^2_{z} }\\
    &\,+  \sum_{  \nu \sim \nu^\prime \gtrsim \rho} 
 (\nu'+\rho) \rho^{\frac{d-3}2 +\frac{3}2 } \left\|     L (P_{\nu} u_l,     P_{\nu^\prime } u_k)   \right\|_{L^4_{t}L^{\frac 43}_y L^1_{z} }\\
 \lesssim &\, \varepsilon^{\frac{3}{4}} \rho^{-s+\frac 94-\sigma} c_\rho (\sigma) +\sum_{  \nu \sim \nu^\prime \gtrsim \rho} 
 \nu' \rho^{\frac{d}2  }  \| P_{  \nu} u_l\|_{L^2_x}   \left\|  P_{\nu^\prime} u_{k}   \right\|_{L^4_{t,y}L^2_{z} }\\
  \lesssim &\, \varepsilon^{\frac{3}{4}} \rho^{-s+\frac 94-\sigma} c_\rho (\sigma) +\varepsilon^{-\frac 14} \sum_{   \nu \sim \nu^\prime \gtrsim \rho} 
   \rho^{\frac{d}2  } \nu^{-s} c_\nu (0)  (\nu')^{-s+\frac 54-\sigma} c_{\nu'} (\sigma) \\
    \lesssim &\, \varepsilon^{\frac{3}{4}} \rho^{-s+\frac 94-\sigma} c_\rho (\sigma) .
\end{align*}
This is exactly \eqref{re3} for $ i= \mathfrak{n} \neq 1 $.

\section{Estimates for \texorpdfstring{$B$}{B}} \label{sec_B}
It remains to  estimate  the unknown function $ B $. As mentioned in Section \ref{sec1},  the scale of $ B $ is different from $ u_i\, (i=1,2,3,4,5) $ and it requires a more refined analysis based on  \eqref{eq:B-schematic}.

Applying  $P_\rho$ to
\eqref{eq:B-schematic} and taking \(L^2_{t,x}\) norm, we obtain
 \begin{align}
    \|\pa_x B_\rho\|_{L^2_{t,x}}  \lesssim&\, \|P_\rho \pa_x(  u_i   u_j)\|_{L_{t,x}^2 }\label{Sp_Bt0} \\
    &\, + \rho^{-1}\|   P_{\sim \rho} (\pa_x u_5\cdot \pa_x  u_4 )\|_{L^2_{t,x}} \label{Sp_Bt1}\\
    &\, + \|P_\rho ( u_3 \cdot B) + \rho^{-1}  P_{\sim \rho} (  u_3\cdot \pa_x B    )\|_{L^2_{t,x}}\label{Sp_Bt2}  \\
    &\, + \|P_\rho (u_i u_j u_l)  + \rho^{-1}  P_{\sim \rho} (  u_i u_j \pa_x u_l + u_i u_j  u_k u_l ) \|_{L^2_{t,x}} \label{Sp_Bt3}.
\end{align}

The high order terms in \eqref{Sp_Bt3} can be estimated in a same manner as the quadratic terms, as the analysis in Remark \ref{Rem_J11t4}. The estimate for \eqref{Sp_Bt0} follows directly from Lemma \ref{lem_l2uiuj}.
 For \eqref{Sp_Bt1}--\eqref{Sp_Bt3}, it suffices to estimate only the terms where the derivative is applied to the product before frequency localization, i.e., terms of the form  $\rho^{-1}P_{\sim\rho}(\partial_x(\cdot)\partial_x(\cdot))$ and $\rho^{-1}P_{\sim\rho}(\cdots\partial_x(\cdot)\cdots)$. Indeed, in the low–high and high–low frequency interaction regimes, the estimates for $P_\rho(\partial_x(uv))$ and $P_\rho(u\partial_x v)$ are equivalent up to constants. In the high–high regime, however, the latter case—where the derivative falls on the product—is more singular and therefore dominates the analysis; controlling it suffices to obtain the desired bounds for the former.

For \eqref{Sp_Bt1}, we estimate 
\begin{align*}
	&\|\rho^{-1}P_{\sim \rho}(\pa_x u_5\pa_x u_4)\|_{L^2_{t,x}}\\
	\lesssim &\,  \rho\left\|L(P_{\ll \rho}   u_4,  P_\rho   u_5)\right\|_{L^2_{t,x}}  + \sum_{\nu^\prime \sim \nu  \gtrsim \rho} \nu^\prime \nu \rho^{\frac d2-1} \left\|	P_{ \nu^\prime}   u_4  P_{ \nu }   u_5\right\|_{L^2_{t }L^1_x}  \\
	  \lesssim& \,   \sum_{\mu \ll \rho }\left\|P_{ \mu}  u_4 \right\|_{L^\infty_x} \rho \left\| P_\rho  u_5\right\|_{L^2_{t,x}} + \sum_{\nu^\prime \sim \nu  \gtrsim \rho} \nu^\prime \nu \rho^{\frac d2 -1} \left\|	P_{ \nu^\prime}   u_4  \right\|_{L^\infty_t L^2_{x}} \left\|P_{ \nu }   u_5\right\|_{L^2_{t,x}}  \\
 \lesssim &\,  \varepsilon^{-\frac 14}\sum_{\mu \ll \rho } \mu^{\frac d2-s} c_{\mu}(0) \rho^{-s-\sigma} c_\rho(\sigma)  + \varepsilon^{-\frac 14} \sum_{\nu^\prime \sim \nu  \gtrsim \rho}    \rho^{\frac d2 -1}  (\nu^\prime)^{- s + 1} c_{\nu^\prime}(0)  \nu^{-s -\sigma} c_{\nu }(\sigma) \\ 
  \lesssim&\, \varepsilon^{\frac 34}\rho^{-s + 1-\sigma} c_\rho(\sigma).
\end{align*}
 
For \eqref{Sp_Bt2}, we estimate 
\begin{align*}
	&\|\rho^{-1}P_{\sim \rho}( u_3 \pa_x B)\|_{L^2_{t,x}} \\
    \lesssim &\,  \|L(P_{\ll \rho}  u_3,    B_{\rho}   )\|_{L^2_{t,x}}  + \|L(P_{\ll \rho}  B,    P_{\rho} u_3 )\|_{L^2_{t,x}}   + \sum_{ \nu\sim \nu^\prime\gtrsim \rho } \nu \rho^{\frac d2 -1 }\| P_{\nu'} u_3 \cdot B_\nu \|_{L^2_{t }L^1_x}\\
   \lesssim&\,    \sum_{\mu \ll \rho }\left\|P_{ \mu}  u_3 \right\|_{L^\infty_x} \left\| B_\rho   \right\|_{L^2_{t,x}} +\sum_{\mu \ll \rho }\left\| B_{ \mu}  \right\|_{L^2_t L^\infty_x}   \left\| P_\rho u_3    \right\|_{L^2_{x}}\\
   &\,+ \sum_{\nu^\prime \sim \nu  \gtrsim \rho}  \nu \rho^{\frac d2 -1} \left\|	P_{ \nu^\prime}   u_3  \right\|_{L^\infty_t L^2_{x}} \left\|B_{ \nu }    \right\|_{L^2_{t,x}}  \\
  \lesssim &\,  \varepsilon^{-\frac 14}\sum_{\mu \ll \rho } \mu^{\frac d2-s} c_{\mu}(0) \rho^{-s-\sigma} c_\rho(\sigma) + \sum_{\mu \ll \rho }\mu^{\frac d2 -s} c_{\mu}(0) \rho^{-s-\sigma} c_\rho(\sigma) \\
  &\,  + \varepsilon^{-\frac 14} \sum_{\nu^\prime \sim \nu  \gtrsim \rho}    \rho^{\frac d2 -1}  (\nu^\prime)^{- s } c_{\nu^\prime}(0)  \nu^{-s + 1-\sigma} c_{\nu }(\sigma) \\
  \lesssim&\,  \varepsilon^{\frac 34}\rho^{-s + 1-\sigma} c_\rho(\sigma).
\end{align*}

Collecting estimates of \eqref{Sp_Bt0}--\eqref{Sp_Bt3},  we  obtain that
\begin{align*}
	\left\|  B_\rho\right\|_{L^2_{t,x}}\lesssim \rho^{-s-\sigma} c_\rho(\sigma),
\end{align*}
which is exactly \eqref{reB2}. 

Finally, combining the estimates obtained in Section \ref{secEn}-\ref{sec_B}, we have finished the proof of Theorem \ref{thm_a_apr}.

\section{Global Well-Posedness in the Critical Besov Space} \label{sec.wpd.5}
In this section, we complete the proof of global well-posedness for the Schr\"odinger--elliptic system
\eqref{gaugsmcf}--\eqref{eq:B-elliptic}. Together with the initial gauge construction and the reconstruction of the geometric flow, this yields global well-posedness for \eqref{smcf1} in the critical Besov space under the harmonic/Coulomb gauge formulation. 
The argument relies on the uniform a priori estimates at the critical regularity obtained in Theorem \ref{thm_a_apr} and a weak difference estimate one derivative below the critical level obtained in this section. To construct solutions for rough initial data, we apply these estimates to a sequence of smooth approximations. The critical estimates provide uniform control and propagate the frequency envelopes, while the weak difference estimate shows that the approximate solutions form a Cauchy sequence in a weaker topology. The same difference estimate also yields uniqueness and continuous dependence on the initial data.

We assume throughout in this section that
\(
s=\frac d2-1,~
d\geq 5,
\)
and  
\(
\|\psi_0\|_{\dot B^s_{2,1}}
\leq \varepsilon,
\)
where \(\varepsilon>0\) is sufficiently small. Let 
\[ \left( \psi^{(a)}, \lambda^{(a)}, h^{(a)}, V^{(a)}, A^{(a)}, B^{(a)} \right), \quad a=1,2, \] 
be two smooth compatible solutions of the Schr\"odinger--elliptic system \eqref{gaugsmcf}--\eqref{eq:B-elliptic} on a time interval \( I=[0,T]. \)
We further set
\[
\left(
u_1^{(a)},
u_2^{(a)},
u_3^{(a)},
u_4^{(a)},
u_5^{(a)}
\right):=\left(
\psi^{(a)},
\lambda^{(a)},
 \partial_x h^{(a)},
V^{(a)},
A^{(a)}
\right)
\]
 and define
\(
\delta u_i:=u_i^{(1)}-u_i^{(2)},~~\delta B:=B^{(1)}-B^{(2)},
\)
for \(i=1,2,3,4,5\). 
 In particular, we write 
 \(
w:=
u_1^{(1)}-u_1^{(2)}=\delta\psi:=\psi^{(1)}-\psi^{(2)}.
\)

Assume that both solutions satisfy the critical a priori bounds with a common
slowly varying frequency envelope \(\{c_\rho\}\) satisfying properties as \eqref{ini1}--\eqref{ini3} and a priori estimates as
\begin{align}
 & \|P_\rho u_i^{(a)}\|_{L^\infty_tL^2_x(I)}
\lesssim
\rho^{-s}c_\rho,
 \label{Ea1}\tag{Ea1}
\\
&\|P_\rho u_{\mathfrak T}^{(a)}\|_{L^2_{t,x}(I)}
\lesssim
\rho^{-s-1}c_\rho,\qquad \mathfrak{T}\in\{3,4,5\}, \label{Ea2}
\tag{Ea2}
\\
&\sup_\tau
\left\|
\|P_\mu u_i^{(a)}\|_{L^2_{\widehat x_\gamma}}
\|P_{\rho,\gamma}u_j^{(a),\tau}\|_{L^2_{\widehat x_\gamma}}
\right\|_{L^2_{t,x_\gamma}(I)}
\lesssim
\mu^{-s}c_\mu\,\rho^{-s-\frac12}c_\rho, \label{Ea3}
\tag{Ea3}
\\
&\|P_\rho u_i^{(a)}\|_{L^4_{t,y}L^2_z(I)}
\lesssim
\rho^{-s+\frac14}c_\rho, \label{Ea4}
\tag{Ea4}   \\
&\|P_\rho B^{(a)}\|_{L^2_{t,x}(I)}
\lesssim
\rho^{-s}c_\rho,  \label{Ea5}
\tag{Ea5}
\end{align}
 where \(\mu\ll\rho; i,j\in\{1,2,3,4,5\};  \mathfrak{T}\in\{3,4,5\}\).

 We measure the difference between the two solutions in a topology one
derivative below the critical regularity, namely in
\(\dot B^{s-1}_{2,1}\). For each dyadic frequency \(\rho\), let
\(D_\rho(I)\) denote the aggregate of the normalized frequency localized
quantities at frequency \(\rho\) arising from the energy, bilinear,
interaction Morawetz, elliptic, and \(B\)-estimates. We then define
the weak difference control norm by
\(
D(I):=\sum_\rho D_\rho(I).
\)
Associated with these dyadic quantities, we define the frequency envelope
by
\(
d_\rho
:=
\frac{1}{D(I)}
\sum_\nu
2^{-\delta\left|\log_2(\rho/\nu)\right|}
D_\nu(I),
\)
where \(0<\delta\ll1\) is fixed sufficiently small. When \(D(I)=0\), we
adopt the convention that \(d_\rho=0\) for every \(\rho\). By construction,
\(
D_\rho(I)\leq D(I)d_\rho.
\)
Moreover,
\[
\sum_\rho d_\rho\lesssim_\delta 1,
\]
and the envelope satisfies the slow-variation property
\begin{align}
d_\mu
\lesssim
\left(
\frac{\max\{\mu,\rho\}}{\min\{\mu,\rho\}}
\right)^\delta
d_\rho,
\qquad
\text{for all dyadic }\mu,\rho.
\label{slow.inva.d}
\end{align}
After a harmless normalization, we may therefore assume that
\[
\sum_\rho d_\rho\leq1.
\]
In particular, if \(\mu\ll\rho\), then
\[
d_\mu
\lesssim
\left(\frac{\rho}{\mu}\right)^\delta d_\rho.
\]

With this choice of frequency envelope, the frequency localized bounds used
in the bootstrap argument take the following form:
\begin{align}
  &  \|w_\rho\|_{L^\infty_tL^2_x(I)}
\leq
D\,\rho^{-s+1}d_\rho,
 \label{UA1} \tag{UA1}  
\\
&\|w_\rho\|_{L^4_{t,y}L^2_z(I)}
\leq
D\,\rho^{-s+\frac54}d_\rho,
 \label{UA2} \tag{UA2}  
\\
&\|P_\rho\delta u_\mathfrak n\|_{L^\infty_tL^2_x(I)}
\leq
D\,\rho^{-s+1}d_\rho,\qquad 
 \label{UA3} \tag{UA3}  
\\
&\|P_\rho\delta u_{\mathfrak T}\|_{L^2_{t,x}(I)}
\leq
D\,\rho^{-s}d_\rho,\qquad  
 \label{UA4} \tag{UA4}  
\\
&\sup_\tau
\left\|
\|P_\mu \delta u_i\|_{L^2_{\widehat x_\gamma}}
\|P_{\rho,\gamma}u_j^{(a),\tau}\|_{L^2_{\widehat x_\gamma}}
\right\|_{L^2_{t,x_\gamma}(I)}
\leq
D\,\mu^{-s+1}d_\mu\,\rho^{-s-\frac12}c_\rho,
 \label{UA5} \tag{UA5}  
\\
&\sup_\tau
\left\|
\|P_\mu u_i^{(a)}\|_{L^2_{\widehat x_\gamma}}
\|P_{\rho,\gamma}\delta u_j^\tau\|_{L^2_{\widehat x_\gamma}}
\right\|_{L^2_{t,x_\gamma}(I)}
\leq
D\,\mu^{-s}c_\mu\,\rho^{-s+\frac12}d_\rho,
 \label{UA6} \tag{UA6}   \\
& \|P_\rho \delta  u_\mathfrak n\|_{L^4_{t,y}L^2_z(I)}
\leq
D\,\rho^{-s+\frac54}d_\rho,
 \label{UA7} \tag{UA7}  
\\
&\|P_\rho\delta B\|_{L^2_{t,x} }
\leq
D\,\rho^{-s+1}d_\rho, 
 \label{UA8} \tag{UA8}  
\end{align}
where \(\mu\ll\rho; i,j\in\{1,2,3,4,5\}; \mathfrak{T}\in\{3,4,5\}; \mathfrak n=2,3,4,5; a=1,2 \) and   \(w_\rho=P_\rho w\).

We write \(D=D(I)\) for the weak difference control norm associated with
\eqref{UA1}--\eqref{UA8}. More precisely, \(D(I)\) is the sum of the
normalized dyadic quantities arising from the energy, bilinear, interaction
Morawetz, elliptic, and \(B\)-field estimates. For example, the components
corresponding to \eqref{UA1} and \eqref{UA3} are, respectively,
\begin{align*}
    D_E:=&\,\sum_\rho
\rho^{2s-2}\|w_\rho\|_{L^\infty_tL^2_x(I)}^2,\\
D_{\mathrm{ell}}(I)
:=&\,
\sum_{\mathfrak{n}=2}^5
\sum_\rho
\rho^{2s-2}
\|P_\rho\delta u_{\mathfrak{n}}\|_{L_t^\infty L_x^2(I)}^2.
\end{align*}
The remaining components are defined analogously by applying the
normalizations appearing in \eqref{UA2} and
\eqref{UA4}--\eqref{UA8}. Consequently, each normalized dyadic quantity is
controlled pointwise by \(D(I)d_\rho\), with the appropriate background
frequency envelope \(c_\rho\) in the mixed bilinear estimates
\eqref{UA5} and \eqref{UA6}.

We measure the initial separation of the two solutions by
 \( D(0) := \|w(0)\|_{\dot B^{s-1}_{2,1}}. \) The corresponding elliptic differences at \(t=0\) are also controlled by \(D(0)\), since the elliptic constraint map is locally Lipschitz in the weak topology.

The uniqueness argument is therefore reduced to a stability estimate for the
full difference norm on \(I\). More precisely, we prove that the
difference remains controlled by its initial size, up to an error carrying a
positive power of the small critical norm of the two background solutions. 
\begin{proposition}[Weak Lipschitz dependence]
\label{prop1}
There exist constants \(C>0\) and \(\kappa>0\), independent of \(T\), such
that
\[
D(I)
\leq
C D(0)
+
C\varepsilon^\kappa D(I).
\]
Consequently, if \(\varepsilon>0\) is sufficiently small, then
\[
D(I)
\lesssim
D(0).
\]
\end{proposition}
The proof follows the frequency localized framework developed for
Theorem~\ref{thm_a_apr}, but is carried out one derivative below the critical
regularity and applied to the difference system. In this setting, the bounds
\eqref{Ea1}--\eqref{UA8} for the background solutions and their differences
play the role of the a priori assumptions \eqref{sup1}--\eqref{supB2}.
After transferring the relevant frequency envelopes to the reference
frequency through \eqref{slow.inva.d}, the dyadic summations are handled by
the sign criterion introduced at the beginning of
Section~\ref{sec_non_est}.

Most nonlinear contributions can therefore be estimated by suitable
difference versions of the arguments already used for
Theorem~\ref{thm_a_apr}. We shall concentrate on the terms that contain a
genuine difference interaction, an additional derivative loss, or a
coefficient difference and hence require a modification of the previous
estimates.

We first derive several consequences of the bootstrap assumptions for the
background and difference variables. These bounds will be used repeatedly in
the energy, bilinear, interaction Morawetz, and elliptic estimates associated
with \eqref{Ea1}--\eqref{UA8}.

\subsection{Preliminary estimates from  assumptions}
We first record the dyadic \(L^\infty_{t,x}\) bounds. The background
variables \(u_i\) carry one frequency power, while the background field
\(B\) carries two. Since the differences are measured one derivative below
the critical regularity, \(\delta u_i\) and \(\delta B\) carry two and three
frequency powers, respectively.
\begin{lemma} Assume \eqref{Ea1}--\eqref{UA8} hold. Then, we have
\begin{align}
       &\|P_\mu u_i^{(a)}\|_{L^\infty_{t,x}}
        \lesssim
        \mu c_\mu,\qquad\|P_\mu B^{(a)}\|_{L^\infty_{t,x}}
        \lesssim
        \mu c_\mu^2,
         \\
          &\|P_\mu \delta u_i\|_{L^\infty_{t,x}}
        \lesssim
        D\mu^2 d_\mu,
        \qquad
        \|P_\mu\delta B\|_{L^\infty_{t,x}}
        \lesssim
        D\mu^3 d_\mu,
\end{align}
for \(a \in \{1,2\}; i \in \{1,2,3,4,5\}. \)
 
\end{lemma}
The next lemma converts the space--time bootstrap bounds into the mixed norms
that will be used to control low-frequency factors in the difference
equations.
\begin{lemma} Assume \eqref{Ea1}--\eqref{UA8} hold. Then, we have
    \begin{align} 
    \|  \partial_x P_{\ll \rho} u_i^{(a)}\|_{L^4_{t,y} L^\infty_z} \lesssim \varepsilon  \rho^{\frac 34} \label{pre.ua.L4.infty1}\\
        \|   w_{\rho}\|_{L^4_{t,y}L^\infty_z} \lesssim \rho^{\frac 34} D d_\rho, \label{pre.uni.L4.es1} \\
        \|\delta B_{\ll \rho} \|_{L^2_{t,y} L^\infty_z} \lesssim \rho^{\frac 12} D d_\rho
 \label{pre.uni.deltB.es1} 
    \end{align}
\end{lemma}

We finally record several anisotropic bilinear estimates. They are the
difference analogues of the corresponding bounds used in the a priori
argument, but one factor is now measured by the difference envelope
\(d_\rho\).

\begin{lemma} \label{lem.bi.uni} Assume \eqref{Ea1}--\eqref{UA8} hold. Then, we have
\begin{align}
   &\sum_{\gamma=1}^d  \big \| \| P_{  \ll \rho} \psi^{(a)} \|_{L^\infty_{\widehat{x}_\gamma}} \| 
 \overline{w}_\rho \|_{L^2_{\widehat{x}_\gamma}} \big\|_{L^2_{ t,x_\gamma}} \lesssim \varepsilon D \sum_{\mu\ll \rho}  \mu^{\frac 12} \rho^{-s+\frac 12} d_\rho, \label{pre.uni.bi.es1}\\
& \sum_{\mu\ll\rho}\sup_\tau \left\|\|L(\delta h_{\ll\mu},P_\mu u_{\mathfrak T}^{(a)}) \|_{L^\infty_z}\|\psi_{\sim\rho,1}^{(b)}\|_{L^2_z}
\|w_{\rho,1}^{\tau}\|_{L^2_z}\right\|_{L^1_{t,y}}  \label{pre.uni.bi.es1+1}\\
 \lesssim  &\,   D^2 \sum_{\mu\ll\rho} \sum_{\omega\ll\mu}  
    \omega^{\frac 12}  \mu^{\frac 12}        \rho^{-2s}  c_\mu c_\rho d_{\omega} d_\rho,  \nonumber \\
&\sum_{\mu\ll\rho}
\sup_\tau 
\left\|
\|L(P_\mu\delta h,P_{\ll\mu}u_{\mathfrak T}^{(a)})\|_{L^\infty_z}
\|\psi_{\sim\rho,1}^{(b)}\|_{L^2_z}
\|w_{\rho,1}^{\tau}\|_{L^2_z}
\right\|_{L^1_{t,y}}  \label{pre.uni.bi.es1+2}\\
   \lesssim &  \, D^2 \sum_{\mu\ll\rho} \sum_{\omega\ll\mu}  
    \omega^{\frac 12}  \mu^{\frac 12}        \rho^{-2s}  c_\mu c_\rho d_{\omega} d_\rho, \nonumber \\
&\sum_{\mu\ll\rho}\sup_\tau \left\|\|L( P_{\ll\mu} \delta u_i,P_\mu u_{j}^{(a)}) \|_{L^\infty_z}\|\psi_{\sim\rho,1}^{(b)}\|_{L^2_z}
\|w_{\rho,1}^{\tau}\|_{L^2_z}\right\|_{L^1_{t,y}}  \label{pre.uni.bi.es2.1}\\
 \lesssim  &\,  \varepsilon^2 D^2 \sum_{\mu\ll\rho} \sum_{\omega\ll\mu}
   \omega^{\frac {3} 2} \mu^{\frac 12}   d_{\omega}    \rho^{-2s} d_\rho,  \nonumber
\end{align}
where \(a\neq b, a,b=1,2; i,j\in\{1,2,3,4,5\}\) and we set \(y_1=x_1\) in \(y=(y_1,y_2,y_3)\).
    
\end{lemma}

The first two lemmas follow directly from the bootstrap assumptions and
Bernstein's inequality. The proof of Lemma~\ref{lem.bi.uni} illustrates the
anisotropic bilinear argument that will be used repeatedly below, so we
include the details.
\begin{proof}
\eqref{pre.uni.bi.es1+1} can be estimated just as (3.10). Indeed
    \begin{align*}
&   \sum_{\mu\ll\rho} \sup_\tau\left\|\|L(\delta h_{\ll\mu},P_\mu u_{\mathfrak T}^{(2)}) \|_{L^\infty_z}\|\psi_{\sim\rho,1}^{(1)}\|_{L^2_z}
\|w_{\rho,1}^{\tau}\|_{L^2_z}\right\|_{L^1_{t,y}}\\
\lesssim &\sum_{\mu\ll\rho} \sum_{\omega\ll\mu}  \omega^{-1} \sup_\tau\left\|\|L(P_\omega \delta u_{3  },P_\mu u_{\mathfrak T}^{(2)}) \|_{L^\infty_z}\|\psi_{\sim\rho,1}^{(1)}\|_{L^2_z}
\|w_{\rho,1}^{\tau}\|_{L^2_z}\right\|_{L^1_{t,y}}\\
\lesssim &\, \sum_{\mu\ll\rho} \sum_{\omega\ll\mu}
  \omega^{\frac {d-3} 2} \mu^{\frac {d-1 }2}     \big \| \|  P_\omega \delta u_3\|_{L^2_{\widehat{x}_1}} \| \psi_{\sim\rho,1}^{(1)}\|_{L^2_{\widehat{x}_1}} \big\|_{L^2_{ t,x_1}} \big \| \| P_{\mu}u_{\mathfrak T}^{(2)}\|_{L^2_{\widehat{x}_1}} \| w_{\rho,1}^{\tau}\|_{L^2_{\widehat{x}_1}} \big\|_{L^2_{ t,x_1}}\\
\lesssim &\, \sum_{\mu\ll\rho} \sum_{\omega\ll\mu}
  \omega^{s-\frac {1} 2} \mu^{s+\frac 12}  D \omega^{-s+1} d_{\omega} \rho^{-s-\frac 12} c_\rho     D \mu^{-s} c_{\mu} \rho^{-s+\frac 12} d_\rho,
\end{align*}
which gives \eqref{pre.uni.bi.es1+1}.

The estimate \eqref{pre.uni.bi.es1+2} is obtained in the same way, but with
the frequency-\(\mu\) difference factor placed in the first bilinear norm
and the lower-frequency background factor placed in the second one. More
precisely,
\begin{align*}
&
\sum_{\mu\ll\rho}
\sup_\tau
\left\|
\|L(P_\mu\delta h,P_{\ll\mu}u_{\mathfrak T}^{(2)})\|_{L^\infty_z}
\|\psi_{\sim\rho,1}^{(1)}\|_{L^2_z}
\|w_{\rho,1}^{\tau}\|_{L^2_z}
\right\|_{L^1_{t,y}}\\
\lesssim &\, \sum_{\mu\ll\rho} \sum_{\omega\ll\mu}
\rho\mu^{\frac {d-3} 2} \omega^{\frac {d-1 }2}  D \mu^{-s+1} d_{\mu} \rho^{-s-\frac 12} c_\rho   D   \omega^{-s} c_{\omega} \rho^{-s+\frac 12} d_\rho 
\end{align*}
The remaining estimates follow from the same anisotropic decomposition.
For \eqref{pre.uni.bi.es1}, one places the low-frequency background factor
in \(L^\infty_{\widehat{x}_\gamma}\) and pairs it with \(w_\rho\) through
the difference bilinear bootstrap bound. For
\eqref{pre.uni.bi.es2.1}, one instead places the low-frequency difference
factor in the corresponding anisotropic norm and pairs the
frequency-\(\mu\) background factor with \(w_\rho\). The additional
smallness follows from the background frequency envelopes. This completes
the proof.

\end{proof}

\subsection{Difference energy estimate} \label{subsec:diff.en}
As \eqref{gaugsmcf}, the equation for \(\psi^{(a)}\), can be written in divergence form as
\[
\sqrt{-1}\partial_t\psi^{(a)}
+\partial_\alpha\big(g^{\alpha\beta,(a)}\partial_\beta\psi^{(a)}\big)
=\mathcal{N}^{(a)},
\]
where \(a=1,2\).
Subtracting the two equations and using \(g^{\alpha\beta,(a)}=\delta^{\alpha\beta}+h^{\alpha\beta,(a)}\), we obtain
\[
\sqrt{-1}\partial_t w
+\partial_\alpha\big(g^{\alpha\beta,(1)}\partial_\beta w\big)
 =
  \delta \mathcal{N}-
\partial_\alpha\big(\delta h^{\alpha\beta}\partial_\beta\psi^{(2)}\big),
\]
where
\(
\delta \mathcal{N}=\mathcal{N}^{(1)}-\mathcal{N}^{(2)}.
\)
Applying \(P_\rho\) and commuting \(P_\rho\) past the low-frequency coefficient yield
 \begin{align}
   \sqrt{-1}\partial_t w_\rho
+\partial_\alpha\big(g^{\alpha\beta,(1)}\partial_\beta w_\rho\big)
= &\,P_\rho\delta \mathcal{N}
-P_\rho\partial_\alpha
\big(\delta h^{\alpha\beta}\partial_\beta\psi^{(2)}\big) \label{eq:diff} \\
&\, -\partial_\alpha\big([P_\rho,h^{\alpha\beta,(1)}]\partial_\beta w\big):=\mathfrak{N}_\rho.\nonumber
 \end{align}

As in Section~\ref{secEn}, we multiply \eqref{eq:diff} by
\(\overline w_\rho\), take the imaginary part, and integrate over
\([0,t]\times\mathbb R^d\). Multiplying the resulting identity by
\(\rho^{2s-2}\) and summing over \(\rho\), we obtain
\begin{align*}
\sum_\rho
\rho^{2s-2}
\|w_\rho(t)\|_{L^2_x}^2
\leq
\sum_\rho
\rho^{2s-2}
\|w_\rho(0)\|_{L^2_x}^2
+
2\sum_\rho
\rho^{2s-2}\mathcal{E}_\rho(t),
\end{align*} 
where
\begin{align}
  \mathcal{E}_\rho(t):=  &\,\Big| \operatorname{Im} \int_0^t\int_{\mathbb R^d}
   P_\rho\delta \mathcal{N}
 \cdot  \overline{w}_\rho
\,\mathrm{d} x \mathrm{d} s \Big| \label{UA1.term1}\\
 &\, +\Big| \operatorname{Im}\int_0^t\int_{\mathbb R^d}
   P_\rho\partial_\alpha
\big(\delta h^{\alpha\beta}\partial_\beta\psi^{(2)}\big)  
 \cdot \overline{w}_\rho
\,\mathrm{d} x \mathrm{d} s \Big| \label{UA1.term2}\\
&\, + \Big|\operatorname{Im} \int_0^t\int_{\mathbb R^d}
 \partial_\alpha\big([P_\rho,h^{\alpha\beta,(1)}]\partial_\beta w\big) \cdot\overline{w}_\rho
\,\mathrm{d} x \mathrm{d} s \Big|. \label{UA1.term3}
\end{align}
The initial term is bounded by \(D(0)^2\). Therefore it remains to prove
\begin{align}
   \sum_\rho
\rho^{2s-2}\mathcal E_\rho(t)
\lesssim
\varepsilon  D^2. \label{Uni.en}  
\end{align}
It is enough to prove that each contribution to \(\mathcal E_\rho\) is either
of the pointwise summable form
\begin{align}
    \mathcal E_\rho (t)
\lesssim
\varepsilon  D^2\rho^{-2s+2}d_\rho^2,  \label{Uni.en1}   
\end{align}
or of the following non-pointwise but summable envelope interaction form:
\begin{align}
\begin{cases}
       \mathcal E_\rho (t)
\lesssim
D^2\rho^{-2s+2}
c_\rho d_\rho
\sum_{\mu\ll\rho}c_\mu d_\mu,  \vspace{0.4cm}\\
       \mathcal E_\rho (t)
\lesssim
D^2\rho^{-2s+2}
c_\rho d_\rho
\sum_{\nu\gtrsim\rho}c_\nu d_\nu, 
 \end{cases}     \label{Uni.en2} 
\end{align}
Both types of terms are summable after multiplication by
\(\rho^{2s-2}\).  Indeed,
the terms as \eqref{Uni.en2} are harmless only after summing in \(\rho\), since
\begin{align*}
 &   \max\left\{\sum_\rho
c_\rho d_\rho
\sum_{\mu\ll\rho}c_\mu d_\mu, \quad\sum_\rho
c_\rho d_\rho
\sum_{\nu\gtrsim\rho}c_\nu d_\nu\right\} \\
&\lesssim \Big(\sum_\rho c_\rho d_\rho\Big)^2
\lesssim\Big(\sum_\rho c_\rho  \Big)^2 \leq
\varepsilon^2.
\end{align*}
Then  \eqref{Uni.en}  follows.

We first treat \eqref{UA1.term3} and \eqref{UA1.term1}. Most of their
contributions can be controlled by suitable adaptations of the corresponding
arguments in the proof of Theorem~\ref{thm_a_apr}.

For \eqref{UA1.term3}, the integration by parts and Bony's decomposition of commutator \eqref{Bony_com} show
\begin{align}
 \eqref{UA1.term3}
        = &\,
        \Big| \int_0^t\int_{\mathbb R^d}
        [P_\rho,h^{\alpha\beta,(1)}]\partial_\beta w
        \,\partial_\alpha\overline{w}_\rho\,\mathrm{d}x\,d\tau \Big | \nonumber\\
\leq  &\,   \Big|   \int_0^T \int_{\mathbb{R}^{d}} L ( \partial_x P_{\ll \rho }  h^{\alpha \beta,(1)},   w_{\sim \rho })\cdot \partial_{\alpha} \overline{w}_{\rho} \,\mathrm{ d} x\mathrm{ d} t \Big|  \label{UA1.term3.t1}\\
	&\quad +  \Big|   \int_0^T \int_{\mathbb{R}^{d}} L (\partial_{ \beta} P_{\ll \rho  }   w, P_{\sim \rho } h^{ \alpha \beta,(1)})\cdot \partial_{\alpha} \overline{w}_{\rho}\, \mathrm{ d} x\mathrm{ d} t \Big|  \label{UA1.term3.t2}\\
	&\quad +\Big| \sum_{\nu^{\prime} \sim \nu \gtrsim \rho }  \int_0^T \int_{\mathbb{R}^{d}} P_{\rho} L( P_{\nu} h^{\alpha \beta,(1)},  \partial _{ \beta}   
 w_{\nu^\prime })\cdot\partial _{\alpha} \overline{w}_{\rho} \,\mathrm{ d} x\mathrm{ d} t\Big|.     \label{UA1.term3.t3}
\end{align}
The \eqref{UA1.term3.t1}--\eqref{UA1.term3.t3} can be bounded just as \eqref{ene_ht1} -\eqref{ene_ht3} by replacing \(\psi\) in \eqref{ene_ht1}--\eqref{ene_ht3} to \(w\)  in  \eqref{UA1.term3.t1}--\eqref{UA1.term3.t3}  and using the assumption \eqref{Ea1}--\eqref{UA8}. For instance, the same
argument used for \eqref{est_J3}, together with
\eqref{pre.uni.L4.es1} and \eqref{slow.inva.d}, yields
\begin{align*}
    \eqref{UA1.term3.t3}  \lesssim &\, \sum_{\nu^{\prime} \sim \nu \gtrsim \rho } \rho \big\|    \partial_x h^{(1)}  \big\|_{L^2_{t,x}}\left\|   w_{\nu^\prime }  \right\|_{L^4_{t,y}L^2_z}  \|   w_{\rho}\|_{L^4_{t,y}L^\infty_z}\\
\lesssim &\,  \varepsilon D^2 \rho^{-2s+2} d_\rho^2.
\end{align*}
Here, the
scale invariant total frequency power of the expression is \(-2s+2\), by \eqref{pre.uni.L4.es1}, the
factor
\(
        \|w_\rho\|_{L^4_{t,y}L^\infty_z}
        \lesssim \rho^{\frac34}D d_\rho
\)
contributes \(\rho^{3/4}\), and excluding the extra \(\rho\), the residual high-frequency power for high frequency \( \nu \sim \nu'\) is
\(-2s+\frac14\), which is negative since \(s=\frac d2-1\ge1\).
 Thus, \( d_\mu \) can be transformed into \( d_\rho\) by \eqref{slow.inva.d} and the frequency can be converted to the \(\rho\), as \( \nu\sim \nu^\prime \gtrsim \rho\). This is precisely an application of the dyadic summation principle stated after Proposition~\ref{prop1}.

We next turn to \eqref{UA1.term1}. We first verify that the Gauss structure
is preserved after taking differences, in the same way as in the proof of
Lemma~\ref{lem_J11t4}. Set
\begin{align}
\delta(\lambda\wedge\lambda)
:=
(\lambda^{(1)}\wedge\lambda^{(1)})
-
(\lambda^{(2)}\wedge\lambda^{(2)}).
\label{def:difference-gauss}
\end{align}
By definition of \eqref{wed.lam}, we have
\begin{align}
\delta(\lambda\wedge\lambda)_{\gamma\sigma}
=&\,
\delta\lambda_{\gamma\sigma}
\overline{\psi^{(1)}}
+
\lambda_{\gamma\sigma}^{(2)}
\overline{w}
 -
\delta\lambda_{\alpha\gamma}
\overline{(\lambda^{(1)})^\alpha_\sigma}
-
\lambda_{\alpha\gamma}^{(2)}
\overline{
\delta\bigl(\lambda^\alpha_\sigma\bigr)
},
\label{difference-gauss-exact}
\end{align}
where
\begin{align}
\delta\bigl(\lambda^\alpha_\sigma\bigr)
=&\,
g^{\alpha\beta,1}\delta\lambda_{\beta\sigma}
+
\delta g^{\alpha\beta}\lambda_{\beta\sigma}^{(2)}.
\label{difference-raised-lambda}
\end{align}
Consequently, modulo the cubic term
\(
\lambda^{(2)}\cdot\delta g\cdot\lambda^{(2)},
\)
the quadratic part of \eqref{difference-gauss-exact} is the polarized
Gauss interaction
\begin{align*}
\delta\lambda_{\gamma\sigma}
\overline{\psi^{(1)}}
+
\lambda_{\gamma\sigma}^{(2)}
\overline{w}
-
\delta\lambda_{\alpha\gamma}
\overline{(\lambda^{(1)})^\alpha_\sigma}
-
\lambda_{\alpha\gamma}^{(2)}
g_1^{\alpha\beta}
\overline{\delta\lambda_{\beta\sigma}}.
\end{align*}

Subtracting the Codazzi equations  \eqref{eq:codazzi-complex} for the two solutions gives
schematically
\begin{align}
\partial_\beta\delta\lambda_{\gamma\sigma}
=&\,
\partial_\sigma\delta\lambda_{\gamma\beta}
+
L\bigl((u_3^{(1)},A^{(1)}),\delta\lambda\bigr)
\label{difference-codazzi} \\
&+
L\bigl((\delta u_3,\delta A),\lambda^{(2)}\bigr)
+
L\bigl(\delta h,\partial_x\lambda^{(2)}\bigr)
+
u_i^{(a)}u_j^{(b)}\delta u_l, \nonumber
\end{align}
while the contracted difference Codazzi equation takes the form
\begin{align}
\partial_\alpha
\delta\bigl(\lambda^\alpha_\sigma\bigr)
-
\partial_\sigma w
=&\,
L\bigl((u_3^{(1)},A^{(1)}),\delta\lambda\bigr)
\label{difference-contracted-codazzi}\\
&+
L\bigl((\delta u_3,\delta A),\lambda^{(2)}\bigr)
+
L\bigl(\delta h,\partial_x\lambda^{(2)}\bigr)
+
u_i^{(a)}u_j^{(b)}\delta u_l, \nonumber
\end{align}
where \(a,b \in \{1,2\}\).

For \(\nu\sim\nu'\gtrsim\rho\), we  decompose the difference Gauss
interaction as
\begin{align}
P_\mu\delta
\bigl(\lambda_\nu\wedge\lambda_{\nu'}\bigr)
=&\,
P_\mu
\left[
(P_\nu\delta\lambda_{\gamma\sigma})
(P_{\nu'}\overline{\psi^{(1)}})
-
(P_\nu\delta\lambda_{\alpha\gamma})
P_{\nu'}\overline{(\lambda^{(1)})^\alpha_\sigma}
\right]
\nonumber\\
&+
P_\mu
\left[
(P_\nu\lambda_{\gamma\sigma}^{(2)})
(P_{\nu'}\overline{w})
-
(P_\nu\lambda_{\alpha\gamma}^{(2)})
P_{\nu'}\overline{
\delta(\lambda^\alpha_\sigma)}
\right]
\nonumber\\
&+
P_\mu L
\bigl(
\delta h,\lambda^{(2)},\lambda^{(2)}
\bigr).
\label{difference-gauss-polarization}
\end{align}
The last term arises from
\[
\delta(\lambda^\alpha_\sigma)
=
g_1^{\alpha\beta}\delta\lambda_{\beta\sigma}
+
\delta g^{\alpha\beta}\lambda_{\beta\sigma}^{(2)}.
\]
Thus, the first two lines of \eqref{difference-gauss-polarization} are the
quadratic polarized Gauss interactions, whereas the contribution involving
\(\delta g\) is cubic.

Using \eqref{Pnu.derivative} 
and the difference Codazzi identities
\begin{align}
\partial_\beta\delta\lambda_{\gamma\sigma}
=&\,
\partial_\sigma\delta\lambda_{\gamma\beta}
+
L\bigl((u_3^{(1)},A^{(1)}),\delta\lambda\bigr)
\nonumber\\
&+
L\bigl((\delta u_3,\delta A),\lambda^{(2)}\bigr)
+
L\bigl(\delta h,\partial_x\lambda^{(2)}\bigr),
\label{difference-codazzi-1}
\\
\partial_\beta\delta\lambda_{\alpha\gamma}
=&\,
\partial_\alpha\delta\lambda_{\beta\gamma}
+
L\bigl((u_3^{(1)},A^{(1)}),\delta\lambda\bigr)
\nonumber\\
&+
L\bigl((\delta u_3,\delta A),\lambda^{(2)}\bigr)
+
L\bigl(\delta h,\partial_x\lambda^{(2)}\bigr),
\label{difference-codazzi-2}
\end{align}
we obtain
\begin{align}
& P_\mu
\left[
(P_\nu\delta\lambda_{\gamma\sigma})
(P_{\nu'}\overline{\psi^{(1)}})
-
(P_\nu\delta\lambda_{\alpha\gamma})
P_{\nu'}\overline{(\lambda^{(1)})^\alpha_\sigma}
\right]
\nonumber\\
=&\,
\nu^{-1}P_\mu
\left[
\partial_\sigma
\bigl(
\widetilde P_{\nu,\beta}
\delta\lambda_{\gamma\beta}
\bigr)
P_{\nu'}\overline{\psi^{(1)}}
-
\partial_\alpha
\bigl(
\widetilde P_{\nu,\beta}
\delta\lambda_{\beta\gamma}
\bigr)
P_{\nu'}\overline{(\lambda^{(1)})^\alpha_\sigma}
\right]
\nonumber\\
&+
\nu^{-1}P_\mu
L\left(
P_\nu
\left[
(u_3^{(1)},A^{(1)})\delta\lambda
+
(\delta u_3,\delta A)\lambda^{(2)}
+
\delta h\,\partial_x\lambda^{(2)}
\right],
P_{\nu'}\lambda^{(1)}
\right) \nonumber\\
=&\,\frac{\mu}{\nu}
P_\mu
L\bigl(
P_\nu\delta\lambda,
P_{\nu'}\lambda^{(1)}
\bigr)
\nonumber\\
&+
\nu^{-1}P_\mu
L\left(
P_\nu
\left[
(u_3^{(1)},A^{(1)})\delta\lambda
+
(\delta u_3,\delta A)\lambda^{(2)}
+
\delta h\,\partial_x\lambda^{(2)}
\right],
P_{\nu'}\lambda^{(1)}
\right)
\nonumber\\
&+
\nu^{-1}P_\mu
L\left(
P_\nu\delta\lambda,
P_{\nu'}
\left[
(u_3^{(1)},A^{(1)})\lambda^{(1)}
\right]
\right).
\label{first-polarized-final}
\end{align}
  
Similarly,   we can use the ordinary Codazzi equations for
\(\lambda^{(2)}\) to  obtain
\begin{align}
& P_\mu
\left[
(P_\nu\lambda_{\gamma\sigma}^{(2)})
(P_{\nu'}\overline{w})
-
(P_\nu\lambda_{\alpha\gamma}^{(2)})
P_{\nu'}\overline{\delta(\lambda^\alpha_\sigma)}
\right]
\nonumber\\
=&\, \frac{\mu}{\nu}
P_\mu
L\bigl(
P_\nu\lambda^{(2)},
P_{\nu'}\delta\lambda
\bigr)
\nonumber\\
&+
\nu^{-1}P_\mu
L\left(
P_\nu
\left[
(u_3^{(2)},A^{(2)})\lambda^{(2)}
\right],
P_{\nu'}\delta\lambda
\right)
\label{second-polarized-final}\\
&+
\nu^{-1}P_\mu
L\left(
P_\nu\lambda^{(2)},
P_{\nu'}
\left[
(u_3^{(1)},A^{(1)})\delta\lambda
+
(\delta u_3,\delta A)\lambda^{(2)}
+
\delta h\,\partial_x\lambda^{(2)}
\right]
\right),
\nonumber
\end{align}
 
Combining
\eqref{difference-gauss-polarization},
\eqref{first-polarized-final}, and
\eqref{second-polarized-final}, we obtain
\begin{align}
P_\mu\delta
\bigl(
\lambda_\nu\wedge\lambda_{\nu'}
\bigr)
=&\,
\frac{\mu}{\nu}
P_\mu
L\bigl(
P_\nu\delta\lambda,
P_{\nu'}\lambda^{(1)}
\bigr)
 +
\frac{\mu}{\nu}
P_\mu
L\bigl(
P_\nu\lambda^{(2)},
P_{\nu'}\delta\lambda
\bigr)
\label{d4.polarized.codazzi}\\
&+
\nu^{-1}P_\mu
L\left(
P_\nu
\left[
(u_3^{(1)},A^{(1)})\delta\lambda
+
(\delta u_3,\delta A)\lambda^{(2)}
+
\delta h\,\partial_x\lambda^{(2)}
\right],
P_{\nu'}\lambda^{(1)}
\right)
\nonumber\\
&+
\nu^{-1}P_\mu
L\left(
P_\nu\delta\lambda,
P_{\nu'}
\left[
(u_3^{(1)},A^{(1)})\lambda^{(1)}
\right]
\right)
\nonumber\\
&+
\nu^{-1}P_\mu
L\left(
P_\nu
\left[
(u_3^{(2)},A^{(2)})\lambda^{(2)}
\right],
P_{\nu'}\delta\lambda
\right)
\nonumber\\
&+
\nu^{-1}P_\mu
L\left(
P_\nu\lambda^{(2)},
P_{\nu'}
\left[
(u_3^{(1)},A^{(1)})\delta\lambda
+
(\delta u_3,\delta A)\lambda^{(2)}
+
\delta h\,\partial_x\lambda^{(2)}
\right]
\right)
\nonumber\\
&+
P_\mu
L\bigl(
\delta h,\lambda^{(2)},\lambda^{(2)}
\bigr). \nonumber
\end{align}
Thus the two quadratic polarized interactions retain the crucial
output-frequency gain \(\mu/\nu\), while every remaining term contains
one difference factor and at least two background factors.

By   \eqref{Gen_N}, we have
\begin{align}
\delta \mathcal{N}= &\,\mathcal{N}^{(1)}-\mathcal{N}^{(2)}   \nonumber \\
=&\,\delta  \mathbf{U} \cdot\partial_x\psi^{(1)}
+\mathbf{U}^{(2)} \cdot\partial_xw+
wB^{(1)}
+\psi^{(2)}\delta B
\nonumber\\
&\quad
+w u_i^{(1)}u_j^{(1)}
+\psi^{(2)}\delta u_i u_j^{(1)}
+\psi^{(2)}u_i^{(2)}\delta u_j, \label{diff_N}
\end{align}
where \(\mathbf{U}\) is given in \eqref{eq_U}.
It suffices to consider the quadratic terms in \eqref{diff_N}, since the cubic terms can be handled by adapting the general argument in Remark~\ref{Rem_J11t4}. Moreover, as in the treatment of \eqref{UA1.term3}, the terms \(\mathbf{U}^{(2)} \cdot\partial_xw\) and \(
wB^{(1)}\) can be estimated by adapting the corresponding arguments in the proof of Theorem~\ref{thm_a_apr}.
For instance, by the fact \(B^{(1)}\) is real, we have
\begin{align*}
   &
        \operatorname{Im}\int_0^t\int_{\mathbb R^d}
        P_\rho\bigl(wB^{(1)} \bigr)
        \overline{w}_\rho\,\mathrm{d} x \mathrm{d} s  =\operatorname{Im}\int_0^t\int_{\mathbb R^d}
       ( [P_\rho,B^{(1)}]w )
        \overline{w}_\rho\, \mathrm{d} x \mathrm{d} s.
\end{align*}
This has the same commutator structure as \eqref{enet3} and can therefore be controlled by the arguments used in \eqref{EnBt1}--\eqref{EnBt3}.

For  \(\delta  \mathbf{U} \cdot\partial_x\psi^{(1)}\), the Bony's decomposition gives
\begin{align}
  P_\rho(\delta  \mathbf{U} \cdot\partial_x\psi^{(1)}) 
=&\,
 \delta \mathbf{U}_{ \ll\rho} \cdot\partial_x\psi^{(1)}_\rho 
+
 \delta \mathbf{U}_{ \rho} \cdot\partial_x\psi^{(1)}_{\ll\rho} \nonumber\\
&\,+
P_\rho\sum_{\nu \sim \nu'\gtrsim\rho}
\delta \mathbf{U}_{ \nu} \cdot\partial_x\psi^{(1)}_{\nu'}.   \label{bony.diff.uI}
\end{align}
The high--low and high--high interactions can be handled by the standard
dyadic summation argument. 
Indeed,  the high--high contribution in \eqref{bony.diff.uI} can be estimated in the same manner as \eqref{UA1.term3.t3}.
While for the high--low contribution, H\"older's inequality, the
\eqref{pre.uni.bi.es1}, \eqref{UA4}, and \eqref{slow.inva.d} imply
\[
\begin{aligned}
 &\left|
\int_0^t\int_{\mathbb R^d}
 \delta \mathbf{U}_{ \rho} \cdot\partial_x\psi^{(1)}_{\ll\rho} \overline{w}_\rho d x \, ds
\right|\\
\lesssim&\, \| P_\rho\delta u_{\mathfrak{T}} \|_{L^2_{t,x}}    \cdot
   \big \| \|  \partial_x P_{ \ll \rho} \psi^{(1)} \|_{L^\infty_{\widehat{x}_\gamma}} \| 
 \overline{w}_\rho \|_{L^2_{\widehat{x}_\gamma}} \big\|_{L^2_{ t,x_\gamma}}
\\
\lesssim &\,  \varepsilon D^2 \rho^{-2s+2} d_\rho^2,
\end{aligned}
\]
for \(\mathfrak{T} =4\) or 5.
Indeed by \eqref{pre.uni.bi.es1},  the index of low frequency \(\mu\) is positive and, hence the sum over \(\mu\ll\rho\) is dominated by its largest admissible frequency  namely \(\rho\), in accordance with the dyadic summation principle stated after Proposition~\ref{prop1}.

For the low--high interactions, we adapt the arguments used in the proof of Lemma~\ref{lem_J1} to close the estimates.  As \eqref{lem_J1.decom}, we can use the decomposition \eqref{decomp_Xi} and only consider the case \(\gamma=1\). Then, we obtain
 \begin{align}
  & \Big |
\int_0^t\int_{\mathbb R^d}
P_\rho
\big(
\delta \mathbf{U} _{\ll\rho}
\partial_x\psi^{(1)}_\rho
\big)
\overline{w}_\rho
\mathrm{d} x \mathrm{d} s 
\Big |   \nonumber\\
\lesssim &\,
\sum_{\mu\ll\rho}
\rho
\sup_\tau
\left\|
\| \delta\mathbf{U} _{\mu}\|_{L^\infty_z}
\|\psi^{(1)}_{\sim\rho,1}\|_{L^2_z}
\|w_{\rho,1}^{\tau}\|_{L^2_z}
\right\|_{L^1_{t,y}}. \label{deltU.lh}
 \end{align}
By \eqref{eq_U}, 
we have
\begin{align*}
    \Delta\delta \mathbf{U}
= &\,
h^{(1)} \cdot\partial_x^2\delta \mathbf{U}
+
\delta h \cdot \partial_x^2 \mathbf{U}^{(2)}
+
u_3^{(1)}\partial_x\delta\mathbf{U}\\
&\,+
\delta u_3\partial_x \mathbf{U}^{(2)}
+
\partial_x\big( \delta (\lambda\wedge\lambda)\big)
+
\delta(u_i u_j u_k).
\end{align*}
Hence, schematically,
\begin{align}
P_\mu\delta \mathbf{U}
&=
L(h^{(1)}_{\ll\mu},P_\mu\delta \mathbf{U})
+
L(\delta h_{\ll\mu},P_\mu \mathbf{U}^{(2)}) \nonumber
\\
&\quad
+
\mu^{-1}L(P_\mu u_3^{(1)},P_{\ll\mu}\delta \mathbf{U})+
\mu^{-1}L(P_\mu\delta u_3,P_{\ll\mu}\mathbf{U}^{(2)})
\label{Bonydelta_uI}\\
&\quad
+
\mu^{-2}
\sum_{ \rho \gg \nu\sim\nu'\gtrsim\mu} (
P_\mu L(P_\nu\delta u_3,P_{\nu'}\partial_xu_{\mathfrak T}^{(2)})+P_\mu L(P_\nu u_3^{(1)},P_{\nu'}\partial_x\delta u_{\mathfrak T})) \nonumber
\\
&\quad
+
 \mu^{-2}
\sum_{\nu\sim\nu'\gtrsim\rho} (
P_\mu L(P_\nu\delta u_3,P_{\nu'}\partial_xu_{\mathfrak T}^{(2)})+P_\mu L(P_\nu u_3^{(1)},P_{\nu'}\partial_x\delta u_{\mathfrak T})) 
\nonumber\\
&\quad
+
\mu^{-1}P_\mu\big(  \delta(\lambda\wedge\lambda)\big)
+
\mu^{-2}P_\mu\delta(u_i u_j u_k). \nonumber
\end{align}
The first term on the right-hand side of \eqref{Bonydelta_uI} is absorbed
by the smallness of \(h^{(1)}\), while the cubic contributions are covered
by Remark~\ref{Rem_J11t4}. For the second term,
\eqref{pre.uni.bi.es1+1} gives
 \begin{align*}
& \rho \sum_{\mu\ll\rho}\sup_\tau\left\|\|L(\delta h_{\ll\mu},P_\mu \mathbf{U}^{(2)}) \|_{L^\infty_z}\|\psi_{\sim\rho,1}^{(1)}\|_{L^2_z}
\|w_{\rho,1}^{\tau}\|_{L^2_z}\right\|_{L^1_{t,y}} 
\lesssim   \varepsilon^2 D^2  \rho^{-2s+2} d_\rho^2.
\end{align*}
 Similarly by \eqref{pre.uni.bi.es1+2}, we can estimate  the  third and fourth  term in \eqref{Bonydelta_uI}.

The fifth and sixth term in \eqref{Bonydelta_uI} can be estimated in the same manner. We therefore estimate only the fifth term. The corresponding bilinear bounds yield
\begin{align*}
    &\sum_{\mu\ll\rho} 
\sum_{ \rho \gg \nu\sim\nu'\gtrsim\mu}  \mu^{-2} \rho \sup_\tau\left\|\|P_\mu L(P_\nu\delta u_3,P_{\nu'}\partial_xu_{\mathfrak T}^{(2)} \|_{L^\infty_z}\|\psi_{\sim\rho,1}^{(1)}\|_{L^2_z}
\|w_{\rho,1}^{\tau}\|_{L^2_z}\right\|_{L^1_{t,y}}\\
\lesssim&     \sum_{\mu\ll\rho} \sum_{ \rho \gg \nu\sim\nu'\gtrsim\mu}  \mu^{d-3} \rho \nu^\prime \big \| \| P_\nu\delta u_3\|_{L^2_{\widehat{x}_1}} \| \psi_{\sim\rho,1}^{(1)}\|_{L^2_{\widehat{x}_1}} \big\|_{L^2_{ t,x_1}} \big \| \| P_{\nu^\prime}u_{\mathfrak T}^{(2)}\|_{L^2_{\widehat{x}_1}} \| w_{\rho,1}^{\tau}\|_{L^2_{\widehat{x}_1}} \big\|_{L^2_{ t,x_1}}\\
\lesssim &\, \varepsilon^2  D^2 \rho^{-2s+2} d_\rho^2,
\end{align*}
where \(d \geq 5\).
Note the fact that scale-invariant total frequency power of the expression is \(-2s+2\), the power of \(\rho \) is \( -2s+1\), we know the sum of power of the  low frequencies is \(1\). The desired estimate then follows from the dyadic summation principle stated after Proposition~\ref{prop1}.

The seventh and eighth term in \eqref{Bonydelta_uI} can be estimated in the same manner and we will only estimate the seventh term as follows:
\begin{align*}
    &\sum_{\mu\ll\rho} 
\sum_{   \nu\sim\nu'\gtrsim\rho}  \mu^{-2} \rho \sup_\tau\left\|\|P_\mu L(P_\nu\delta u_3,P_{\nu'}\partial_xu_{\mathfrak T}^{(2)} \|_{L^\infty_z}\|\psi_{\sim\rho,1}^{(1)}\|_{L^2_z}
\|w_{\rho,1}^{\tau}\|_{L^2_z}\right\|_{L^1_{t,y}}\\
\lesssim&  \sum_{\mu \ll \rho} \sum_{\nu\sim \nu^\prime \gtrsim \rho }    \rho  \mu^{d-2} \nu^\prime  \left\|  P_\nu  \delta  u_3 \right\|_{L^2_{t,x}} \left\|  P_{\nu^\prime}   u_{\mathfrak{T}} \right\|_{L^2_{t,x}}   \cdot \left\|  \psi_{\sim \rho } \right\|_{L^\infty_t L^2_x}   \left\|    w_{\rho} \right\|_{L^\infty_t L^2_x}  \\
\lesssim &   \sum_{\mu \ll \rho} \sum_{\nu\sim \nu^\prime \gtrsim \rho }    \rho  \mu^{d-2} \nu^\prime D\nu^{-s} d_\nu (\nu^\prime)^{-s-1} c_{\nu^\prime} \rho^{-s} c_\rho  D \rho^{-s+1} d_\rho \\
\lesssim &\, \varepsilon^2 D^2 \rho^{-2s+2} d_\rho^2.
\end{align*}

For the term \(\mu^{-1}P_\mu\big(  \delta(\lambda\wedge\lambda)\big)\),  the Bony's decomposition gives
\begin{align}
 P_\mu\big(  \delta(\lambda\wedge\lambda)\big)=  &  L(P_{\ll \mu} \lambda^{(2)}, P_\mu \delta \lambda)+ L(P_{ \mu} \lambda^{(2)}, P_{\ll \mu} \delta \lambda) \label{bon.delta.lam} \\
  &+ \sum_{  \rho \gg \nu\sim\nu'\gtrsim\mu}  P_\mu( L(\lambda_\nu^{(2)}, \delta\lambda_{ \nu^\prime} )+ \sum_{   \nu\sim\nu'\gtrsim\rho}  P_\mu\delta
\bigl(
\lambda_\nu\wedge\lambda_{\nu'} \nonumber
\bigr).
\end{align}
 The first three terms can be estimated by \eqref{pre.uni.bi.es2.1}. For the
high--high contribution, it is enough to estimate the second quadratic term
in \eqref{d4.polarized.codazzi}; the first one is identical, while the cubic
terms are covered by Remark~\ref{Rem_J11t4}. Thus
\begin{align*}
     &\sum_{   \nu\sim\nu'\gtrsim\rho} 
   \mu^{-1}   \frac{\mu}{\nu}  \sup_\tau\left\|\|P_\mu L (P_\nu  \lambda^{(2)}, P_{ \nu^\prime} \delta \lambda)\|_{L^\infty_z}\|\psi_{\sim\rho,1}^{(1)}\|_{L^2_z}
\|w_{\rho,1}^{\tau}\|_{L^2_z}\right\|_{L^1_{t,y}}\\
\lesssim &\, \sum_{\nu\sim \nu^\prime   \gtrsim \rho}	  \rho   \mu^{d-3} \nu^{-1}  \| P_{ \nu } \lambda^{(2)} \|_{L^4_{t,y}L^2_z}     \|P_{\nu^\prime } \delta  \lambda  \|_{L^4_{t,y}L^2_z}   \left\|   \psi_{\sim \rho,1} \right\|_{L^4_{t,y}L^2_z}   \|        w_{\rho,1} \|_{ L^4_{t,y}L^2_z}    \nonumber\\
 	\lesssim&\,  \sum_{\nu^\prime \sim \nu  \gtrsim \rho}\sum_{\rho^\prime \sim  \rho} 	 \rho \mu^{d-3}  \nu^{-1} \cdot\nu^{-s+\frac 14}c_{\nu  }  D (\nu^\prime)^{-s+\frac 54}d_{\nu^\prime }   \cdot  (\rho^\prime)^{-s+\frac 14 }c_{\rho^\prime}  \cdot D\rho^{-s+\frac 54 }d_{\rho}  \nonumber\\
 	\lesssim&\, \varepsilon^2 D^2 \rho^{-2s +2 }  d_{\rho }^2,\nonumber
\end{align*}
where \(d \geq 5.\)
This completes the estimate for the term
\(\delta \mathbf{U} \partial_x\psi^{(1)}\).

Next, we consider the term \(\psi^{(2)}\delta B\) in \eqref{diff_N}. Bony's decomposition gives
\[
        P_\rho(\psi^{(2)}\delta B)
        =
       \psi^{(2)}_{\ll\rho}\delta B_\rho
        +
       \psi^{(2)}_\rho\delta B_{\ll\rho}
        +
        P_\rho
        \sum_{\nu\sim\nu'\gtrsim\rho}
        \psi^{(2)}_\nu\delta B_{\nu'} .
\]
 For the low--high interaction, H\"older's inequality and \eqref{pre.uni.bi.es1} yield
\begin{align}
       \Big | 
        \int_0^t\int
       \psi^{(2)}_{\ll \rho }\delta B_\rho
        \overline{w}_\rho\,\mathrm{d} x \mathrm{d} s \Big |
        &\lesssim
        \|\delta B_\rho\|_{L^2_{t,x}}
        \| \psi^{(2)}_{\ll \rho }  w_\rho \|_{L^2_{t,x}}                                     \notag\\
        &\lesssim
           \varepsilon    D^2\rho^{-2s+2} d^2_\rho
\end{align}
For the high--low interaction, H\"older's inequality and \eqref{pre.uni.deltB.es1} give
\begin{align}
        \Big|   
        \int_0^t\int
         \psi^{(2)}_\rho\delta B_{\ll\rho}
        \overline{w}_\rho\,\mathrm{d} x \mathrm{d} s \Big|
        &\lesssim
        \|\delta B_{\ll \rho} \|_{L^2_{t,y} L^\infty_z}
        \|\psi^{(2)}_\rho\|_{L^4_{t,y} L^2_z}
        \|w_\rho\|_{L^4_{t,y} L^2_z}                                      \notag\\
       &\lesssim   \varepsilon D^2\rho^{-2s+2}d_\rho^2 .                     
\end{align}
The high--high interaction can be estimated by the same dyadic summation argument used above:
\begin{align*}
       & \Big| \int_0^t\int
        \sum_{\nu\sim\nu'\gtrsim\rho}  P_\rho(
        \psi^{(2)}_\nu\delta B_{\nu'} )
        \overline{w}_\rho\,\mathrm{d} x \mathrm{d} s \Big|\nonumber\\
        &\lesssim
        \sum_{\nu\sim\nu'\gtrsim\rho}
        \|\delta B_{\nu'}\|_{L^2_{t,x}}
        \|\psi^{(2)}_\nu\|_{L^4_{t,y}L^2_z}
        \|w_\rho\|_{L^4_{t,y}L^\infty_z}                         \notag\\
                &\lesssim
        \varepsilon D^2\rho^{-2s+2}d_\rho^2.          
\end{align*}
This completes the estimate for the term
\(\psi^{(2)}\delta B\). Collecting the estimates obtained above, we  conclude that
\begin{align*}
     \eqref{UA1.term1} \lesssim \varepsilon D^2\rho^{-2s+2}d_\rho^2. 
\end{align*}

We now turn to \eqref{UA1.term2},   where \( \partial_\alpha
\big(\delta h^{\alpha\beta}\partial_\beta\psi^{(2)}\big)\) incurs an additional loss of one derivative and has no counterpart in the original equation \eqref{gaugsmcf}. Consequently, some parts of the argument require the weaker estimate \eqref{Uni.en2}, rather than the stronger pointwise bound \eqref{Uni.en1}.
Integration by parts gives
\[
\eqref{UA1.term2}=\left|
\int_0^t\int_{\mathbb R^d}
P_\rho
(\delta h^{\alpha\beta}\partial_\beta\psi^{(2)})\cdot
\partial_\alpha\overline{w}_\rho\,\mathrm{d} x \mathrm{d} s 
\right|.
\]
Bony's decomposition yields
\[
P_\rho (\delta h\,\partial_x\psi^{(2)})
=
 \delta h_{\ll\rho}\partial_x\psi_\rho^{(2)}
+
\delta h_\rho\partial_x\psi_{\ll\rho}^{(2)}
+
\sum_{\nu \sim \nu^\prime\gtrsim \rho} P_\rho L( 
\delta u_{3,\nu},\psi_{  \nu^\prime}^{(2)}),
\]
The high--low and high--high contributions can be estimated in the same manner as \eqref{est_J2} and \eqref{est_J3}, respectively, using the corresponding weak energy, elliptic, and interaction Morawetz assumptions. More precisely,
\begin{align*}
   &\,
\left|
\int_0^t\int_{\mathbb R^d}
\big(
\delta h_\rho^{\alpha\beta}
\partial_\beta\psi_{\ll\rho}^{(2)}
\big)
\partial_\alpha\overline{w}_\rho
\right|\\
\lesssim &\,    
\|\delta u_{3,\rho}\|_{L^2_{t,x}}
\|  \partial_\beta P_{\ll \rho} \psi^{(2)}\|_{L^4_{t,y} L^\infty_z}
\|w_\rho\|_{L^4_{t,y} L^2_z}
\\
\lesssim &\,  \varepsilon D^2 \rho^{-2s+2} d_\rho^2.
\end{align*}
and
\begin{align*}
  & \sum_{\nu \sim \nu^\prime\gtrsim\rho}
\left|
\int_0^t\int_{\mathbb R^d}
 P_\rho L( 
\delta u_{3,\nu},\psi_{  \nu^\prime}^{(2)})
 \partial_\alpha \overline{w}_\rho
\right| \\
\lesssim &   \sum_{\nu \sim \nu^\prime\gtrsim\rho} \rho \|\delta u_{3,\nu}\|_{L^2_{t,x}}
\| \psi_{  \nu^\prime}^{(2)}\|_{L^4_{t,y} L^2_z}
\|w_\rho\|_{L^4_{t,y} L^\infty_z} \\
\lesssim & \, \varepsilon D^2 \rho^{-2s+2} d_\rho^2.
\end{align*}

It remains to treat the low--high contribution.  Decomposing the low frequency coefficient dyadically gives
\begin{align}
    &\rho 
\int_0^t\int
L
(\delta h_{\ll\rho}^{\alpha\beta}, \psi_\rho^{(2)})
\partial_\alpha\overline{w}_\rho\,\mathrm{d} x \mathrm{d} s \nonumber\\
   = &\,  \rho \sum_{\mu \ll \rho}
\int_0^t\int
L
(\delta h_{\mu}^{\alpha\beta}, \psi_\rho^{(2)})
\partial_\alpha\overline{w}_\rho\,\mathrm{d} x \mathrm{d} s.  \label{UA1.term2.t1}
\end{align}
As in the low--high analysis of \eqref{bony.diff.uI}, we use the elliptic
equation for \(\delta h\) to recover the apparently lost derivative.
Taking the difference of the schematic metric equation \eqref{ellis2}, we
obtain
\[
\begin{aligned}
\Delta\delta h
&=
h^{(1)}\partial_x^2\delta h
+
\delta h \partial_x^2h^{(2)}
+
(u_3^{(1)}+u_3^{(2)})\delta u_3
+
\delta  (\lambda\wedge\lambda )
 + \delta ( h\cdot\partial_x h\cdot \partial_x h)
\end{aligned}
\]
Applying \(P_\mu\Delta^{-1}\), we get
\begin{align}
   \delta h_\mu
=&\,
L(h^{(1)}_{\ll\mu},\delta h_\mu)
+L(\delta h_{\ll \mu},  h^{(2)}_{\mu})+\mu^{-1} L( u_{3,\ll \mu}^{(a)}, \delta h_{\mu} ) \label{Pmuh}
\\
&\,+\mu^{-2}L(h^{(1)}_{\mu},P_{\ll\mu}\partial_x^2 \delta h)
+\mu^{-2} L(\delta h_{  \mu},  \partial_x^2 h^{(2)}_{\ll\mu}) +\mu^{-1} L( h_{ \mu}^{(a)}, \delta u_{3,\ll\mu} )\nonumber\\
&+\mu^{-2}\sum_{\mu\lesssim \nu\sim \nu^\prime \ll \rho} P_\mu L(u^{(a)}_{3,\nu}, \delta u_{3,\nu^\prime}) +
\mu^{-2}\sum_{\nu\sim \nu^\prime \gtrsim \rho} P_\mu L(u^{(a)}_{3,\nu}, \delta u_{3,\nu^\prime})\nonumber\\
&+ 
\delta  (\lambda\wedge\lambda )
+
 \delta ( h\cdot\partial_x h\cdot \partial_x h), \nonumber
\end{align}
where \(a=1,2\).
The contribution of \(L(h^{(1)}_{\ll\mu},P_\mu\delta h)\) can be absorbed, since
\(
\|h^{(1)}\|_{L^\infty}\lesssim \varepsilon. 
\) 
The cubic terms can be handled by the general argument in Remark~\ref{Rem_J11t4} and will therefore be omitted.

The second through sixth terms on the right-hand side of \eqref{Pmuh} are
handled similarly by \eqref{pre.uni.bi.es1+1}--\eqref{pre.uni.bi.es2.1}.
We illustrate the argument with the third term. Substituting it into
\eqref{UA1.term2.t1} and applying \eqref{pre.uni.bi.es1+2}, we obtain
\begin{align}
  &\,    \rho \mu^{-1} \left| \sum_{\mu \ll \rho}
\int_0^t\int
L
(u^{(a)}_{3, \ll\mu}, \delta h_{ \mu}, \psi_\rho^{(2)})
\partial_\alpha\overline{w}_\rho\,\mathrm{d} x \mathrm{d} s  
\right|  \label{Mrho.term2.1} \\
\lesssim &  \sum_{\mu\ll\rho}   
\rho^2\mu^{-1}  \big \| \|   L(  \delta h_{\mu},    u^{(a)}_{3, \ll \mu})\|_{L^\infty_{\widehat{x}_1}} \| w_{\rho,1}^{\tau}\|_{L^2_{\widehat{x}_1}} \|  \psi_{\sim\rho,1}^{(2)}\|_{L^2_{\widehat{x}_1}} \big\|_{L^2_{ t,x_1}} \nonumber\\
\lesssim &\,  D^2\rho^{-2s+2}
c_\rho d_\rho
\sum_{\mu\ll\rho}c_\mu d_\mu, \nonumber
\end{align}
which is of the form \eqref{Uni.en2}. Similarly, the second term on the right-hand side of \eqref{Pmuh} is
controlled by \eqref{pre.uni.bi.es1+1}:
\begin{align*}
  &\,    \rho  \left| \sum_{\mu \ll \rho}
\int_0^t\int
L
(\delta h_{\ll \mu}, h^{(2)}_{ \mu}, \psi_\rho^{(2)})
\partial_\alpha\overline{w}_\rho\,\mathrm{d} x \mathrm{d} s  
\right|  \\
\lesssim &\,  D^2\rho^{-2s+2}
c_\rho d_\rho
\sum_{\mu\ll\rho}c_\mu d_\mu.
\end{align*}
The seventh term on the right-hand side of \eqref{Pmuh} is estimated as
follows: 
\begin{align*}
     &\sum_{\mu\ll\rho}   \sum_{\mu \lesssim\nu\sim\nu'\ll \rho}   \rho \mu^{-2}\left|  
\int_0^t\int
L
(P_\mu(
  u_{3,\nu}^{(a)},  \delta u_{3, \nu'}), \psi_\rho^{(2)})
\partial_\alpha\overline{w}_\rho\,\mathrm{d} x \mathrm{d} s 
\right|   \\
\lesssim &\,  \sum_{\mu\ll\rho}   \sum_{\mu \lesssim\nu\sim\nu'\ll \rho}   \rho \mu^{-2} \sup_{\tau,\tau'} \big \| \|   P_\mu(
  u_{3,\nu}^{(a)},  \delta u_{3, \nu'})\|_{L^\infty_{\widehat{x}_1}} \| w_{\rho,1}^{\tau'}\|_{L^2_{\widehat{x}_1}} \|  (\psi_{\sim\rho,1}^{(2)})^{\tau} \|_{L^2_{\widehat{x}_1}} \big\|_{L^2_{ t,x_1}}\\
  \lesssim &\,  \sum_{\mu\ll\rho}   \sum_{\mu \lesssim\nu\sim\nu'\ll \rho}   \rho \mu^{d-3} \sup_{\tau,\tau'} \big \| \|   P_\mu(
  u_{3,\nu}^{(a)},  \delta u_{3, \nu'})\|_{L^1_{\widehat{x}_1}} \| w_{\rho,1}^{\tau'}\|_{L^2_{\widehat{x}_1}} \|  (\psi_{\sim\rho,1}^{(2)})^{\tau} \|_{L^2_{\widehat{x}_1}} \big\|_{L^2_{ t,x_1}}\\
  \lesssim  &\,  \sum_{\mu\ll\rho}   \sum_{\mu \lesssim\nu\sim\nu'\ll \rho}   \rho \mu^{d-3} \nu^{-s} c_\nu D \rho^{-s+\frac 12} d_\rho \cdot (\nu')^{-s+1} D d_{\nu'} \rho^{-s-\frac 12} c_\rho \\
\lesssim &\,  D^2\rho^{-2s+2}
c_\rho d_\rho
\sum_{\nu\ll\rho}c_\nu d_\nu.
\end{align*}

For the eighth term on the right-hand side of \eqref{Pmuh}, we have
\begin{align}
     &\,   \sum_{\mu\ll\rho}  \sum_{\nu\sim\nu'\gtrsim\rho}   \rho \mu^{-2}  \left|  
\int_0^t\int
L
(P_\mu(
  u_{3,\nu}^{(a)},  \delta u_{3, \nu'}), \psi_\rho^{(2)})
\partial_\alpha\overline{w}_\rho\,\mathrm{d} x \mathrm{d} s 
\right|   \label{Mrho.term2.3} \\
\lesssim & \sum_{\mu\ll\rho}   \rho \mu^{-2} \sum_{\nu\sim\nu'\gtrsim\rho}   \| 
P_\mu(
  u_{3,\nu}^{(a)},  \delta u_{3, \nu'})\|_{L^\infty_x L^1_{t}}   \|\psi_\rho^{(2)}\|_{L^2_{x}} \|\partial_\alpha\overline{w}_\rho\|_{L^2_{x}} \nonumber\\
\lesssim &  \sum_{\mu\ll\rho}  \rho^2 \mu^{d-2}    \sum_{\nu\sim\nu'\gtrsim\rho}    \| 
\delta u_{3,\nu} \|_{L^2_{t,x}}   
 \| 
  u^{(2)}_{3,\nu'}\|_{L^2_{t,x}}   \|\psi_\rho^{(2)}\|_{L^2_{x}} \| \overline{w}_\rho\|_{L^2_{x}} \nonumber\\
\lesssim & \, \varepsilon^2 D^2 \rho^{-2s+2} d_\rho^2. \nonumber
\end{align}

It remains to estimate
\(\delta(\lambda\wedge\lambda)\) in \eqref{Pmuh}, whose Bony
decomposition is given in \eqref{bon.delta.lam}. When substituted into
\eqref{UA1.term2.t1}, the first three terms in
\eqref{bon.delta.lam} are controlled by \eqref{pre.uni.bi.es2.1}, exactly
as for the second through seventh terms on the right-hand side of
\eqref{Pmuh}.

For the fourth term in \eqref{bon.delta.lam}, we exploit the Gauss
structure in \eqref{d4.polarized.codazzi}. It suffices to estimate the
second quadratic contribution, since the first is identical:
\begin{align}
     &\,   \sum_{\mu\ll\rho} \sum_{\nu\sim\nu'\gtrsim\rho} \rho \mu^{-2} \frac{\mu}{\nu}  \left|  
\int_0^t\int
L
(P_\mu(
 \lambda_\nu, \delta\lambda_{ \nu^\prime }), \psi_\rho^{(2)})
\partial_\alpha\overline{w}_\rho\,\mathrm{d} x \mathrm{d} s  
\right|  \label{uni.critic.d=5} \\
\lesssim &\,  \sum_{\mu\ll\rho} \sum_{\nu\sim\nu'\gtrsim\rho} \rho^2 \mu^{-1} \nu^{-1}  \|P_\mu(
 \lambda_\nu, \delta\lambda_{ \nu^\prime })\|_{L^2_t L^\infty_x} \left\|   \psi_{  \rho} \right\|_{L^4_{t,y}L^2_z}   \|        w_{\rho}  \|_{ L^4_{t,y}L^2_z}    \nonumber
 \\
\lesssim &  \sum_{\mu\ll\rho}  \rho^2 \mu^{ d-4} \nu^{-1} \sum_{\nu\sim\nu'\gtrsim\mu}   \| \lambda_\nu \|_{L^4_{t,y}L^2_z}     \|\delta\lambda_{ \nu^\prime }  \|_{L^4_{t,y}L^2_z}   \left\|   \psi_{  \rho} \right\|_{L^4_{t,y}L^2_z}   \|        w_{\rho}  \|_{ L^4_{t,y}L^2_z}   \nonumber\\
\lesssim &  \sum_{\mu\ll\rho}     \sum_{\nu\sim\nu'\gtrsim\rho}     \rho^2 \mu^{d-4}  \nu^{-s-\frac 34} c_\nu \cdot D (\nu^\prime)^{-s+\frac 54} d_{\nu^\prime} \cdot \rho^{-s+\frac 14} c_\rho \cdot D \rho^{-s+\frac 54} d_\rho    \nonumber\\
\lesssim & \, \varepsilon^2 D^2 \rho^{-2s+2} d_\rho^2, \nonumber
\end{align}
provided that \( d \geq 5\). 

 Finally, combining the estimates for
\eqref{UA1.term1}--\eqref{UA1.term3}, we obtain \eqref{Uni.en}.

\subsection{Difference interaction Morawetz and bilinear estimates}
The difference interaction Morawetz estimate
\[
\|w_\rho\|_{L_{t,y}^4L_z^2(I)}
\lesssim
\rho^{-s+\frac54}d_\rho
\bigl(D(0)+\varepsilon D(I)\bigr)
\]
follows by adapting the argument of Section~\ref{secInter_ac} to the
frequency localized difference equation
\begin{align*}
\sqrt{-1}\partial_t w_\rho
+\partial_\alpha
\bigl(g_{\ll\rho}^{\alpha\beta,(1)}
\partial_\beta w_\rho\bigr)
&=
\mathfrak{N}_\rho
-\partial_\alpha
\bigl(g_{\gtrsim\rho}^{\alpha\beta,(1)}
\partial_\beta w_\rho\bigr)
:=\mathcal{I}_\rho^{h,(1)}.
\end{align*}
Indeed, repeating the argument of Appendix~\ref{appMora},  with
\(
\psi_\rho,\, g,\, h,\, \mathcal{I}_\rho^h
\)
in \eqref{intaet0}--\eqref{intaet6} replaced respectively by
\(
w_\rho, \,g^{(1)}, \, h^{(1)},\,
\mathcal{I}_\rho^{h,(1)},
\) yields the corresponding interaction Morawetz inequality for \(w_\rho\). The resulting terms are controlled by the estimates
established in Section~\ref{secInter_ac} and
Subsection~\ref{subsec:diff.en}, which yields the desired bound. We omit
the details.

 We next turn to the difference bilinear estimates. We shall establish the following two types of estimates:
\begin{align}
&\sup_\tau
\left\|
\|w_\mu \|_{L_{\widehat x_\gamma}^2}
\| \psi^{(a),\tau}_{\rho,\gamma}\|_{L_{\widehat x_\gamma}^2}
\right\|_{L_{t,x_\gamma}^2(I)}
\nonumber\\
&\qquad\lesssim
\mu^{-s+1}d_\mu
\rho^{-s-\frac12}c_\rho
\left(
D(0)+\varepsilon D(I)
\right),
\label{weak-bilinear-one}
\end{align}
and
\begin{align}
&\sup_\tau
\left\|
\|\psi_\mu^{(a)}\|_{L_{\widehat x_\gamma}^2}
\|  w^\tau_{\rho,\gamma}\|_{L_{\widehat x_\gamma}^2}
\right\|_{L_{t,x_\gamma}^2(I)}
\nonumber\\
&\qquad\lesssim
\mu^{-s}c_\mu
\rho^{-s+\frac12}d_\rho
\left(
D(0)+\varepsilon D(I)
\right),
\label{weak-bilinear-two}
\end{align}
 where \( a=1,2; \gamma=1,\cdots,d.\)  As in the proof of the a priori bilinear estimates in
Section~\ref{sectBi}, our strategy is to reduce the nonlinear error terms
arising from the div--curl estimates to those already controlled by
the difference energy estimates in Subsection~\ref{subsec:diff.en} or by
the corresponding a priori estimates in Section~\ref{sectBi}.

We first prove \eqref{weak-bilinear-one}. As in \eqref{phiNeq}, the
frequency localized equation for \(\psi_{\rho,1}^{(a)}\) is
\begin{align*} 
	&\sqrt{-1} \pa_t \psi^{(a)}_{\rho,1} +\sum_{\alpha,\beta=1}^d \partial_{\alpha} (g^{\alpha \beta} \partial_{ \beta}  \psi^{(a)}_{\rho,1} )=P_{\rho,1} \mathcal{N}^{(a)}- \partial_{\alpha}  ( [P_{\rho,1}, h^{\alpha \beta} ]  \partial_\beta  \psi^{(a)} ):= \mathcal {G}^{(a)}_\rho,	 
\end{align*}
Replacing \(\rho\) by \(\mu\) in \eqref{eq:diff}, we similarly obtain
 \begin{align*} 
	  \sqrt{-1}\partial_t w_{\mu}
+\partial_\alpha\big(g^{\alpha\beta}_{1 }\partial_\beta w_{\mu}\big)
 = \mathfrak{N}_\mu,
\end{align*}
%The momentum balance law for \(\psi_{\rho,1}^{(a)}\) and the mass balance
%law for \(w_\mu\) are, respectively,
%\begin{equation}  \begin{aligned}
%		&-\partial_t \int_{\mathbb{R}^{d-1}}  \mathrm{Im}\left( \psi^{(a)}_{\rho,1 } \partial_{1} \overline{\psi}^{(a)}_{\rho,1 }\right) \mathrm{d} \widehat{x}_1 \\
%	&+\partial_{1}\Big(    \int_{\mathbb{R}^{d-1}} 2\mathrm{Re} (\partial_{1} \psi^{(a)}_{\rho,1 }\,   g^{1\beta} \partial_{\beta} \overline{\psi}^{(a)}_{\rho})-\partial_{1 }\mathrm{Re} ( \psi^{(a)}_{\rho,1 }\,    g^{1\beta} \partial_{\beta} \overline{\psi}^{(a)}_{\rho} ) \mathrm{d} \widehat{x}_1 \Big)\\
%	=&\, -\int_{\mathbb{R}^{d-1}}  \mathrm{Re} (\partial_{\alpha} \psi^{(a)}_{\rho,1 }   \partial_1 h^{\alpha \beta} \partial_\beta \overline{\psi}^{(a)}_{\rho,1} )  \mathrm{d} \widehat{x}_1  \\
%	&\,  + \int_{\mathbb{R}^{d-1}}  \mathrm{Re} (\mathcal{G}_\rho^{(a)} \partial_{1} \overline{\psi}^{(a)}_\rho) -\mathrm{Re} ( \overline{\psi}^{(a)}_{\rho}\partial_{1} \mathcal{G}_\rho^{(a)} )\mathrm{d} \widehat{x}_1 ,     
%\end{aligned}
%\end{equation}
%and
%\begin{equation} 
%\begin{aligned}
%	&\frac 12\partial_t \int_{\mathbb{R}^{d-1}}   \left|P_\mu w\right|^2  	\mathrm{d} \widehat{x}_1+  \partial_{1} \int_{\mathbb{R}^{d-1}} \mathrm{Im}(   \overline{ w}_\mu g^{1\beta}   \partial_{\beta} w_\mu  ) \mathrm{d} \widehat{x}_1  \\
%	=&\int_{\mathbb{R}^{d-1}} \mathrm{Im} (  \mathfrak{N}_\mu  \overline{ w_\mu }) \mathrm{d} \widehat{x}_1.
%\end{aligned}
%\end{equation}
As in Section \ref{sectBi},  applying Lemma~\ref{lemdiv-curl} to  the  momentum balance law for \(\psi_{\rho,1}^{(a)}\) and the mass balance
law for \(w_\mu\), we obtain
\begin{align}
& \rho^2 \int_{[0,T]\times \mathbb{R}} \int_{\mathbb{R}^{d-1}} |  w_\mu   |^2 	\mathrm{d} \widehat{x}_1 \int_{\mathbb{R}^{d-1}} |  \psi_{\rho,1}^{(a)} |^2 \mathrm{d} \widehat{x}_1 \mathrm{ d} x_1 \mathrm{ d} t \\
\lesssim&\int_{[0,T]\times \mathbb{R}} \int_{\mathbb{R}^{d-1}} |  P_\mu  w |^2 	\mathrm{d} \widehat{x}_1 \int_{\mathbb{R}^{d-1}} \mathrm{Re} (\partial_{1} \psi_\rho^{(a)}   g^{1\beta} \partial_{\beta} \overline{\psi}_{\rho}^{(a)} ) \mathrm{d} \widehat{x}_1 \mathrm{ d} x_1 \mathrm{ d} t \nonumber\\
&+\int_{[0,T]\times \mathbb{R}} \int_{\mathbb{R}^{d-1}} \mathrm{Re} (    P_\mu w \partial_{1} P_\mu  \overline{w}   )	\mathrm{d} \widehat{x}_1 \int_{\mathbb{R}^{d-1}} \mathrm{Re} (\psi_\rho^{(a)} \,   g^{1\beta} \partial_{\beta} \overline{\psi}_{\rho}^{(a)} ) \mathrm{d} \widehat{x}_1 \mathrm{ d} x_1 \mathrm{ d} t\nonumber\\
&-\int_{[0,T]\times \mathbb{R}}  \int_{\mathbb{R}^{d-1}}  \mathrm{Im} ( \psi_{\rho }^{(a)} \partial_{1} \overline{\psi}_\rho^{(a)} ) \mathrm{d} \widehat{x}_1 \int_{\mathbb{R}^{d-1}}\mathrm{Im} ( P_\mu \overline{w} g^{1\beta}    \partial_{\beta}P_\mu  w   ) \mathrm{d} \widehat{x}_1 \mathrm{ d} x_1 \mathrm{ d} t \nonumber \\
\lesssim &\,\rho \max_{0\leq t \leq T} \| \psi_\rho^{(a)} (t, \cdot) \|_{  L_x^2}^2  \max_{0\leq t \leq T} \|    w_\mu(t,\cdot )\|_{ L_x^2}^2 \label{diff.1biNt}\\
	&\,+\Big|\int_{[0,T]\times \mathbb{R}}  \left( \int_{-\infty}^{x_1} \int_{\mathbb{R}^{d-1}}  \mathrm{Im}\left( \psi_{\rho,1 }^{(a)} \partial_{1} \overline{\psi}_{\rho,1}^{(a)}\right) 
 \mathrm{d} \widehat{x}_1 \mathrm{ d} y\right) \label{non.uni.bit1}\\
 &\qquad\qquad\qquad  \cdot\int_{\mathbb{R}^{d-1}} \mathrm{Im} \left( \overline{ w}_\mu  \mathfrak{N}_\mu    \right) \mathrm{d} \widehat{x}_1 \mathrm{ d} x_1 \mathrm{ d} t\Big| \nonumber\\
	&\,+\Big|\int_{[0,T]\times \mathbb{R}}  \left( \int^{+\infty}_{x_1}   \left\|w_\mu\right\|^2_{L^2_{\widehat{x}_1}}   \mathrm{ d} y\right )\label{non.uni.bit2}\\
 &\qquad\qquad\qquad  \cdot\int_{\mathbb{R}^{d-1}}  \mathrm{Re}\left(\partial_{\alpha} \psi_{\rho,1}^{(a)}   \partial_1 h^{\alpha \beta} \partial_\beta \overline{\psi}^{(a)}_{\rho,1}\right)  \mathrm{d} \widehat{x}_1\mathrm{ d} x_1 \mathrm{ d} t \Big|\nonumber\\
	&\,+\Big|\int_{[0,T]\times \mathbb{R}}  \left( \int^{+\infty}_{x_1}   \left\|w_\mu\right\|^2_{L^2_{\widehat{x}_1}}   \mathrm{ d} y\right) 
 \Big|\label{non.uni.bit3}\\
	&\qquad\qquad\qquad  \cdot \int_{\mathbb{R}^{d-1}} \mathrm{Re} \left(\overline{\psi}_{\rho,1}^{(a)} \partial_{1} \mathcal{G}_\rho^{(a)} - \mathcal{G}_\rho^{(a)} \partial_{1} \overline{\psi}_{\rho,1}^{(a)}\right)   \mathrm{d} \widehat{x}_1 \mathrm{ d} x_1 \mathrm{ d} t \Big| \nonumber\\
    \lesssim&\, \varepsilon^2 \rho^{-2s+1} \mu^{-2s+2} c_\rho^2 d_\mu^2 ( D^2(0)+ D^2(I)).\nonumber
\end{align}
Indeed, for the first nonlinear contribution, we have
\begin{align*}
    \eqref{non.uni.bit1} \lesssim  &\, \rho \left\|\psi_{\rho } \right\|_{L^\infty_t L^2_x }^2  \left\|\mathrm{Im} \left(  \mathfrak{N}_\mu \overline{ w}_\mu\right) \right\|_{L^1_{t,x}},
\end{align*}
where the estimates
\[
\left\|
\operatorname{Im}
\left(
\mathfrak N_\mu\overline{w_\mu}
\right)
\right\|_{L_{t,x}^1} \lesssim \varepsilon  D^2 \mu^{-2s+2} d_\mu^2
\]
has already been obtained in Subsection~\ref{subsec:diff.en}. The
nonlinear contributions in \eqref{non.uni.bit2} and
\eqref{non.uni.bit3} can be estimated by adapting the corresponding
arguments from the proof of the a priori bilinear estimates in
Section~\ref{sectBi}. This proves \eqref{weak-bilinear-one}.

Similarly, for \eqref{weak-bilinear-two}, we have 
\begin{align*} 
	&\sqrt{-1} \pa_t \psi^{(a)}_\mu +\sum_{\alpha,\beta=1}^d \partial_{\alpha} (g^{\alpha \beta} \partial_{ \beta}  \psi^{(a)}_\mu )=P_\mu \mathcal{N}- \partial_{\alpha}  ( [P_\mu, h^{\alpha \beta} ]  \partial_\beta  \psi^{(a)} ):= \mathcal{I}_\mu^{(a)} 
\end{align*}
and
 \begin{align*} 
	  \sqrt{-1}\partial_t w_{\rho,1}
+\partial_\alpha\big(g^{\alpha\beta}_{1 }\partial_\beta w_{\rho,1}\big)
 = \mathfrak{N}_\rho,  
\end{align*}
where \(\mathfrak{N}_\rho\) is given in \eqref{eq:diff}.
Thus the same process as those in Section \ref{sectBi}, we have
\begin{align}
& \rho^2 \int_{[0,T]\times \mathbb{R}} \int_{\mathbb{R}^{d-1}} |  \psi_\mu^{(a)}   |^2 	\mathrm{d} \widehat{x}_1 \int_{\mathbb{R}^{d-1}} |  w_{\rho,1} |^2 \mathrm{d} \widehat{x}_1 \mathrm{ d} x_1 \mathrm{ d} t \nonumber \\
\lesssim &\,\rho \max_{0\leq t \leq T} \| \psi_\mu^{(a)} (t, \cdot) \|_{  L_x^2}^2  \max_{0\leq t \leq T} \|    w_\rho(t,\cdot )\|_{ L_x^2}^2 \label{1biNt0}\\
	&\,+\Big|\int_{[0,T]\times \mathbb{R}}  \left( \int_{-\infty}^{x_1} \int_{\mathbb{R}^{d-1}}  \mathrm{Im}\left( w_{\rho,1 }  \partial_{1} \overline{w}_{\rho,1} \right) 
 \mathrm{d} \widehat{x}_1 \mathrm{ d} y\right) \label{1non.uni.bit1}\\
 &\qquad\qquad\qquad  \cdot\int_{\mathbb{R}^{d-1}} \mathrm{Im} \left( \overline{ \psi}_\mu  \mathcal{I}_\mu^{(a)}    \right) \mathrm{d} \widehat{x}_1 \mathrm{ d} x_1 \mathrm{ d} t\Big| \nonumber\\
	&\,+\Big|\int_{[0,T]\times \mathbb{R}}  \left( \int^{+\infty}_{x_1}   \big\|\psi_\mu^{(a)}\big\|^2_{L^2_{\widehat{x}_1}}   \mathrm{ d} y\right )\label{1non.uni.bit2}\\
 &\qquad\qquad\qquad  \cdot\int_{\mathbb{R}^{d-1}}  \mathrm{Re}\left(\partial_{\alpha} w_{\rho,1}    \partial_1 h^{\alpha \beta} \partial_\beta \overline{w}_{\rho,1}\right)  \mathrm{d} \widehat{x}_1\mathrm{ d} x_1 \mathrm{ d} t \Big|\nonumber\\
	&\,+\Big|\int_{[0,T]\times \mathbb{R}}  \left( \int^{+\infty}_{x_1}   \big\|\psi_\mu^{(a)}\big\|^2_{L^2_{\widehat{x}_1}}   \mathrm{ d} y\right) 
 \Big|\label{1non.uni.bit3}\\
	&\qquad\qquad\qquad  \cdot \int_{\mathbb{R}^{d-1}} \mathrm{Re} \left(\overline{w}_{\rho,1}  \partial_{1} \mathfrak{N}_\rho  - \mathfrak{N}_\rho  \partial_{1} \overline{w}_{\rho,1} \right)   \mathrm{d} \widehat{x}_1 \mathrm{ d} x_1 \mathrm{ d} t \Big|. \nonumber
\end{align}
The \eqref{1non.uni.bit1} and \eqref{1non.uni.bit2} can be estimated similar to \eqref{non.bit1} and \eqref{non.bit2} respectively. While by \eqref{diff_N}, we can use the same analysis for \eqref{momen_B_3} to obtain
\begin{align}
    \eqref{1non.uni.bit3} \lesssim &\,     \big\|\psi_\mu^{(a)}\big\|^2_{L^\infty_t L^2_{x}}    \cdot \Big(  \rho\,  \| L( w_{\rho,1}, P_\rho\delta \mathcal{N} -wB^{(1)}+ [P_{\rho,1} , B^{(1)}]w )\|_{L^1_{t,x}} \nonumber \\
 & + \|  \pa_x w_{\rho,1} \cdot P_\rho\partial_\alpha
\big(\delta h^{\alpha\beta}\partial_\beta\psi^{(2)}\big) \|_{L^1_{t,x}}  \label{1non.uni.bit3.t1}\\
 &+ \|  \pa_x w_{\rho,1} \cdot \partial_\alpha\big([P_\rho,h^{\alpha\beta,(1)}]\partial_\beta w\big)\|_{L^1_{t,x}}\Big) \nonumber\\
    & \label{1non.uni.bit3.t2}\, +\Big| \int_{[0,T]\times \mathbb{R}}  \left( \int^{+\infty}_{x_1}   \big\|\psi_\mu^{(a)}\big\|^2_{L^2_{\widehat{x}_1}}   \mathrm{ d} y\right) \\
  &\qquad \qquad\qquad \cdot  \int_{\mathbb{R}^{d-1}}  \pa_{1}  \big(w_\rho \cdot  ([P_{\rho,1}, B^{(1)}]w)\big) -\pa_{1}B^{(1)} \cdot |w_{\rho,1} |^2 \mathrm{d} \widehat{x}_1 \mathrm{ d} x_1 \, \mathrm{ d} t\Big|. \nonumber
\end{align}
The  \eqref{1non.uni.bit3.t1} can be estimated as those in Subsection \ref{subsec_non_est1}, while \eqref{1non.uni.bit3.t2} can be estimated similar to \eqref{non.bit2_4}. This proves \eqref{weak-bilinear-two}.

\subsection{Difference Elliptic and \texorpdfstring{\(\delta B\)}{delta B} Estimates} \label{subsec:uniB}
Unlike the difference Schr\"odinger equation \eqref{eq:diff}, the difference elliptic equations and the difference \(B\)-equation, derived respectively from \eqref{Gen_u2}, \eqref{Gen_ellev}, and \eqref{eq:B-schematic}, contain no derivative-losing term analogous to \( P_\rho\partial_\alpha \bigl( \delta h^{\alpha\beta}\partial_\beta\psi^{(2)} \bigr). \) Consequently, the estimates for the difference elliptic variables \(\delta u_i\), \(i=2,3,4,5\), and for \(\delta B\) follow by adapting the corresponding arguments in Sections~\ref{sec_ell} and~\ref{sec_B}. We therefore omit the details.

Combining the estimates established in Subsections~\ref{subsec:diff.en}--\ref{subsec:uniB}, we complete the proof of Proposition~\ref{prop1}.

\subsection{Proof of Theorem \ref{thm1}}
Let \(\mathcal J_{1/n}\) be a standard smooth Fourier multiplier and define
\(
\psi_0^n
=
\mathcal J_{1/n}\psi_0.
\)
Then
\(
\|\psi_0^n\|_{\dot B^s_{2,1}}
\lesssim
\|\psi_0\|_{\dot B^s_{2,1}},
\)
and \(
\psi_0^n
\longrightarrow
\psi_0
~~
\text{in }\dot B^s_{2,1}.
\)
Moreover, for every integer \(m\geq0\),
\[
\|\psi_0^n\|_{\dot B^{s+m}_{2,1}}
\lesssim
n^m\|\psi_0\|_{\dot B^s_{2,1}}.
\]

For each \(n\), let
\(
\mathfrak S_0^n:
=
\mathfrak S(\psi_0^n)
=
\left(
\lambda_0^n,h_0^n,V_0^n,A_0^n,B_0^n
\right)
\)
be the solution of the elliptic constraint system in the harmonic/Coulomb
gauge. Then
\(
(\psi_0^n,\mathfrak S_0^n)
\)
is a smooth compatible datum. By the high regularity local well-posedness
theory in Huang--Tataru \cite{HT22}, there exists a maximal high regularity solution
\(
\left(
\psi^n,
\lambda^n,
h^n,
V^n,
A^n,
B^n
\right)
\)
on a maximal interval \(
[0,T_n^*).
\)
We claim that
\(
T_n^*=\infty.
\)
Indeed, the Theorem~\ref{thm_a_apr}  gives
\begin{align}
   \sup_{t<T_n^*}
\sum_{i=1}^5
\|u_i^n(t)\|_{\dot B^s_{2,1}}
& \lesssim
\varepsilon, \nonumber\\
\sup_{t<T_n^*}
\|u_i^n(t)\|_{\dot B^{s+2}_{2,1}}
&=
\sup_{t<T_n^*}
\sum_\rho
\rho^{s+2}
\|P_\rho u_i^n(t)\|_{L_x^2}
\nonumber\\
&\lesssim
\|\psi_0^n\|_{\dot B^{s+2}_{2,1}} \nonumber\\
&\lesssim  
n^2\|\psi_0\|_{\dot B^s_{2,1}}
=:C_n
<\infty.
\label{eq:approximate-high regularity}
\end{align}
The constant \(C_n\) may depend on \(n\), but it is uniform in
\(t<T_n^*\). Furthermore, since
\(
\|h^n\|_{L_{t,x}^\infty([0,T_n^*))}
\lesssim
\varepsilon,
\)
the metric
\(
g_n=I+h^n
\)
remains uniformly elliptic. The continuation criterion for high regularity
solutions therefore excludes a finite maximal time. Hence
\(
T_n^*=\infty.
\)

We have constructed a sequence of global smooth compatible solutions
satisfying
\[
\sup_n
\sup_{t\geq0}
\sum_{i=1}^5
\|u_i^n(t)\|_{\dot B^s_{2,1}}
\lesssim
\varepsilon,
\]
together with the uniform energy,  bilinear, interaction Morawetz, elliptic and
\(B\) estimates of Theorem~\ref{thm_a_apr}.

Fix \(T>0\). For \(n,m\geq1\), set
\(
w^{n,m}
=
\psi^n-\psi^m.
\)
Applying Proposition~\ref{prop1} on \([0,T]\), we
obtain
\[
\|\psi^n-\psi^m\|_{\mathcal X_T^{s-1}}
\lesssim
\|\psi_0^n-\psi_0^m\|_{\dot B^{s-1}_{2,1}},
\]
where \(\mathcal X_T^{s-1}\) denotes the collection of difference weak energy,
bilinear, interaction Morawetz,   elliptic and \(B\) norms occurring in
\(D([0,T])\).

It is evident that the right-hand side tends to zero. Indeed, since
\(\widehat{\mathcal J}\) is smooth and equals one at the origin, on the
support of \(P_\rho\) we have
\begin{align*}
 \|\psi_0^n-\psi_0^m\|_{\dot B^{s-1}_{2,1}}
\lesssim &\,
\left(
\frac1n+\frac1m
\right)
\|\psi_0\|_{\dot B^s_{2,1}}
\end{align*}
and hence \(\{\psi^n\}\) is Cauchy in \(\mathcal X_T^{s-1}\).

There exists a limit \(\psi\) such that
\(
\psi^n
\longrightarrow
\psi
~~
\text{in}\mathcal X_T^{s-1}.
\)
The difference elliptic estimates imply, at the same time, that
\[
\lambda^n\longrightarrow\lambda,
\quad
\partial_x h^n\longrightarrow\partial_x h,
\quad
V^n\longrightarrow V,
\quad
A^n\longrightarrow A
\]
in the corresponding weak topologies, and
\(
B^n\longrightarrow B
\)
in the weak dyadic \(B\) topology.

We next show that the convergence is strong in
\(
C([0,T];\dot B^s_{2,1}).
\)
Define
\(
c_\rho
=
\sum_\nu
2^{-\delta|\log_2(\rho/\nu)|} \nu^s\|P_\nu\psi_0\|_{L_x^2},
\)
where \(\delta>0\) is sufficiently small. Then
\[
\sum_\rho c_\rho
\lesssim
\|\psi_0\|_{\dot B^s_{2,1}}.
\]
Since \(\mathcal J_{1/n}\) is a uniformly bounded Fourier multiplier, the
same envelope controls every approximate datum:
\[
\rho^s\|P_\rho\psi_0^n\|_{L_x^2}
\lesssim
c_\rho.
\]
The frequency envelope version of Theorem~\ref{thm_a_apr} propagates this
bound, so that
\[
\sup_n
\rho^s
\|P_\rho\psi^n\|_{L_t^\infty L_x^2([0,T])}
\lesssim
c_\rho.
\]
Passing to the weak limit at each fixed frequency gives
\[
\rho^s
\|P_\rho\psi\|_{L_t^\infty L_x^2([0,T])}
\lesssim
c_\rho.
\]

Fix a dyadic number \(R>1\) and decompose
\[
\psi^n-\psi
=
P_{\leq R}(\psi^n-\psi)
+
P_{>R}(\psi^n-\psi).
\]
For the low-frequency part,
\begin{align}
\|P_{\leq R}(\psi^n-\psi)\|_{L_t^\infty\dot B^s_{2,1}([0,T])}
&=
\sum_{\rho\leq R}
\rho^s
\|P_\rho(\psi^n-\psi)\|_{L_t^\infty L_x^2([0,T])}
\nonumber\\
&\leq
R
\sum_{\rho\leq R}
\rho^{s-1}
\|P_\rho(\psi^n-\psi)\|_{L_t^\infty L_x^2([0,T])}
\nonumber\\
&\leq
R
\|\psi^n-\psi\|_{L_t^\infty\dot B^{s-1}_{2,1}([0,T])}.
\label{eq:low-frequency-critical-upgrade}
\end{align}
For fixed \(R\), the last expression tends to zero as \(n\to\infty\).

For the high-frequency part,
\begin{align}
\|P_{>R}(\psi^n-\psi)\|_{L_t^\infty\dot B^s_{2,1}([0,T])}
&\leq
\sum_{\rho>R}
\rho^s
\left(
\|P_\rho\psi^n\|_{L_t^\infty L_x^2([0,T])}
+
\|P_\rho\psi\|_{L_t^\infty L_x^2([0,T])}
\right)
\nonumber\\
&\lesssim
\sum_{\rho>R}c_\rho.
\label{eq:high-frequency-critical-tail}
\end{align}
Since \(\{c_\rho\}\in\ell^1\), the right-hand side tends to zero as
\(R\to\infty\), uniformly in \(n\). We first choose \(R\) sufficiently large
and then let \(n\to\infty\). It follows that
\[
\psi^n
\longrightarrow
\psi
\quad
\text{in }C([0,T];\dot B^s_{2,1}).
\]

The critical elliptic stability estimates similarly imply
\[
\lambda^n\longrightarrow\lambda,
\quad
\partial_x h^n\longrightarrow\partial_x h,
\quad
V^n\longrightarrow V,
\quad
A^n\longrightarrow A
\]
in the corresponding critical Besov topologies. The variable \(B^n\)
converges in the dyadic fixed-time and space--time topology associated with
the estimates for \(B\).

It remains to show that the limit satisfies the full system. We write the
Schr\"odinger equation in the form
\[
\sqrt{-1}\,\partial_t\psi^n
+
\Delta\psi^n
=
-\partial_\alpha
\left(
h_n^{\alpha\beta}\partial_\beta\psi^n
\right)
+
\mathcal N^n.
\]
For the principal quasilinear term,
\begin{align}
h_n^{\alpha\beta}\partial_\beta\psi^n
-h^{\alpha\beta}\partial_\beta\psi
=&\,
\left(
h_n^{\alpha\beta}-h^{\alpha\beta}
\right)
\partial_\beta\psi^n
\nonumber\\
&+
h^{\alpha\beta}
\partial_\beta(\psi^n-\psi).
\label{eq:principal-limit-decomposition}
\end{align}
The mixed coefficient estimates used in the proof of
Proposition~\ref{prop1} show that the right-hand
side of \eqref{eq:principal-limit-decomposition} tends to zero in the
corresponding nonlinear distribution space.

For the lower-order terms, we use
\[
u_{\mathfrak T}^n\partial_x\psi^n
-u_{\mathfrak T}\partial_x\psi
=
(u_{\mathfrak T}^n-u_{\mathfrak T})\partial_x\psi^n
+u_{\mathfrak T}\partial_x(\psi^n-\psi),
\]
\begin{align}
\psi^n u_i^n u_j^n
-\psi u_i u_j
=&\,
(\psi^n-\psi)u_i^n u_j^n
+
\psi(u_i^n-u_i)u_j^n
\nonumber\\
&+
\psi u_i(u_j^n-u_j),
\end{align}
and
\[
B^n\psi^n-B\psi
=
(B^n-B)\psi^n
+B(\psi^n-\psi).
\]
The nonlinear estimates from the preceding sections, together with the
strong critical convergence and the uniform space--time bounds, imply that
all these differences tend to zero in distributions on \([0,T]\times
\mathbb R^d\).

The elliptic equations pass to the limit by the same argument. Therefore
\(
\left(
\psi,
\lambda,
h,
V,
A,
B
\right)
\)
satisfies
\eqref{gaugsmcf}-\eqref{eq:B-elliptic}
in the distributional sense. Since
\(
\psi^n(0)=\psi_0^n
\longrightarrow
\psi_0
~~
\text{in }\dot B^s_{2,1},
\)
and
\(
\psi^n
\longrightarrow
\psi
~~
\text{in }C([0,T];\dot B^s_{2,1}),
\)
we obtain
\(
\psi(0)=\psi_0.
\)

Since \(T>0\) was arbitrary and all a priori estimates are uniform in
\(T\), the same construction gives a global solution satisfying
\(
\psi
\in
C([0,\infty);\dot B^s_{2,1}),
\)
and
\[
\sup_{t\geq0}
\sum_{i=1}^5
\|u_i(t)\|_{\dot B^s_{2,1}}
\lesssim
\varepsilon,
\]
together with the global space--time, bilinear, interaction Morawetz, and
\(B\)-estimates.

Now we prove uniqueness. Let
\[
\left(
\psi^{(a)},
\lambda^{(a)},
h^{(a)},
V^{(a)},
A^{(a)},
B^{(a)}
\right),
\qquad a=1,2,
\]
be two solutions in the critical solution class with the same initial datum:
\[
\psi^{(1)}(0)
=
\psi^{(2)}(0)
=
\psi_0.
\]
The weak difference estimate was proved first for smooth solutions. By
applying it to smooth approximations of the two solutions and passing to the
limit, it remains valid for critical solutions. Hence, on every finite
interval \(I=[0,T]\),
\[
D(I)
\lesssim
\|\psi^{(1)}(0)-\psi^{(2)}(0)\|_{\dot B^{s-1}_{2,1}}.
\]
The right-hand side is zero. Therefore
\(
D(I)=0.
\)
In particular,
\(
\psi^{(1)}
=
\psi^{(2)}
~~
\text{on }[0,T].
\)
The elliptic difference estimates then imply
\[
\lambda^{(1)}=\lambda^{(2)},
~~
h^{(1)}=h^{(2)},
~~
V^{(1)}=V^{(2)},
~~
A^{(1)}=A^{(2)},
~~
B^{(1)}=B^{(2)}.
\]
Since \(T>0\) is arbitrary, the two global solutions coincide on
\([0,\infty)\). This proves uniqueness. The continuous dependence on the initial data can be proved in the same manner. 

This completes the proof of Theorem \ref{thm1}.
 
\subsection{Proof of Corollary \ref{cor1}}

We divide the proof into three steps: the construction of the initial harmonic/Coulomb gauge, the reconstruction of the geometric flow from the good gauge variables, and the application of Theorem~\ref{thm1}. We begin with the construction of harmonic coordinates at the initial time.
\begin{proposition}[Critical Besov harmonic coordinates]
\label{prop:critical-besov-harmonic-coordinates}
Let $d\ge4$, $s=d/2-1$, and let $F:\mathbb{R}_x^d\to\mathbb{R}^{d+2}$ be an immersion with induced metric $g=\mathrm{Id}+h$. Assume
\[
 \|\nabla h\|_{\dot {B}^s_{2,1}}\le \varepsilon,
~~ h(x)\to0\quad (|x|\to\infty),
\]
where $\varepsilon>0$ is sufficiently small. Then there exists a unique change of coordinates
\[
 y=x+\varphi(x),~~ \varphi(x)\to0\quad (|x|\to\infty),
\]
with $\|\nabla\varphi\|_{L^\infty}\ll1$, such that the transformed metric $\widetilde g$ satisfies
$\widetilde g^{\alpha\beta}\widetilde\Gamma^\gamma_{\alpha\beta}=0$. Moreover,
\[
 \|\nabla^2\varphi\|_{\dot {B}^s_{2,1}}
 \lesssim \|\nabla h\|_{\dot {B}^s_{2,1}},
 ~~
 \|\nabla_y(\widetilde g-\mathrm{Id})\|_{\dot {B}^s_{2,1}}
 \lesssim \|\nabla_x(g-\mathrm{Id})\|_{\dot {B}^s_{2,1}}.
\]
If, in addition, $\|H\|_{\dot {B}^s_{2,1}}\le\varepsilon$, then, for $\widetilde F(y)=F(x(y))$,
\[
 \|\nabla_y^2\widetilde F\|_{\dot {B}^s_{2,1}}
 \lesssim \|\nabla h\|_{\dot {B}^s_{2,1}}+\|H\|_{\dot {B}^s_{2,1}}.
\]
\end{proposition}

\begin{proof}
We follow the harmonic coordinate construction of Huang--Tataru \cite[Proposition~8.2]{HT22}, replacing the Sobolev estimates there by the critical Besov embedding
$\dot {B}^{d/2}_{2,1}\hookrightarrow L^\infty$ and the corresponding paraproduct estimates. Writing $y^\gamma=x^\gamma+\varphi^\gamma$ and imposing $\Delta_g y^\gamma=0$ gives schematically
\[
 \Delta\varphi=\partial h+Q(h,\partial h)+Q(h,\partial^2\varphi)
      +Q(\partial h,\partial\varphi).
\]
Hence Riesz-transform bounds and critical Besov product estimates yield
\[
 \|\nabla^2\varphi\|_{\dot {B}^s_{2,1}}
 \lesssim \|\nabla h\|_{\dot {B}^s_{2,1}}
 +\varepsilon\|\nabla^2\varphi\|_{\dot {B}^s_{2,1}}.
\]
The last term is absorbed for small $\varepsilon$, and a contraction argument gives existence and uniqueness. The resulting map is bi-Lipschitz because
$\|\nabla\varphi\|_{L^\infty}\lesssim\|\nabla^2\varphi\|_{\dot {B}^s_{2,1}}$; standard paracomposition estimates then give the bound for $\widetilde g$.

The key geometric structure is that in harmonic coordinates the metric equation is genuinely elliptic,
\[
 \widetilde g^{\alpha\beta}\partial_{\alpha\beta}\widetilde h
   =Q(\partial\widetilde h,\partial\widetilde h)+Q(\widetilde\lambda,\widetilde\lambda),
\]
where the Gauss equation has been used and no second derivative occurs on the right-hand side. Likewise,
$\widetilde g^{\alpha\beta}\partial_{\alpha\beta}\widetilde F=H$.
Perturbative elliptic estimates therefore give the asserted bound for $\nabla_y^2\widetilde F$.
\end{proof}

Having fixed harmonic coordinates, we next impose the Coulomb condition on
the normal frame. As in Huang--Tataru~\cite{HT22}, we rotate an
asymptotically constant normal frame according to $m'=e^{i\theta}m$ and solve
$\Delta_g\theta=\nabla^\alpha A_\alpha$, with $\theta(\infty)=0$. The same critical Besov elliptic estimates give
$\|\nabla\theta\|_{\dot {B}^s_{2,1}}\lesssim\|A\|_{\dot {B}^s_{2,1}}$.
Thus the harmonic coordinate and Coulombframe constructions place the
initial immersion in the good gauge. In particular, the resulting initial
complex mean curvature satisfies
\[
\|\psi_0\|_{\dot B^s_{2,1}}
\lesssim
\|\nabla(g_0-\mathrm{Id})\|_{\dot B^s_{2,1}}
+
\|H_0\|_{\dot B^s_{2,1}}
\lesssim
\varepsilon.
\]

We now explain how to reconstruct the geometric flow from the global
good gauge solution. 
Let $(\psi,\lambda,h,V,A,B)$ be the global good gauge solution. The reconstruction is algebraically identical to Section~8 of Huang--Tataru~\cite{HT22}. One solves the spatial and temporal frame equations
\[
 \partial_\alpha F_\beta
 =\Gamma^\gamma_{\alpha\beta}F_\gamma
   +\operatorname{Re}(\lambda_{\alpha\beta}\bar m),
 \qquad
 \partial_\alpha^A m=-\lambda_\alpha^{\gamma}F_\gamma,
\]
\[
 \partial_tF_\alpha
 =-\operatorname{Im}\!\left((\partial_\alpha^A\psi-i\lambda_{\alpha\gamma}V^\gamma)\bar m\right)
 +\left(\operatorname{Im}(\psi\bar\lambda_\alpha^{\gamma})+\nabla_\alpha V^\gamma\right)F_\gamma,
\]
with the analogous equation for $\partial_t^B m$.
 
 The Gauss--Codazzi--Ricci identities provide the compatibility conditions
for the spatial frame equations, while the evolution equations for
\(g\), \(\lambda\), and \(A\) provide the space--time compatibility
conditions. Frobenius' theorem therefore yields the moving frame
\(
(F_1,\ldots,F_d,m).
\)
The immersion is then recovered by integrating
\(
\partial_\alpha F=F_\alpha.
\)
By construction,
\[
\partial_tF
=
-\operatorname{Im}(\psi\overline m)
+
V^\gamma F_\gamma,
\]
and hence
\[
(\partial_tF)^\perp
=
-\operatorname{Im}(\psi\overline m)
=
JH.
\]
Therefore \(F\) solves the skew mean curvature flow \eqref{smcf1}.

The geometric reconstruction is the same as in
Huang--Tataru~\cite{HT22}; the only difference is analytic. Here the
good gauge solution is global and is constructed at the critical Besov
regularity. To justify the reconstruction at this regularity, we first
apply the preceding argument to the smooth approximate solutions. The
critical a priori estimates and the weak difference estimates allow us to
pass to the limit in the frame equations, the constraint identities, and
the evolution equations. Consequently, the limiting immersion satisfies
the same geometric identities and solves \eqref{smcf1} globally.

We can now complete the proof. By
Proposition~\ref{prop:critical-besov-harmonic-coordinates} and the Coulomb
frame construction, the initial immersion can be placed in the good gauge,
with
\[
\|\psi_0\|_{\dot B^s_{2,1}}
\lesssim
\varepsilon.
\]
Theorem~\ref{thm1} then yields a unique global solution of the
Schr\"odinger--elliptic system satisfying all constraint and evolution
identities. The reconstruction argument above produces a corresponding
global solution of the skew mean curvature flow. Uniqueness and continuous
dependence follow from the corresponding properties of the good gauge
system once the normalization at infinity, and hence the residual
rigid-motion freedom, has been fixed. This completes the proof of
Corollary~\ref{cor1}.

\appendix
\section{Interaction Morawetz Estimates}  \label{appMora}
In this appendix, we establish interaction Morawetz estimates for the following quasilinear Schrödinger equation:
\begin{equation} \tag{QNLS}	 \label{qNLS}
	\begin{cases}
		\sqrt{-1} u_t+\sum_{\alpha,\beta=1}^d \partial_{\alpha }(g^{\alpha \beta} \partial_{\beta} u)=\mathcal{N}(u),\qquad x\in\mathbb{R}^d\  (d\ge 3),\\
		u|_{t=0}=u_0,
	\end{cases}
\end{equation}
where   the metric $ g $ satisfies the evolution equation  
\begin{equation}\label{dt_append}
  \pa_t g^{\alpha \beta} = H^{\alpha \beta}
\end{equation}
and we denote $ g^{\alpha \beta}=\delta_{\alpha \beta}+h^{\alpha \beta}$, with $h^{\alpha \beta}$ satisfying
\[
 \|h^{\alpha \beta}\|_{L^\infty} \ll 1.
\]

We begin by    introducing some notations.  Let  $\pa^\alpha := g^{\alpha \beta}\pa_\beta$ and
$$ x=(y,z)=(y_1,y_2,y_3,z_4,\dots,  z_{d}) \in \mathbb{R}^d, y \in \mathbb{R}^3, z \in \mathbb{R}^{d-3}. $$
For simplicity, we limit the Roman indices $i,j,k,l$ to range from $1$ to $3$,  corresponding to the variables \(y = (y_1, y_2, y_3), y^\prime = (y_1^\prime, y_2^\prime, y_3^\prime)\).

The mass balance law for    \eqref{qNLS}  is
\begin{align} 
	&\,\pa_t   \left| u \right|^2  +\partial_{\alpha}\mathrm{Im} ( \overline{u}    \partial^{\alpha}u  )=\mathrm{Im} \left( \overline{u}\mathcal{N}(u) \right)  ,  \label{mass_qNLS}
\end{align}
while the momentum balance law reads
\begin{equation}
  \begin{aligned}  
	&\partial_t \mathrm{Im} ( \overline{u} \partial_\gamma u ) + \partial_{\alpha}\mathrm{Re} (2\partial_\gamma u    \partial^{\alpha} \overline{u} - \frac{1}{2}  \pa_\gamma \pa^\alpha|u|^2) \\
	&=- \partial_\gamma\mathrm{Re} ( \bar{u} \mathcal{N}(u)  )+2\mathrm{Re}  (\mathcal{N}(u)\partial_\gamma \overline{u} )-\mathrm{Re} (\partial_{\alpha} u  \partial_\gamma h^{\alpha \beta} \partial_{\beta} \overline{u}).
\end{aligned} \label{mom_qNLS_pre}
\end{equation}
Using \eqref{dt_append}, we raise the index $\gamma$ in \eqref{mom_qNLS_pre} and obtain
\begin{equation}
  \begin{aligned}  
	&\partial_t \mathrm{Im} (  \overline{u} \partial^\gamma u ) + \partial_{\alpha}\mathrm{Re} (2\partial^\gamma u    \partial^{\alpha} \overline{u} - \frac{1}{2}  \pa^\gamma \pa^\alpha|u|^2) \\
	&=     \mathrm{Re}  (\mathcal{N}(u)\partial^\gamma \overline{u} - \bar{u} \partial^\gamma\mathcal{N}(u) ) + H \cdot u \cdot \pa_x u  + \pa_x h \cdot[ (\pa_x u)^2+u\cdot \pa_x^2 u]\\
    & = \mathcal{M}^\gamma(u) + H \cdot u \cdot \pa_x u  + \pa_x h \cdot[ (\pa_x u)^2+u\cdot \pa_x^2 u].
\end{aligned} \label{mom_qNLS}
\end{equation}
Let $a(y) = |y|$. Then there holds
\begin{align*}
	&a_{k}(y) := (\pa_k a)(y) = \frac{y_k}{|y|},\\
	&a_{kl}(y) := (\pa_{kl} a)(y) = \frac{\delta_{kl}}{|y|}-\frac{y_k y_l}{|y|^3},\\
	&(\Delta a)(y) = \frac{2}{|y|},\\
	& (\Delta^2 a)(y) = -8\pi \delta(y).
\end{align*}  
Moreover, for $m\ge 1$, we can interpret $   \pa^m_y a  $
  as   kernels for appropriate multipliers, which act like the composition of  $|\nabla_y |^{4-m}$ and Riesz type operators $R_y$. 

We define the \textit{spatially localized Morawetz interaction potential}
\begin{align} \label{interapoten}
	-\int_{\mathbb{R}^3\times \mathbb{R}^3\times  \mathbb{R}^{d-3} \times \mathbb{R}^{d-3}} a_k(y-y')   |u(y,z)|^2      \mathrm{Im} (\bar{u} \pa^k u) (y',z') \, \mathrm{d}y\mathrm{d}y'\mathrm{d}z\mathrm{d}z'
\end{align}
For convenience, we omit both the integration domain and the explicit volume elements in what follows. Unless otherwise stated, we set
$$ (y,z),  (y',z') \in \mathbb{R}^3\times \mathbb{R}^{d-3}, $$
and it is understood that the variables $y, y', z, z'$ implicitly include their corresponding measures $\mathrm{d}y,\mathrm{d}y',\mathrm{d}z,\mathrm{d}z'$, which are suppressed for brevity.

Differentiating the   potential \eqref{interapoten} in time and integrating by parts, we apply \eqref{mass_qNLS} and \eqref{mom_qNLS} to obtain
\begin{align}
	-\frac{ \mathrm{d}}{\mathrm{d}t} \int a_k(y-y')  |u(y,z)|^2      \mathrm{Im} (\bar{u} \pa^k u) (y',z') \, = \mathrm{J}(u) + \text{Err}_1+ \text{Err}_2,  \label{interapoten_t1}
\end{align}
where
\begin{align*}
   \mathrm{J}(u):=&\,\int a_{ kl}(y-y')\,  |u(y,z)|^2 \, [2\mathrm{Re} (\pa^k \overline{u}  \pa^l u) -   \frac{1}{2}\pa^k \pa^l|u|^2 ](y',z')\nonumber \\
	 &- \int a_{ kl}(y-y')\,    2 \mathrm{Im}(\overline{u}  \pa^l u)(y,z)  \mathrm{Im}(\overline{u}  \pa^k u)  (y',z')]\,  \nonumber \\
   \text{Err}_1:=&\,- \int a_k(y-y')\, \Big[ |u(y,z)|^2 \, \operatorname{Re}\big(\mathcal{N}(u)\partial^k \overline{u} - \bar{u} \partial^k\mathcal{N}(u)\big)(y',z') \nonumber\\
    &\qquad\qquad\qquad\qquad+ \mathrm{Im} \left( \overline{u}\mathcal{N}(u) \right)(y,z) \, \mathrm{Im}(\overline{u}  \pa^k u)(y',z') \Big]\,  \nonumber\\
	\text{Err}_2:=&\,\int a_k(y-y')\, |u(y,z)|^2 \, \{H \cdot u \cdot \pa_x u  + \pa_x h \cdot[ (\pa_x u)^2+u\cdot \pa_x^2 u]\}(y',z').\,   
\end{align*}

We now simplify \(\mathrm{J}(u)\). Integration by parts yields
\begin{align} \label{interapoten_c1}
	&-\int a_{kl}(y-y')\, |u(y,z)|^2 \,     \frac{1}{2}\pa^k \pa^l|u|^2 (y',z') \nonumber\\
	 =&\,\int a_{kl}(y-y')\,   \frac{1}{2}\pa^k|u(y,z)|^2\, \cdot \pa^l|u(y',z')|^2+\text{Err}_3
\end{align}
where
\begin{align*}
   \text{Err}_3:=  \int (\pa^2 a)(y-y')\big[ (\pa_x h  \cdot \pa_x|u|^2 )(y,z)\, |u|^2 (y',z')  + (\pa_x h \cdot |u|^2 )(y,z)\, \pa_x|u(y',z')|^2 \big].\nonumber
\end{align*}
 %Therefore
 %\begin{align*}
%	J(u)= &  2\int a_{kl}(y-y')\, \big[|u(y,z)|^2 \, \mathrm{Re} (\pa^k \overline{u}  \pa^l u  )(y',z')\nonumber\\
 %&\qquad\qquad\qquad\qquad\quad -    \mathrm{Im}(\overline{u}  \pa^l \lambda_k)(y,z)  \mathrm{Im}(\overline{u}_k \pa^k u)  (y',z')\big]\, \nonumber\\
 %  &\quad + \frac{1}{2}\int a_{kl}(x-y)\pa^k|u(y,z)|^2\,\pa^l|u(y',z')|^2  + \text{Err}_3\nonumber. 
 %\end{align*}
Observing the identity
 \begin{align}
	&\mathrm{Re}[( u(y,z) \pa^k u(y',z')-\pa^k u(y,z) u(y',z'))\cdot ( \bar{u}(y,z) \pa^l \bar{u}(y',z')-\pa^l \bar{u}(y,z) \bar{u}(y',z'))]\nonumber\\
	& = |u(y,z)|^2\mathrm{Re}(\pa^k \bar{u} \pa^l u)(y',z') +|u(y',z')|^2\mathrm{Re}(\pa^k \bar{u} \pa^l u)(y,z)\nonumber\\
	&\quad -2  \mathrm{Re}(\pa^k u(y,z) u(y',z') \cdot \bar{u}(y,z) \pa^l \bar{u}(y',z') )\nonumber\\
	& =|u(y,z)|^2\mathrm{Re}(\pa^k \bar{u} \pa^l u)(y',z') +|u(y',z')|^2\mathrm{Re}(\pa^k \bar{u} \pa^l u)(y,z) \nonumber\\
	&\quad-2\,\mathrm{Im}(\bar{u} \pa^k u)(y,z)\, \mathrm{Im}(\bar{u} \pa^l u)(y',z') - \frac{1}{2}\pa^k|u(y,z)|^2 \pa^l|u(y',z')|^2 \nonumber,
 \end{align}
together with the symmetry of the kernel $a_{kl}$ and the convexity of $a$, we use \eqref{interapoten_c1} to obtain
\begin{align} \label{interapoten_t2}
	  \mathrm{J}(u) &= \int a_{kl }(y-y')[F^k \bar{F}^l + \pa^l|u(y,z)|^2 \pa^l|u(y',z')|^2 ] + \text{Err}_3\\
	&\ge  \int a_{kl }(y-y') \pa^k|u(x)|^2 \pa^l|u(y)|^2  + \text{Err}_3,\nonumber
\end{align}
where $F^k = u(y,z) \pa^k u(y',z')-\pa^k u(y,z) u(y',z') $.
The dominant term is handled by further integration by parts:
\begin{align}
	 &\int a_{kl }(y-y') \pa^l|u(y,z)|^2 \pa^k|u(y',z')|^2 
	  =  8\pi \|u\|_{L^4_y L^2_z}^4 + \text{Err}_4, \label{interapoten_t3}
\end{align}
with
\begin{align*}
   \text{Err}_4:=& \,\int (\pa^2 a)(y-y')( \pa_x h \cdot |u|^2 )(y,z)\, \pa_x|u(y',z')|^2  \nonumber\\
	 & + \int (\pa^3 a)(y-y')( \pa_x h \cdot |u|^2 )(y,z)\,  |u(y',z')|^2 \nonumber\\
	 & + \int (\pa^4 a)(y-y')(   h \cdot |u|^2 )(y,z)\,  |u(y',z')|^2.
\end{align*}

Substituting \eqref{interapoten_t2} and \eqref{interapoten_t3} into \eqref{interapoten_t1} and integrating over \([0,T]\), we finally obtain the following theorem.
\begin{theorem}  \label{thm_im}
	If  $ u $ is the solution of \eqref{qNLS} and the metric $g$ satisfies the smallness condition $\|g-I\|_{L^\infty_x}\ll 1$ and the evolution equation 
	\begin{equation}
	  \pa_t g^{\alpha \beta} = H^{\alpha \beta},
	\end{equation}
	 then the following estimate holds:
	\begin{align}
	\|u\|_{L^4_{t,y}L^2_z}^4 \lesssim \|u\|_{L^\infty_t L^2_x}^3  \|\pa_x u\|_{L^\infty_t L^2_x} + \sum_{j = 1}^4 \|\mathrm{Err}_i\|_{L^1_t},
	\end{align}
	where
    \begin{align}
		\|\mathrm{Err}_{1,2}\|_{L^1_t}&\lesssim \|u\|_{L^\infty_t L^2_x}   \|\pa_x u\|_{L^\infty_t L^2_x} \|\operatorname{Im}(\bar{u}\mathcal{N}(u))\|_{L^1_t L^1_x}\\
		&\quad + \|u\|_{L^\infty_t L^2_x}^2(\|   H \cdot u \cdot \pa_x u,\  \pa_x h \cdot( (\pa_x u)^2+ u\cdot \pa_x^2 u)\|_{L^1_{t,x}}),\nonumber\\
        &\quad+ \Big| \int_0^T\int a_k(y-y')\,   |u(y,z)|^2 \, \operatorname{Re}\big(\mathcal{N}(u)\partial^k \overline{u} - \bar{u} \partial^k\mathcal{N}(u)\big)(y',z')  \mathrm{d} t\Big|\nonumber\\
		\|\mathrm{Err}_{3,4}\|_{L^1_t}&\lesssim   \Big|\int(|\nabla_y|^{-2}\pa_y  |u(y,z)|^2  \cdot (\pa_x h \cdot |u|^2)(y,z')\, \Big|\nonumber\\
		&\quad    +\Big|\int   |\nabla_y|^{-2} |u(y,z)|^2 \cdot (\pa_x h  \cdot \pa_x|u|^2)(y,z')\, \Big|   \nonumber\\
	&\quad    + \Big|\int(|\nabla_y|^{-1}|u(y,z)|^2 \cdot( \pa_x h  \cdot |u|^2)(y,z')\Big|\nonumber\\
	&\quad  + \Big|  \int R^4_y |u(y,z)|^2 \cdot (h  \cdot|u|^2)(y,z'))\Big| \nonumber.
\end{align}
\end{theorem}

\makeatletter
{\renewcommand{\addcontentsline}[3]{}
\section*{Acknowledgements}} % 不在目录中
\makeatother
  Y. Zhou was partially supported by NSFC
(Nos.\,12171097, 12571231). N. Lai  was partially supported by NSFC
(Nos.\,12271487, W2521007),	
J. Shao  was partially supported by NSFC
(No.\,12501315), Hunan Provincial Natural Science Foundation of China (No.\,2025JJ60042) and Key Laboratory of Computing and Stochastic Mathematics (Ministry of Education), Hunan Normal University.\vspace{0.3cm}\\
 \textbf{Data Availability } Data sharing is not applicable to this article as no datasets were generated or analyzed
during the study.
\vspace{0.1cm}\\
\textbf{Conflict of interest } The authors declare that there is no Conflict of interest in this article.
\vspace{0.6cm}
\bibliographystyle{plain}

\end{document}